%% file: finaldis.tex
\documentclass{report}

\usepackage{amsfonts,amsmath,amssymb}
\usepackage[dvips]{graphics}

\newtheorem{theorem}{Theorem}[section]
          
\newtheorem{lemma}[theorem]{Lemma}

\newtheorem{corollary}[theorem]{Corollary}
\newtheorem{definition}[theorem]{Definition}

\newcommand{\rem}[1]{{\bf Remark:}}
\newenvironment{proof}{\noindent {\bf Proof: }}{\QED\medskip}
\def\QED{{\hspace*{\fill}{\vrule height .5ex width 1ex }\quad} 
    \vskip 0pt plus20pt}
\newcommand{\be}{\begin{equation}}
\newcommand{\ee}{\end{equation}}
\newcommand{\bea}{\begin{eqnarray}}
\newcommand{\eea}{\end{eqnarray}}
\newcommand{\beann}{\begin{eqnarray*}}
\newcommand{\eeann}{\end{eqnarray*}}
\newcommand{\nn}{\nonumber}
\newcommand{\V}{\Vert}
\newcommand{\vt}{\vert}

\begin{document}
\pagenumbering{roman}
\begin{center}
{\Large \bf Local Solutions of the Dynamic Programming Equations and the Hamilton Jacobi Bellman PDE}\\
By\\
\large{\textbf{Carmeliza Luna Navasca}}\\[0.5cm]
\large{B.A. University of California, Berkeley 1996}\\[0.5cm]
\large{DISSERTATION}\\[0.5cm]
\large{Submitted in partial satisfaction of the requirements for the
degree of}\\[0.5cm]
\large{DOCTOR OF PHILOSOPHY}\\
\large{in}\\
\large{APPLIED MATHEMATICS}\\
\large{in the}\\
\large{OFFICE OF GRADUATE STUDIES}\\
\large{of the}\\
\large{UNIVERSITY OF CALIFORNIA, DAVIS}
\vskip 0.1in
\leftline{\large Approved:}
\underbar{\hskip 2.5in}
\vskip 0.1in
\underbar{\hskip 2.5in}
\vskip 0.1in
\underbar{\hskip 2.5in}
\vskip 0.1in
\large Committee in Charge
\vskip 0.1in
\large 2002
\end{center}
\large
\listoffigures
\newpage
\large
\tableofcontents


\newpage
\large
{\Large \bf ACKNOWLEDGEMENTS} \\
First and foremost, I would like to thank Professor Arthur J. Krener, without whose provisions of knowledge, experience, and encouragement this dissertation would not exist.  He introduced me to the rich field of control theory and to the world of mathematical research.  I am fortunate for having a major advisor who has given me his time and attention to provide the necessary tools to do mathematics.  In addition, he shared many of his creative ideas and techniques as well as his philosophy on mathematics.

I would also like to express my gratitude to the professors who served in my committees and who had been my teachers as well.  I had the pleasure of taking several optimization and stochastic dynamics courses from Professor Roger Wets.  Professor E. Gerry Puckett taught me some important numerical techniques in his Navier-Stokes courses.  I appreciate all the comments and suggestions made by Professor Roger Wets and Professor Jesus Deloera in the earlier version of this dissertation.

I would like to thank my parents, Carmelita and Jovino; they had made a decision to start their life over 17 years ago so that their children would have a better education.  This dissertation is one of the proof that they made the right the decision.  I'm also thankful for having a wonderful sister and two supportive brothers (Joanna, Joe and Joel) who are always there for me.  I would like to thank my friend Heath who made my stay at Davis memorable and fun.  I wish to thank the rest of my family and friends for the love and support:  the Smith family (Cora, Laura-Ashley --the next mathematician in the family, and Barry), Kris, Max and family, Manuel, Rose, Anthony, Bingki, the Navasca family in New Jersey, the Luna family in Hawaii, my relatives in the Philippines, the Starr family and all my other friends.

I met two important people at Davis and learned to important things.  One is my major advisor, Professor Krener, who taught me control theory.  The other is my fiance, Shannon Starr, who taught me love.  I'm grateful to have found a friend, a colleague, and a fiance all in one.  


\newpage
\begin{center}
{\Large \bf \underline{Abstract}}
\end{center}
\medskip
We present a method for the solution of the Dynamic Programming Equation (DPE) that arises in an infinite horizon optimal control problem.  The method extends to the discrete time version of Al'brecht's procedure for locally approximating the solution of the Hamilton Jacobi Bellman PDE.  Assuming that the dynamics and cost are $\mathcal{C}^{r-1}(\mathbb{R}^{n+m})$ and $\mathcal{C}^r(\mathbb{R}^{n+m})$ smooth, respectively, we explicitly find expansions of the optimal control and cost in a term by term fashion.  These finite expansions are the first $(r-1)$th and $r$th terms of the power series expansions of the optimal control and the optimal cost, respectively.

Once the formal solutions to the DPE are found, we then prove the existence of smooth solutions to the DPE that has the same Taylor series expansions as the formal solutions.  The Pontryagin Maximum Principle provides the nonlinear Hamiltonian dynamics with mixed direction of propagation and the conditions for a control to be optimal.   We calculate the forward Hamiltonian dynamics, the dynamics that propagates forward in time.  We learn the eigenstructure of the Hamiltonian matrix and symplectic properties which aid in finding the graph of the gradient of the optimal cost.  Furthermore, the Local Stable Manifold Theorem, the Stokes' Theorem, and the Implicit Function Theorem are some of the main tools used to show the optimal cost and the optimal control do exist and satisfy the DPE.

Assuming that there is a forward Hamiltonian dynamics is to assume that the Hamiltonian matrix is invertible.  If $0$ is eigenvalue of the Hamiltonian matrix, then $0$ is also a closed loop eigenvalue.  We consider $0$ as a possible closed loop eigenvalue since the optimal cost calculated term by term is true for all closed loop eigenvalues with magnitude less than 1.  We prove that there exists a local stable manifold for the bidirectional nonlinear Hamiltonian dynamics.  

We also present a method for numerically solving the Hamilton Jacobi Bellman (HJB) PDE that arises in the infinite horizon optimal control problem.  We compute smooth solutions of the optimal control and the optimal cost up to some degree $r-1$ and $r$, respectively.  The first step is to approximate around $0$ using Al'brecht's local solutions.  Then, using a Lyapunov criteria we find a point on a level set where we truncate the approximated solutions and begin another polynomial estimates. We generate the polynomial solutions in the same fashion described in the Cauchy-Kovalevskaya Theorem.  This process is repeated over at other points as the smooth solutions are patched together like a circular quilt. 


\newpage
\pagestyle{myheadings} 
\pagenumbering{arabic}
\markright{  \rm \normalsize CHAPTER 1. \hspace{0.5cm}
Introduction}
\large 
\chapter{Introduction}
\thispagestyle{myheadings}
\input{chapter1_3.tex}

\newpage
\pagestyle{myheadings}
\markright{  \rm \normalsize CHAPTER 2. \hspace{0.5cm} 
Formal Solution of the DPE}
\chapter{Formal Solution of the Dynamic Programming Equations in Discrete-Time} 
\thispagestyle{myheadings}
\input{chapter3_46.tex}

\newpage
\pagestyle{myheadings} 
\markright{  \rm \normalsize CHAPTER 3. \hspace{0.5cm} 
Existence of the Local Solution of DPE}
\chapter{Existence of the Local Solution of the DPE}
\thispagestyle{myheadings}
\input{chapter4_39.tex}

\newpage
\pagestyle{myheadings} 
\markright{  \rm \normalsize CHAPTER 4. \hspace{0.5cm} 
Local Stable Manifold for the Bidirectional Discrete-Time Dynamics}
\chapter{Local Stable Manifold Theorem for the Bidirectional Discrete-Time Dynamics}
\thispagestyle{myheadings}
\input{chapter5_5.tex}

\newpage
\pagestyle{myheadings} 
\markright{  \rm \normalsize CHAPTER 5. \hspace{0.5cm} 
Solution of the Hamilton Jacobi Bellman Equations}
\chapter{Solution of the Hamilton Jacobi Bellman Equations}
\thispagestyle{myheadings}

\input{chapter6_6.tex}


\newpage
\pagestyle{myheadings} 
\markright{  \rm \normalsize APPENDIX. \hspace{0.5cm}}
\appendix
\chapter{Codes}
\section{Driver 1}
{\it driver1.m}
{\small
\begin{verbatim}
% generates polynomials around 0 via albrecht
%
% Prager Example
%    Dynamics:        xdot=xu + u; x(0)=x0
%    Cost:            (ln x+1)^2 + u^2
%    Taylorized Cost: x^2 -x^3+(11/12)x^4-...+u^2
%-------------------------------------------------------------------
% Dimension Guide
% d=1,n=1,m=1: input f(1,2), l(1,3) (linear, starts with quadratic)
%              output ka=(1,1), py=(1,1) (linear,st w quadratic)        
% d=2,n=1,m=1: input f(1,5), l(1,7) (up to 2nd, up to 3rd)
%              output ka=(1,2), py=(1,2) (up to 2nd, up to 3rd)
% d=3,n=1,m=1: input f(1,9), l(1,12) (up to 3rd, up to 4th)
%              output ka=(1,3), py=(1,3) (up to 3rd, up to 4th)
% d=4,n=1,m=1: input f(1,14), l(1,18) (up to 4th, up to 5th)
%              output ka=(1,4), py=(1,4) (up to 4th, up to 5th)
%--------------------------------------------------------------------
%Subroutines:
%           poly1.m
%           endpoint.m
%           theG.m
%           theF.m
%Find coeffients of Taylor''s Expansion

%d -- degree up to
%xbar -- Taylor expand around xbar

function driver1(d)

mesh=256;

[x1,pie,U,dom1,dom2,rrealu,rrealpy]=poly1(d);
[root]=endpoint(pie);

subplot(2,1,1),plot(x1,pie,'m');
hold;
plot(x1,rrealpy,'c--');
legend('\pi^0(x)','\pi^{*}(x)');
hold;
title('\pi(x) of degree d+1=4');
xlabel('x');
ylabel('\pi');
axis([0 4 -2 20]);

subplot(2,1,2),plot(x1,U);
hold;
plot(x1,rrealu,'c--');
legend('\kappa^0(x)','\kappa^{*}(x)');
hold;
title('\kappa(x) of degree d=3');
xlabel('x');
ylabel('\kappa');
axis([0 4 -10 10]);
\end{verbatim}
}

{\bf Subroutines}

{\it poly1.m}

{\small
\begin{verbatim}
% make polynomial functions
function [x1,pie,U,dom1,dom2,rrealu,rrealpy]=poly1(d)

n=1;
if (d==1)
  f=[0 1];
  l=[1 0 1];
elseif (d==2)
  f=[0 1 0 1 0];
  l=[1 0 1 -1 0 0 0];
elseif (d==3)
  f=[0 1 0 1 0 zeros(1,4)];
  l=[1 0 1 -1 0 0 0 11/12  zeros(1,4) ];
else
  f=[0 1 0 1 0 zeros(1,9)];
  l=[ 1 0 1 -1 0 0 0 11/12 zeros(1,4) -5/6 zeros(1,5)];
  % l=[ 1 0 1 -1 0 0 0 -1 zeros(1,10) ];
end% ifloop

[ka,fk,py,lk]= hjb(f,l,1,1,d);

%generating vectors X, Y (py & ka)
xslots=nchoosek(n+(1+1)-1,1+1);
yslots=nchoosek(n+(1)-1,1);
X=zeros(xslots,1);
Y=zeros(yslots,1);

if (d>=2)
  for i=2:d
   xslots=nchoosek(n+(d+1)-1,d+1);
   yslots=nchoosek(n+(d)-1,d);
   X=[[X] zeros(xslots,1)];
   Y=[[Y] zeros(yslots,1)];
  end% dloop
end %ifloop
X=X';
Y=Y';
[xsize,dummy]=size(X);
[ysize,dummy]=size(Y);

aa=0;
bb=4;
mesh=256;
dx=(bb-aa)/mesh;
dy=dx;

for i=1:mesh
  x1(i)= aa + i*dx;
  a=x1(i);
  X=[a^2; a^3; a^4];
  pie(i)=py*X;
  PY=pie(i);
  realpy(i)=x1(i)^2 -x1(i)^3 +(11/12)*x1(i)^4;
  rrealpy(i)=(log(x1(i)+1))^2;
  P=py*X;
  Y=[a; a^2; a^3];
  U(i)=ka*Y;
  realu(i)=-(x1(i)-(1/2)*x1(i)^2+(1/3)*x1(i)^3);
  rrealu(i)=-log(1+x1(i)); 
  FK=pragerf(x1(i),U(i));
  LK=pragercost(x1(i),U(i));
  gradP=0;

for k=1:d
 gradP=gradP + (k+1)*py(k)*x1(i)^(k);
 gP(i)=gradP;
end %gradPloop

end %for loop
\end{verbatim}
}

{\it endpoint.m}
{\small
\begin{verbatim}
%finds pt on levelset
function [root]=endpoint(pie)

C=fmincon('theG',2,0,0,0,0,0,4,'theF');

%finding root
root=100;
m=100;
mesh=256;

for i=1:mesh
     if (abs(C-pie(i))==0)
     root=i;
     end
end %for loop
\end{verbatim}
}

{\it theG.m}
{\small
\begin{verbatim}
[ka,py]=justgivemepk;
G=-py*[x^2; x^3; x^4];
\end{verbatim}
}

{\it theF.m}
{\small
\begin{verbatim}
%contraint func for fmincon
function [F,eqF]=theF(x)

[ka,py]=justgivemepk;
u=ka*[x; x^2; x^3];
FK=pragerf(x,u);     
LK=pragercost(x,u);
gradP=0;
d=3;

for k=1:d
 gradP=gradP + (k+1)*py(k)*x^(k);
end %gradPloop

  eps=2^(-6);
  F1=gradP*FK+(1-eps)*LK;
  F2=-(gradP*FK+(1-eps)*LK);
  F=[F1;F2];
  eqF=0;
\end{verbatim}
}

{\it justgivepk.m}
{\small
\begin{verbatim}
%output coefficients
function [ka,py]=justgivemepk

d=3;
n=1;
if (d==1)
  f=[0 1];
  l=[1 0 1];
elseif (d==2)
  f=[0 1 0 1 0];
  l=[1 0 1 -1 0 0 0];
elseif (d==3)
  f=[0 1 0 1 0 zeros(1,4)];
  l=[1 0 1 -1 0 0 0 11/12  zeros(1,4) ];
else
%d=4
  f=[0 1 0 1 0 zeros(1,9)];
  l=[ 1 0 1 -1 0 0 0 11/12 zeros(1,4) -5/6 zeros(1,5)];
end% ifloop

%hjb(f,l,n,m,d,f_,n_,m_) 
[ka,fk,py,lk]= hjb(f,l,1,1,d);
\end{verbatim}
}

\section{driver 2}
{\it Driver2.m}
{\small
\begin{verbatim}
% generate polynomial around the x_0
% improves albrecht approx
%
% Prager Example
%    Dynamics:        xdot=xu + u; x(0)=x0
%    Cost:            (ln x+1)^2 + u^2
%    Taylorized Cost: x^2 -x^3+(7/12)x^4-...+u^2
%-------------------------------------------------------------------
% Dimension Guide
% d=1,n=1,m=1: input f(1,2), l(1,3) (linear, starts with quadratic)
%              output ka=(1,1), py=(1,1) (linear,st w quadratic)        
% d=2,n=1,m=1: input f(1,5), l(1,7) (up to 2nd, up to 3rd)
%              output ka=(1,2), py=(1,2) (up to 2nd, up to 3rd)
% d=3,n=1,m=1: input f(1,9), l(1,12) (up to 3rd, up to 4th)
%              output ka=(1,3), py=(1,3) (up to 3rd, up to 4th)
% d=4,n=1,m=1: input f(1,14), l(1,18) (up to 4th, up to 5th)
%              output ka=(1,4), py=(1,4) (up to 4th, up to 5th)
%--------------------------------------------------------------------
%Subroutines:
%           kovalesky.m
%           poly2.m
%           glue.m
%
%d -- degree up to
%xbar -- Taylor expand around xbar

function driver2(d)

mesh=256;
[x1,pie,U,gP]=poly1(d);

for j=1:4

root=endpoint(pie);
ii=root;
[ka,py]=kovalesky(d,xbar,k1,pi1,pi2);
[x1,pie2,U2,rrealu,rrealpy,gP]=poly2(d,ka,py,xbar,x1);
[newp,newu]=glue(pie2,U2,pie,U,istar,x1,rrealu,rrealpy);

figure;
subplot(2,1,1);
hold;
plot(x1,rrealpy,'c--');
plot(x1,newp,'k');
hold;
legend('\pi_{*}(x)','\pi_{new}(x)');
title('\pi(x) of degree d+1=4');
xlabel('x');
ylabel('\pi');
axis([0 4 -2 20]);

subplot(2,1,2);
hold;
plot(x1,rrealu,'c--');
plot(x1,newu,'k');
hold;
legend('\kappa_{*}(x)','\kappa_{new}(x)');
title('\kappa(x) of degree d=3');
xlabel('x');
ylabel('\kappa');
axis([0 4 -10 10]);

pie=newp;
U=newu;

end %for
\end{verbatim}
}

{\bf Subroutines}

{\it kovalesky.m}
{\small
\begin{verbatim}
%generates coefficients
function [ka,py]=kovalesky(d,xbar,k1,pi1,pi2)

%py -- vector containing Taylor coefficients of py
%ka -- vector containing Taylor coefficients of ka

ka=zeros(1,d+1);
py=zeros(1,d+2);
%degree=d;

     if (d >= 1)
        [ka1,py1]=Coeffd1(xbar,ka,py,d,k1,pi1,pi2);
     end% if loop

     if (d >= 2)
     [ka2,py2]=Coeffd2(xbar,ka1,py1);
     end% if loop

     if (d >= 3)
        [ka,py]=Coeffd3(xbar,ka2,py2);
     end% if loop
\end{verbatim}
}

{\it poly2.m}
{\small
\begin{verbatim}
%set up polynomial functions
function [x1,pie,U,rrealu,rrealpy,gP]=poly2(d,ka,py,xbar,x1)

mesh=256;

for i=1:mesh
  a=x1(i);
  X=[1; (a-xbar);(1/factorial(2))*(a-xbar)^2;(1/factorial(3))*(a-xbar)^3;(1/factorial(4))*(a-xbar)^4]; 
  pie(i)=py*X;
  PY=pie(i);
  rrealpy(i)=(log(x1(i)+1))^2;
  Y=[1; (a-xbar);(1/factorial(2))*(a-xbar)^2;(1/factorial(3))*(a-xbar)^3]; 
  U(i)=ka*Y;
  rrealu(i)=-log(1+x1(i)); 
  FK=pragerf(x1(i),U(i));
  LK=pragercost(x1(i),U(i));
  gradP=0;

for k=1:d+1
   gradP=gradP + (k)*(1/factorial(k))*py(k+1)*(x1(i)-xbar)^(k-1);
   gP(i)=gradP;
end %gradPloop

end %for loop i
\end{verbatim}
}

{\it glue.m}
{\small
\begin{verbatim}
% attach new polynomial in appropriate interval
function [newp,newu]=glue(p,u,pie,U,istar,x1,rrealu,rrealpy)

mesh=256;
[dum n]=size(p);
ii=istar;

for i=1:ii-1
  newp(i)=pie(i);
  newu(i)=U(i);
end

newp(ii:n)=p(ii:n);
newu(ii:n)=u(ii:n);

\end{verbatim}
}

{\it coeffd1.m}
{\small
\begin{verbatim}
%Calculate the deg=1 coefficients of Prager''s Example
function [k,p]=Coeffd1(a,k,p,deg,k1,pi1,pi2)

[PP,KK]=Tconst(a,deg);
p(1)=pi1;
k(1)=k1;
p(2)=pi2;
l=(log(a+1))^2 + k(1)^2;
f=a*k(1) + k(1);
p(3)=(-1/f)*( p(2)*( k(2) + k(1) +a*k(2) ) +2*log(a+1)/(a+1) +2*k(1)*k(2));
k(2)= (1/2)*(- p(2) - p(3) - a*p(3));
\end{verbatim}
}

{\it coeffd2.m}
{\small
\begin{verbatim}
%Calculate the deg=2 coefficients of Prager''s Example
function [k,p]=Coeffd2(a,k,p)

l=pragercost(a,k(1));
f=pragerf(a,k(1));
p(4)=(-1/f)*( 2*p(3)*(k(2)+ k(1) + a*k(2)) + p(2)*(k(3)+2*k(2)+a*k(3)) +2/(a+1)^2 -2*log(a+1)/(a+1)^2 + 2*k(2)^2 +2*k(1)*k(3));
k(3)=(-1/2)*(p(3) + p(4) +a*p(4));
\end{verbatim}
}

{\it coeff3.m}
{\small
\begin{verbatim}
%Calculate the deg=3 coefficients of Prager''s Example
function [k,p]=Coeffd3(a,k,p)

l=pragercost(a,k(1));
f=pragerf(a,k(1));
p(5)=(-1/f)*( 3*p(4)*(k(2) + k(1) + a*k(2)) + 3*p(3)*(k(3) + 2*k(2) + a*k(3)) 
+ p(2)*(k(4)+3*k(3)+a*k(4)) -6/(a+1)^3 + 4*log(a+1)/(a+1)^3 + 6*k(2)*k(3) +2*k(1)*k(4));
k(4)=(-1/2)*(3*p(4) + p(5) + a*p(5)); 
\end{verbatim}
}


\newpage
\pagestyle{myheadings} 
\markright{  \rm \normalsize BIBLIOGRAPHY. \hspace{0.5cm}}

\end{document}

%% file: chapter1_3.tex
We solve two equations: the Hamilton-Jacobi-Bellman (HJB) PDE and the Dynamic Programming Equations (DPE).  Both equations are associated with the infinite horizon optimal control problem of minimizing a running cost subject to a nonlinear control dynamics.  Since the optimal control problems arise ubiquitously in engineering, economics and biological science where models of these types surface naturally, this is the motivation to conjure up methods in solving them.  It is actually an economist who first studied optimal control problems in modelling capital accumulation \cite{Ra28}.  However, it is the advent of modern control theory that had considerable impact on the treatment of these problems.  The classic optimal control problem, the nonlinear regulator problem, is what we consider in this dissertation.  We seek a feedback or a control which minimizes the running cost under the nonlinear control system for the regulator problem.  The optimal control maintains the dynamics close to $0$ while keeping the expenditure of the cost at a minimum.  The Dynamic Programming technique is used to derive the HJB PDE and the DPE from the infinite horizon optimal control problem.  Thus, our main reason for solving the HJB PDE and DPE is to find a stabilizing feedback and a Lyapunov function to confirm local asymptotic stability of the given nonlinear control system.  The HJB PDE corresponds to the optimal control problem with continuous dynamics and cost while the DPE is the discrete-time analogue.  

Around 1961 one of the earliest results was made by Albrecht \cite{Al61} on the nonlinear regulator problem.  He discovered the formal power series solutions to the HJB equations while studying analytical nonlinear control systems.  In 1969, Lukes proves the existence and uniqueness of these power series representations of the optimal control feedback and the optimal cost.  To our knowledge there has been no extension of these results to the discrete-time case.  So we look at the discrete-time version of the optimal control problem with the discrete-time nonlinear dynamics and cost.  Using the method of Albrecht, we find the formal power series solutions to the DPE.  Then we prove the existence and uniqueness of these power series solutions in Theorem \ref{thebigkahuna}.  The Pontryagin Maximum Principle (PMP) \cite{AM90}, \cite{CC89} gives the necessary condition for a control to be optimal.  The PMP also provides the Hamiltonian dynamics satisfied by the optimal state and the costate trajectories.  The existence of the optimal cost is proved in two cases.  One case is when we assume that the Hamiltonian map is invertible.  The invertibility of the Hamiltonian matrix allows us to rewrite the dynamics where both the state and the costate dynamics propagate in the direction where time approaches infinity.  The other case is when the Hamiltonian map is not a diffeomorphism; i.e., we have a nonlinear bidirectional Hamiltonian discrete-time dynamics.  Finding the stable manifold for the Hamiltonian dynamics in both cases is key to finding the existence of the optimal cost.  For invertible maps, Hartmann \cite{Ha82} shows the existence of a stable manifold by the method of successive approximations on the implicit functional equation.  There is the technique of using the Contraction Mapping Theorem on a complete space developed by Kelley \cite{Ke66} in 1966.  Applications of this method can be found in the book by Carr \cite{Ca81} and the paper of Krener \cite{Kr01}.  There is another technique by Irwin \cite{Ir70} based on an application of the inverse function theorem on a Banach space of sequences.  One of our main results is Theorem \ref{nondiffthm}; we show that there exists a local stable manifold for the bidirectional discrete map with a hyperbolic fixed point.  After a two-step process of diagonalizing the dynamics we apply the technique of Kelly \cite{Ke66} on a complete space of Lipschitz functions with the supremum norm.

Computing the solutions of the HJB equations in higher dimension is major challenge.  A standard approach requires the temporal and spacial discretizations of the optimal control problem and then solves the corresponding nonlinear program, see  \cite{KD92} and the appendix by Falcone in \cite{BCD97}.  Other methods for solving the HJB PDE are similar to those for conservation laws \cite{OS91} and marching methods \cite{Se99}.  Although these standard approaches work for such equations for low dimensional systems, these numerical methods become infeasible for real world problems which are higher dimensional systems.  Our method is a higher order approach which requires few discretizations at each dimension.  We start with Albrecht's solutions as the initial approximations for the HJB PDE.  We then truncate the computed solutions at the point where the optimal cost satisfies the optimality and stability constraints and begin new approximations at the same point. The idea is to patch together successive approximations to obtain a larger domain for which the numerically computed optimal cost is still a Lyapunov function.  The numerical procedure is described in \cite{Na00}.

The outline of the dissertation is as follows:
In Chapter 2, we derive the DPE.  We first discuss the linear-quadratic regulator problem followed by the description of the method of finding formal power series solutions of the DPE.  The coefficients of the polynomials are expressed in linear equations and the discrete-time algebraic Riccati equation (DTARE).  We discuss the solvability of these equations as well.  In Chapter 3, we prove the existence of the optimal solutions of the DPE.  We show some properties of the Hamiltonian dynamics that are of importance in the proof of the existence of the optimal cost.  Chapter 4 is where we prove the existence of a stable manifold for the noninvertible Hamiltonian map.  It is actually the generalization of the results in Chapter 3.  In Chapter 5, we summarize our numerical approach for solving the HJB PDE.  We include numerical results from a 1-d example and describes the actual algorithm in more detail.

%% file: chapter3_46.tex
We present a method for the solution of the Dynamic Programming Equations, the discrete-time Hamilton Jacobi Bellman PDE that arises in an infinite time optimal control problem.  The method extends to the discrete time Al'brecht's procedure for locally approximating the solution.  Assuming that the dynamics and cost are $\mathcal{C}^{k-1}(\mathbb{R}^{n+m})$ and $\mathcal{C}^k(\mathbb{R}^{n+m})$ smooth, we explicitly find expansions of the optimal control and cost in a term by term fashion.  These finite expansions are the first few terms the power series expansions of the optimal control and cost. 

\section{Preliminaries}
First we discuss some definitions and notions about control systems.
Consider the $n$-dimensional input and output state equations
\beann
x^+&=& Ax + Bu\\
y&=& Cx + Du.
\eeann
where $A,~ B,~ C,$ and $D$ are  $n \times n,~n \times m,~l \times  n,~l \times  m$ matrices with $x \in \mathbb{R}^n$, $y \in \mathbb{R}^l$, and $u \in \mathbb{R}^m$.  We use the notation $x^+=x_{k+1}$ and $x=x_k$.

\begin{definition} \label{stab}
A system  is \emph{controllable} when any inital state $x_0$ can be driven to the final state $x_F$ in finite number of steps.  Equivalently, if the $n \times  nm$ controllability matrix
\[\mathcal{C}=[B, AB, \ldots, A^{n-1}B]\] 
has rank $n$, then $(A,B)$ is a controllable pair.  Otherwise (A,B) is said to be uncontrollable.  The pair (A,B) is stabilizable if all uncontrollable modes are asymptotically stable.
\end{definition}
If a state equation is controllable, then all eigenvalues can be assigned arbitrarily by introducing a state feedback.  Moreover, every uncontrollable system can be diagonalized into
\bea \label{uncontrol}
\left[ \begin{array}{c}
x_c^+ \\
x_u^+ \\
\end{array}\right]=
\left[ \begin{array}{cc}
A_c & \bar{A} \\
0   & A_u\\
\end{array}\right]
\left[ \begin{array}{c}
x_c \\
x_u \\
\end{array}\right]+
\left[ \begin{array}{c}
b_c \\
0   \\
\end{array}\right]
\eea
where $(A_c,b_c)$ is controllable.  If $A_u$ is stable and $(A_c,b_c)$ is controllable, then the system (\ref{uncontrol}) is stabilizable.  As in the definition above, we refer to the eigenvalues of $A_u$ as the uncontrollable modes.

Analogously, we define observability and dectectibility, the duals of controllabilily and stabilizability, respectively.
\begin{definition} \label{dect}
A sytem is observable if for any unknown $x_0$, $\exists~T$ s.t. the input $u_k$, and output $y_k$ uniquely determine $x_0$.  Equivalently, if the $nl \times n$ observability matrix
\[\mathcal{O}=[C, CA, \ldots, CA^{n-1}]^T\]
has rank $n$.  Otherwise, $(A,C)$ is said to be unobservable.
The pair $(A,C)$ is detectable if the unobservable eigenvalues are stable.
\end{definition}
 
Given a control dynamics
\beann
x^+&=&Ax+Bu\\
x(0)&=&x_0
\eeann
where the state $x \in \mathbb{R}^n$, the control $u \in \mathbb{R}^m$, $A$, $B$ are $n \times n$ and $n \times m$ matrices, respectively.  
We want to find $u=Kx$ such that the system is driven to 0; i.e., the system
\[x^+=(A+BK)x\]
is asymptotically stable around 0; i.e., the spectrum of $A+BK$ lies inside the unit circle.  We call such control a stabilizing feedback.  One way to solve this stabilization problem is to set up an optimal control problem.  This will be discuss in the next sections in details.

\section{Discrete-Time Optimal Control Problem}
We formulate a discrete in time infinite horizon optimal control problem of minimizing the cost functional, 
\[ \min_u \sum_{k=0}^{\infty} l(x_k,u_k)\]
subject to the dynamics
\beann
x^+&=&f(x,u)\\
x(0)&=&x_0
\eeann
where the state vector $x\in \mathbb{R}^n$, the control $u\in \mathbb{R}^m$, and 
\beann
f(x,u)&=&Ax + Bu + f^{[2]}(x,u) + f^{[3]}(x,u) + \ldots \\
l(x,u)&=&\frac{1}{2}x'Qx + x'Su + \frac{1}{2}u'Ru + l^{[3]}(x,u) + \ldots
\eeann
Here we denote $f^{[m]}(x,u)$ and $l^{[m]}(x,u)$ as homogeneous polynomials in $x$ and $u$ of degree $m$.

To solve the optimal control problem is to look for an optimal feedback $u^*=\kappa(x)$ such that the cost functional is kept at its minimum, namely the optimal cost $\pi(x_0)$ of starting the system at $x_0$, while driving the dynamics to $0$.

Here are the assumptions:
We assume $l(x,u)$ is convex in $x$ and $u$ so that
\beann
\left[ \begin{array}{cc}
Q   & S \\
S^* & R \\
\end{array}\right] \geq 0 
\eeann
and let $R >0.$
In addition, it is assumed that the pair $(A,B)$ is stabilizable and the pair $(A,Q^{1/2})$ is detectable.



\section{Derivation of the Dynamic Programming Equations}
Given that $x(0)=x_0$ the optimal value function $\pi(x)$ is defined by
\[\pi(x_0)=\min_u \sum_{k=0}^{\infty} l(x_k,u_k).\]
This value function satisfies a functional equation, called the dynamic programming equation.  The optimal feedback $\kappa(x)$ is constructed from the dynamic programming equation \cite{BCD97}.  First, we state the optimality principle:
\begin{theorem}
Discrete-Time Optimality Principle:  
\bea \label{optcond}
\pi(x)=\min_{u} \{\pi(f(x,u)) + l(x,u)\}
\eea
\end{theorem}
\begin{proof}
We have that
\beann
\pi(x_0)&=&\min_{u_0} \{\sum_{k=0}^{\infty} l(x_k,u_k)\}\\
        &=&\min_{u_0} \{l(x_0,u_0) +\sum_{k=1}^{\infty} l(x_k,u_k) \}\\
        &=&\min_{u_0} \{l(x_0,u_0) + \pi(x_1)\}
\eeann
Generalizing the optimality principle at the $k^{th}$-step, we have
\bea \label{dhjb1}
\pi(x)=\min_{u} \{\pi(x^+) +  l(x,u)\}.
\eea
\end{proof}
The optimality equation (\ref{optcond}) is the first equation of the dynamic programming equations.  An optimal policy $u^*=\kappa(x)$ must satisfy
\beann
\pi(x) -  \pi(f(x,u^*)) - l(x,u^*)  = 0
\eeann
if we assume convexity of the LHS of (\ref{dhjb1}).  We can find $u^*$ through
\[\frac{\partial (\pi(x) - \pi(f(x,u)) - l(x,u)}{\partial u}=0\]
which by the chain rule becomes 
\[\frac {\partial \pi}{\partial x}(f(x,u)) \frac{\partial f}{\partial u}(x,u) +  \frac{\partial l}{\partial u}(x,u)=0\]
Thus, $\pi(x)$ and $\kappa(x)$ satisfy these equations, the Dynamic Programming Equations (DPE):
\bea
\pi(x) -  \pi(f(x,u)) - l(x,u) ) &=& 0 \label{dpe1}\\
\frac {\partial \pi}{\partial x}(f(x,u)) \frac{\partial f}{\partial u}(x,u) + \frac{\partial l}{\partial u}(x,u) &=& 0 \label{dpe2}
\eea
We now introduce a method for solving the dynamic programming equations for $\kappa(x)$ and $\pi(x)$.



\section{Power Series Expansion}
Our method of solving the DPE is an extension of Al'brecht idea for continuous-time systems \cite{Al61}.  We require that $f(x,u) \in \mathcal{C}^{k-1}(\mathbb{R}^{n+m})$ and $l(x,u) \in \mathcal{C}^k(\mathbb{R}^{n+m})$ are expressible in Taylor's form around 0 in $N_{\varepsilon}(0)$.  Then, all the series representation of the given $f(x,u),\;l(x,u)$ and the unknown $\kappa(x), \pi(x)$ are substituted into the DPE.  With the exception of the first level, gathering the terms of the same degree in each equation will yield linear equations.  At the first level, we obtain the discrete-time algebraic Riccati equation.  The advantage of our process is the reduction of DPE which is a nonlinear system of equation into a system of one Riccati equation and many linear equations. We will see that the set of linear equations has a triangular structure within the levels of degree, a structure that allows the system to be easily solvable. 

On some neighborhood around 0, $N_{\varepsilon}(0) \subset \mathbb{R}^{n+m}$, the dynamics and cost assume the following form:
\beann
f(x,u)&=&Ax + Bu + f^{[2]}(x,u) \\
      &&+ f^{[3]}(x,u) + \ldots + f^{[d]}(x,u)
\eeann
\beann
l(x,u)&=&\frac{1}{2}x'Qx + x'Su + \frac{1}{2}u'Ru \\
      &&+ l^{[3]}(x,u) + \ldots + l^{[d+1]}(x,u)
\eeann
where $f \in \mathcal{C}^{k-1}$ and $l \in \mathcal{C}^k$. 

In consequence, we expect the unknowns to have power series expansions,
\bea 
\pi(x)&=& \frac{1}{2}x'Px + \pi^{[3]}(x) + \ldots \label{whatlook41}\\
\kappa(x)&=& Kx + \kappa^{[2]}(x) + \ldots \label{whatlook42}
\eea
as well.  
The known matrices are $A,\;B,\;f^{[2]},\;\ldots$ and $Q,\;S,\;R,\; l^{[3]},\; \ldots$ while $P,\;\pi^{[3]},\;\ldots$ and $K,\;\kappa^{[2]},\;\ldots$ are the unknowns.

\subsection{Special Case: Linear-Quadratic Regulator (LQR)}
The linear-quadratic regulator is an infinite-time horizon
optimal control problem with
\bea \label{lindynamics}
f(x,u)&=&Ax + Bu 
\eea
\bea \label{quadcost}
l(x,u)&=&\frac{1}{2}x'Qx + x'Su + \frac{1}{2}u'Ru; 
\eea
i.e., the higher degree homogeneous polynomials are zero.  Also, in a linear-quadratic regulator, we expect that the optimal cost and the optimal control will be quadratic and linear, respectively; i.e.,
\beann
\pi(x)&=& \frac{1}{2}x'Px \\
\kappa(x)&=& Kx.   
\eeann
Thus, we look for P and K.  After substituting all the expansions in (\ref{dpe1})
and collecting quadratic terms, we obtain
\bea \label{preric}
\frac{1}{2}x'\Big[P - A'PA - K'B'PA - A'PBK - K'B'PBK - Q - 2SK- K'RK \Big]x=0 \nonumber \\
\eea
Meanwhile gathering linear terms of (\ref{dpe2}), we get
\bea \label{prekappa1}
(Ax + Bu)'PB + x'S + u'R = 0
\eea
that reduces to
\bea \label{prekappa2}
K= -(B'PB + R)^{-1}(A'PB+ S)'.
\eea
It follows that (\ref{preric}) can be simplified to
\bea \label{finric}
\frac{1}{2}x'\Big[P - A'PA + (A'PB+ S)(B'PB+ R)^{-1}(A'PB+ S)' -Q \Big]x=0
\eea
by (\ref{prekappa2}).
Thus, (\ref{finric}) and (\ref{prekappa1}) are the pair of equations obtained by collecting the quadratic terms of (\ref{dpe1}) and the linear terms of (\ref{dpe2}):
\bea
0 &=& P - A'PA + (A'PB+ S)(B'PB+R)^{-1}(A'PB+S)' -Q \label{ric1}\\    
K &=& -(B'PB + R)^{-1}(A'PB+ S)' \label{kappa1}
\eea
Equation (\ref{ric1}) is known as the Discrete Algebraic Riccati Equation (DTARE).  Since $(B'PB + R)$ is positive definite, the matrix $K$ is well-defined once P is known.  The next theorem gives the necessary condition for the DTARE to have a unique positive definite solution P.
\begin{theorem} \label{specthm}
If the pair $(A,B)$ is stabilizable and the pair $(A,Q^{1/2})$ is detectable, then there exists a unique positive definite matrix P which satisfies the Riccati equation and 
\[\vert \sigma(A+BK) \vert < 1. \]
where $\sigma(A+BK)$ is spectrum of (A+BK).
\end{theorem}
Moreover, the resulting feedback $u=Kx$ is asymptotically stabilizing for the system
\[x^+=(A+BK)x.\]
The above theorem can be found in \cite{AM90}.

\subsection{First Level}
When the problem has nonlinear dynamics and cost with higher order terms, the first step is to
look for the first terms of (\ref{whatlook41}) and (\ref{whatlook42}),i.e.,
\beann
\pi(x)&=& \frac{1}{2}x'Px + \ldots\\
u(x)&=& Kx + \ldots.
\eeann
The homogeneous polynomial of higher degrees of (\ref{whatlook41}) and (\ref{whatlook42}) are eliminated as the only terms gathered are the quadratic terms in (\ref{dpe1}) and the linear terms in (\ref{dpe2}) when (\ref{whatlook41}) and (\ref{whatlook42}) are substituted in the DPE.  Then, the first level is exactly the linear-quadratic regulator case.  Hence, we get the solutions (\ref{ric1}) and (\ref{kappa1}).

\subsection{Higher Levels:  Higher Degree Terms}
In the second level, we look for the second terms of the expansions,
\beann
\pi(x)&=& \frac{1}{2}x'Px + \pi^{[3]}(x) + \ldots\\
u(x)&=& Kx + \kappa^{[2]}(x) + \ldots,   
\eeann
namely $\pi^{[3]}(x)$ and $\kappa^{[2]}(x)$.  Before the method is applied, we reexpressed the $f(x,u)$ and $l(x,u)$ with the first level feedback solution; i.e.,
\beann
\bar{f}(x,u)&=&f(x,Kx + u)\\
\bar{l}(x,u)&=&l(x,Kx + u)
\eeann
These functions have power series representation through terms of $3^{rd}$ and $4^{th}$ of the form
\beann
\bar{f}(x,u)&=&(A+BK)x + Bu + \bar{f}^{[2]}(x,u) + \ldots \\
\bar{l}(x,u)&=&\frac{1}{2}(x'Qx + 2x'SKx + x'K'RKx) + x'Su + u'Ru + \bar{l}^{[3]}(x,u) + \ldots 
\eeann
Again repeating the same process but this time grouping all the cubic terms of (\ref{dpe1}) gives
\bea \label{pi3}
\pi^{[3]}(x)-\pi^{[3]}((A+BK)x) = \frac{1}{2} x'(A+BK)'\bar{f}^{[2]}(x,0) \nonumber\\
  + \frac{1}{2} \bar{f}^{[2]'}(x,0)P(A+BK)x + \bar{l}^{[3]}(x,0)
\eea
where the terms
\bea \label{cancel1}
x'Su + \frac{1}{2}x'K'Ru + \frac{1}{2}u'RKx +\frac{1}{2}x'(A+BK)'PBu) +\frac{1}{2}u'B'P(A+BK)x
\eea
drops out of (\ref{pi3}) as it is zero by (\ref{prekappa2}) and quadratic terms of (\ref{dpe2}) yields
\bea \label{kappa2}
\kappa^{[2]}(x)&=&-(B'PB + R)^{-1}\cdot\\
            &&\left[ \frac{\partial \bar{f}^{[2]'}} {\partial u}(x,0)P(A+BK)x + B'P\bar{f}^{[2]}(x,0)+B' \frac{\partial \pi^{[3]}}{\partial f}((A+BK)x)\right. \nonumber \\
&&\left. + \frac{\partial l^{[3]}}{\partial x}(x,0)\right] \nonumber
\eea

Cancelling these terms (\ref{cancel1}) is crucial because it introduces a block triangular structure within the levels.  If $\kappa^{[2]}$ appears in (\ref{pi3}), then the system is coupled at the second level.  Similar condition applies for any $d$ level.  The block triangular form facilitates solving the system.

The $d^{th}$ Level:

We discuss the system of equations at the $d^{th}$ level.  In the $d^{th}$, the two linear systems are obtained by collecting the $d+1$ degree terms of the \ref{dpe1} and $d$ degree terms of the \ref{dpe2}.  Suppose we have solved through the $d-1^{th}$ level, we incorporate the feedback up to the degree $d-1$ into the dynamics and cost to obtain the updated $\bar{f}(x,u)$ and $\bar{l}(x,u)$; i.e.,
\beann
\bar{f}(x,u)&=& f(x,Kx + \kappa^{[2]}(x) + \kappa^{[3]}(x) + \ldots + \kappa^{[d-1]}(x) + u)\\
\bar{l}(x,u)&=& l(x,Kx + \kappa^{[2]}(x) + \kappa^{[3]}(x) + \ldots + \kappa^{[d-1]}(x) + u).
\eeann

  In the $d^{th}$ level, the corresponding cancellation is
\bea \label{cancel2}
x'Su + \frac{1}{2}x'K'Ru + \frac{1}{2}u'RKx 
+\frac{1}{2}x'(A+BK)'PBu +  \frac{1}{2}u'B'P(A+BK)x=0
\eea
After the substitution and collection of all the $d+1$ degree terms of (\ref{dpe1}) and $d$ degree terms of (\ref{dpe2}), we have the equations
\bea \label{pid+1}
\pi^{[d+1]}(x)-\pi^{[d+1]}((A+BK)x) &=& \frac{1}{2} \Bigg[ x'(A+BK)'\bar{f}^{[d]}(x,0) \\
&&  + \sum_{j=2}^{d-1} [\kappa^{[j]'}B' + \bar{f}^{[j]'}(x,0)  ]P[B\kappa^{[d-1-j]}+ \bar{f}^{[d-1-j]}(x,0)] \nonumber \\
&& + \bar{f}^{[d]'}(x,0)P(A+BK)x \Bigg]+ \tilde{\pi}^{[d+1]}(\bar{f})+ \bar{l}^{[d+1]}(x,0) \nonumber
\eea
and 
\bea \label{kappad}
\kappa^{[d]}(x)&=&-(B'PB + R)^{-1}\cdot \nonumber\\
            &&\Bigg[ \frac{\partial f^{[d]'}} {\partial u}(x,0)P(A+BK)x \\
&& + \sum_{j=2}^d \frac{\partial \bar{f}^{[j]}}{\partial u}(PB\kappa^{[d-j]}+\frac{\partial \pi^{[d+1-j]}}{\partial f} \nonumber + \frac{\partial \bar{l}^{[d+1]}}{\partial u}(x,0) \Bigg] \nonumber
\eea
where
$\tilde{\pi}^{[d+1]}(\bar{f})$ are terms of degree $d+1$ of the homogeneous polynomials $\pi^{[m]}$ where $2 \leq m \leq d+1$.  These terms are generated by the input of the feedback control into the DPE. 

\subsection{Solvability of the System of Equations}
We know discuss the solutions of the system of equations obtained from level $\geq 2$.  
For the case $d=2$, we show that $\pi^{[3]}$ exists as a solution for (\ref{pid+1}).  First, the equation (\ref{pi3}) is linear since
\bea \label{pimap}
\pi^{[3]}(x) \mapsto \pi^{[3]}(x) + \pi^{[3]}((A+BK)x)
\eea
is linear
\begin{lemma} \label{speclem}
Given that $\vert \sigma(A+BK) \vert < 1$.  There exists a unique solution of the homogeneous polynomial $\pi^{[3]}$.
\end{lemma}
\begin{proof}
For simplicity, we assume that $(A+BK)$ has simple eigenvalues $\mu_i$, i.e.
\[w_i(A+BK)=\mu_i w_i.\]
where $w_i \in \mathbb{R}^{1 \times n}$.
Since the system is stable,
\bea \label{stabcond}
\vert \mu_i \vert < 1.
\eea
The polynomial $\pi^{[3]}(x)$ has the following form:
\[\sum_{i_1}^n \sum_{i_2}^n \sum_{i_3}^n c_{i_1,i_2,i_3} (w_{i_1}x_{i_1}) (w_{i_2}x_{i_2}) (w_{i_3}x_{i_3}).\]
where $(w_{i_1}x_{i_1}) (w_{i_2}x_{i_2}) (w_{i_3}x_{i_3})$ are the basis for all cubic polynomials.
It follows that 
\bea \label{lem1}
\pi^{[3]}(x)- \pi^{[3]}((A+BK)x)=& \nonumber\\
                                &\sum_{i_1}^n \sum_{i_2}^n \sum_{i_3}^n c_{i_1,i_2,i_3}(1-\mu_{i_1} \mu_{i_2} \mu_{i_3})(w_{i_1}x_{i_1}) (w_{i_2}x_{i_2}) (w_{i_3}x_{i_3}) \nonumber\\
\eea
The RHS of (\ref{pi3}) is
\bea \label{lem2}
 \sum_{i_1}^n \sum_{i_2}^n \sum_{i_3}^n d_{i_1,i_2,i_3} (w_{i_1}x_{i_1})  (w_{i_2}x_{i_2}) (w_{i_3}x_{i_3}).
\eea
Together with (\ref{lem1}) and (\ref{lem2}), we express the coefficients as
\[c_{i_1,i_2,i_3}= \frac {d_{i_1,i_2,i_3}}{(1-\mu_{i_1} \mu_{i_2} \mu_{i_3})}\]
where
\[(1-\mu_{i_1} \mu_{i_2} \mu_{i_3}) \neq 0\]
by (\ref{stabcond}).  Hence, there exists a unique homogeneous polynomial $\pi^{[3]}(x)$.
\end{proof}
It follows that the linear map (\ref{pimap}) is invertible.  For any level d, the homogeneous polynomial $\pi^{[d]}(x)$ exists.  The proof on the existence for higher order levels follows in the same fashion as above. 

For any $d$ level, the RHS of equation (\ref{kappad}) contains a multiplicative term
\[(B'PB + R)^{-1}.\]
The inverse matrix exists because
\[(B'PB + R)\]
is positive definite since $R$ is positive definite.  The rest of the terms on the RHS of (\ref{kappad}) are known from the previous levels.  Thus, $\kappa^{[d]}$ exists for any level $d$.  Therefore, all of the linear equations for degree 3 or higher are solvable.



%% file: chapter4_39.tex
In the previous chapter, we found polynomials of degree $(r-1)$ and $r$, namely the optimal control and optimal cost, respectively, that satisfy the DPE.  We will prove the existence of smooth solutions to the DPE which have the same Taylor series expansions as of the formal solutions previously discussed in Chapter $3$.

\section{The Main Result}
Here we state our main result.  
\begin{theorem} \label{thebigkahuna}
Suppose that the dynamics and cost are $C^{r-1}$ and $C^{r}$, respectively and the linear part of the nonlinear Hamiltonian system (\ref{nonham}) is stabilizable and detectable around zero.  Then there exists $C^{r-1}(\mathbb{R}^n)$ optimal control and $C^{r}(\mathbb{R}^n)$ optimal cost solutions to the DPE locally around zero.
\end{theorem}
If we let $k \longrightarrow \infty$ in Theorem (\ref{thebigkahuna}), we then have the following corollary.
\begin{corollary}
Suppose that the dynamics and cost are real analytic functions on $\mathbb{R}^{2n}$ and the linear part of the nonlinear Hamiltonian system (\ref{nonham}) is stabilizable and detectable around zero.  Then the power series solutions converges to solution of the DPE locally around zero.
\end{corollary}

\section{PMP and Hamiltonian Dynamics}
We look at the Pontryagin Maximum Principle (PMP) as it presents the Hamilton difference equations satisfied by the optimal trajectories.  The PMP gives a necessary condition for a control to be optimal.  Along with the optimal control problem formulation, there is an associated nonlinear Hamiltonian,
\bea \label {theham}
H(x,u,\lambda^+)=\lambda^{+'}f(x,u) + l(x,u)
\eea
where $\lambda^+=\lambda_{k+1}$.
The PMP states the following:
\begin{theorem}
If $x_k$ and $u_k$ are optimal for $k \in 0, 1, 2, \ldots$, then there
exists $\lambda_k \neq 0$ for $k \in 0, 1, 2, \ldots$ such that
\bea
x^+&=&\frac{\partial H}{\partial \lambda^+}(x,u,\lambda^+) \label{hode1} \\
\lambda&=& \frac{\partial H}{\partial x}(x,u,\lambda^+) \label{hode2} \\
u^*&=&\mbox{arg}\min_{u}H(x,u,\lambda^+). \label{hode3}
\eea
\end{theorem}
Thus, the minimizer of the nonlinear Hamiltonian evaluated at the optimal
$x$ and $\lambda^+$ is the optimal control $u$ amongst all admissible
controls $v$.  Note that $u^*(x,\lambda^+) \in \mathcal{C}^{r-1}$ since $f \in \mathcal{C}^{r-1}$ and $l \in \mathcal{C}^r$ in (\ref{theham}).  We assume that $H$ is convex in u to guarantee a unique optimal control.  So the PMP gives the existence of the optimal control.  However, the equation (\ref{hode3}) is an optimal control that has $x$ and $\lambda$ as its independent variables.  We must show that $\lambda$ is function of $x$ to complete this part of the proof.
\section{Forward Hamiltonian Dynamics}
If we take the (\ref{hode1}), (\ref{hode2}), and (\ref{hode3}) of the PMP and linearize the equations around zero.  The linearized sytem that we obtain is exactly the linear Hamiltonian system.  The corresponding Hamiltonian is
\[H(x,\lambda^+,u)=\lambda^{+'}(Ax + Bu) + \frac{1}{2} x^{'} Q x + x^{'}S u
+ \frac{1}{2}u^{'}R u\]
and the system is
\bea \label{hamsys}
\left[ \begin{array}{c}
x^+ \\
\lambda \\
\end{array}\right]=\mathbb{H}
\left[\begin{array}{c}
x \\
\lambda^+ \\
\end{array}\right]
\eea
where
\beann
\mathbb{H}=
\left[ \begin{array}{cc}
A - BR^{-1}S'  & -BR^{-1}B' \\
Q -SR^{-1}S'   &  A' - SR^{-1}B'\\
\end{array}\right]
\eeann
is the associated Hamiltonian matrix. 

Notice the opposing directions of the propagation of the state and costate dynamics in (\ref{hamsys}).  This presents difficulty in studying some Hamiltonian properties, but it is addressed in Chapter 4.  By assuming invertibility of the matrix $A-SR^{-1}B'$, we can rewrite the Hamilton equations (\ref{hamsys}) so that both equations propagate in the same direction, in particular the direction of increasing time.
The forward recursion,
\bea \label{hamfor}
\left[ \begin{array}{c}
x^+ \\
\lambda^+ \\
\end{array}\right]= \mathbb{H}^F
\left[ \begin{array}{c}
x \\
\lambda \\
\end{array}\right]
\eea
where
\bea \label{fhmat}
\mathbb{H}^F=
\left[ \begin{array}{cc}
\scriptstyle (A-BR^{-1}S')-BR^{-1}B'(A'-SR^{-1}B')^{-1}(Q-SR^{-1}S')  &
\scriptstyle BR^{-1}B'(A'-SR^{-1}B')^{-1} \\
\scriptstyle -(A'-SR^{-1}B')^{-1}(Q-SR^{-1}S')   & \scriptstyle
(A'-SR^{-1}B')^{-1}\\
\end{array}\right]
\eea
is directly derived from the system (\ref{hamsys}).  The forward recursion describes the flow of the state and costate as $t$ approaches to infinity.  We refer to $\mathbb{H}^F$ as the forward Hamiltonian matrix.  Thus, the dynamics (\ref{hamfor}) rewritten in the forward direction relies on the invertibility of $A'-SR^{-1}B'$.

The existence of $(A'-SR^{-1}B')^{-1}$ has the following consequences:

\begin{lemma} \label{nozero}
If $A'-SR^{-1}B'$ is invertible, then $\mathbb{H}^F$ exists and is invertible.
\end{lemma}
\begin{proof} \label{noopnocl}
First, we observe in (\ref{fhmat}) that $(A'-SR^{-1}B')^{-1}$ appears.  Then, if $(A'-SR^{-1}B')^{-1}$ exists so does $\mathbb{H}^F$.  Explicit calculation gives

\bea \label{bhmat}
(\mathbb{H}^F)^{-1}=
\left[ \begin{array}{cc}
\scriptstyle (A-BR^{-1}S')^{-1} & \scriptstyle (A-BR^{-1}S')BR^{-1}B'\\
\scriptstyle (Q-SR^{-1}S')(A-BR^{-1}S')^{-1} & \scriptstyle
(Q-SR^{-1}S')(A-BR^{-1}S')BR^{-1}B' + (A-BR^{-1}S')'\\
\end{array}\right]
\eea
The matrix (\ref{bhmat}) satisfies the backward hamiltonian dynamics.  Thus, $\mathbb{H}^B=(\mathbb{H}^F)^{-1}$ where 
\bea \label{hamback}
\left[ \begin{array}{c}
x \\
\lambda \\
\end{array}\right]= \mathbb{H}^B
\left[ \begin{array}{c}
x^+ \\
\lambda^+ \\
\end{array}\right].
\eea
We see that $\mathbb{H}^B$ exists since $(A-BR^{-1}S')$ is invertible.
Hence, $\mathbb{H}^F$ is invertible.
\end{proof}
It follows that $\mathbb{H}^F$ does not have zero as an eigenvalue if
$A'-SR^{-1}B'$ is invertible.  Also, the Lemma (\ref{nozero}) implies that if $\mathbb{H}^F$ has a zero eigenvalue, then so does the spectrum of $A-SR^{-1}B'$.  In \cite{AM97}, if $(A,B)$ is stabilizable and $(Q^{1/2},A)$ is detectable, then the system is hyperbolic; i.e., none of the eigenvalues lie on the unit circle.  If zero is an eigenvalue, then infinity is also an eigenvalue because of the hyperbolicity of the system.  We refer to these eigenvalues as infinite eigenvalues.  This issue of infinite eigenvalues is discussed in Chapter 4.  For the time being we will require that the forward Hamiltonian matrix exists and be invertible.  For the rest of the chapter, we assume the invertibility of $A'-SR^{-1}B'$ since the existence of $\mathbb{H}^F$ depends on the inverse of $A'-SR^{-1}B'$.

The forward Hamiltonian matrix has $2n$ eigenvalues; n are stable, i.e., these eigenvalues lie inside the unit
circle.  The other $n$ eigenvalues are unstable and are positioned outside
the unit circle.  We will see in Section 3.4.3 that the eigenstructure is hyperbolic.

Recall the LQR problem in Chapter $2$ of minimizing

\bea \label{lqrprob}
&&\min_u \sum_{j=0}^{\infty}  \frac{1}{2}x_j'Qx_j + x_j'Su_j +
\frac{1}{2}u'_jRu_j \nonumber \\
&&\mbox{subject to the dynamics}\\
&&x^+=Ax + Bu \nonumber \\
&&x_0=x(0). \nonumber
\eea

\begin{lemma} \label{myzero}
Suppose $u=Kx$ is the optimal control.  If zero is an eigenvalue of $(A-BR^{-1}S')$, then zero is in $\sigma(A+BK)$. 
\end{lemma}
\begin{proof}
With the transformation $u=Lx + v$ the optimal control problem is changed to the following problem:
\bea \label{woxterm}
&&\min_u \sum_{j=0}^{\infty}  \frac{1}{2}x_j'(Q + L'RL + S)x_j +
\frac{1}{2}v'_jRv_j \nonumber\\
&&\mbox{subject to the dynamics}\\
&&x^+=(A-BR^{-1}S')x + Bv \nonumber \\
&&x_0=x(0) \nonumber
\eea
where $L=-R^{-1}S'$.  The transformation removes the cross term in the cost.  For the modified problem (\ref{woxterm}), we have $v=\bar{K}x$ is the optimal control, then $K=L+\bar{K}$.  The transformation removes the cross term in the cost.  Then, we choose the eigenvector $x_0 \neq 0$ of $A-BR^{-1}S'$ such that its corresponding eigenvalue is zero; i.e.,
\bea \label{zeroeig}
(A-BR^{-1}S')x_0=0
\eea
The cost starting from $x_0$ is
\bea \label{startcost}
\pi(x_0)=x_0'(Q+L'RL + S)x_0.
\eea
No additional control $v$ is exercised since the dynamics reaches zero
at one time step.  Therefore, $v=\bar{K}x_0=0$ and the optimal cost is the cost (\ref{startcost}) at the zero time.  We have that
\bea \label{lek} 
Lx_0=Kx_0
\eea
for $x_0 \neq 0$ that satisfies (\ref{zeroeig}) since $Kx_0=(L + \bar{K})x_0$.

Returning to (\ref{lqrprob}), the closed-looped spectrum is $\sigma(A+BK)$.   Then, from (\ref{lek})
\bea \label{lek2}
(A+BL)x_0=(A+BK)x_0.
\eea
Since the LHS of (\ref{lek2}) is equal to 0, we have
\[(A+BK)x_0=0\]
for $x_0 \neq 0$.  Thus, $0 \in \sigma(A+BK)$.
Hence, zero is a closed-loop eigenvalue.
\end{proof}
By Lemma (\ref{myzero}), zero is not a closed-loop eigenvalue.  This implies that zero is not an eigenvalue of $A-BR^{-1}S'$ if $u=Kx$ is assumed to be the optimal control. 

\section{Properties of the Nonlinear Hamiltonian Dynamics}
It is imperative to look at some of the properties that the nonlinear Hamiltonian dynamics possesses as it is needed in the proof of the existence of the optimal cost.  We look at its tangent dynamics, sympletic form, and eigenstructure.

Recall that the PMP gives the nonlinear Hamilton difference equations (\ref{hode1}), (\ref{hode2}) satisfied by the optimal trajectories.  These equations are the nonlinear Hamiltonian dynamics corresponding to the optimal control problem.  We calculate the forward nonlinear Hamiltonian dynamics from (\ref{hode1}), (\ref{hode2}). Since the forward linear Hamiltonian dynamics exists at zero and by the Implict Function Theorem, the forward nonlinear dynamics also exists in the neigborhood around 0.  The forward nonlinear Hamiltonian dynamics,
\bea \label {nonhamforG}
      x^+&=&G_1(x,\lambda)\\
\lambda^+&=&G_2(x,\lambda)
\eea
evolves in the direction of increasing time.

\subsection{Tangent Dynamics}
We now introduce the idea of tangent vectors to $\mathcal{M}=\{(x,\lambda)\vert x\in \mathbb{R}^n,~\lambda+ \in \mathbb{R}^n)\}$; i.e., $\mathcal{M}=\mathbb{R}^{2n}$.  First, we denote the tangent vectors 
\beann
v=\left[ \begin{array}{c}
\delta x\\
\delta \lambda \\
\end{array}\right].
\eeann
The set of tangent vector to $\mathcal{M}$ at $(x,\lambda)$ forms a vector space $T_{x,\lambda}\mathcal{M}$, the tangent space to $\mathcal{M}$ at $(x,\lambda) \in \mathcal{M}$.  In addition, the tangent bundle of $\mathcal{M}$, denoted by $T\mathcal{M}$, is a differential manifold where
\bea \label{tanbundle}
T\mathcal{M}=\bigcup_{x,\lambda^+\in \mathcal{M}}T_{x,\lambda^+}\mathcal{M}.
\eea
The local coordinate system on $T\mathcal{M}$ is $4n$ numbers $x_1,\ldots,x_n,\lambda_1,\ldots,\lambda_n$ and\\
$\delta x_1,\ldots,\delta x_n,\delta\lambda_1,\ldots,\delta \lambda_n$ where the former is the local coordinate on $\mathcal{M}$ and the latter is the components of the tangent vector.

The tangent dynamics describes the flow of the tangent vector at each sequential point of a map.
\begin{definition} \label{Dtangentdynamics}
Suppose $x_{k+1}=f(x_k)$ with known $x_0$.  Let $\delta x_0$ be a tangent vector at $x_0$.  Then we define the tangent dynamics around the trajectory $x_k$ as
\[\delta x_{k+1}=\frac{\partial f}{\partial x}(x_k)\delta x_k.\]
\end{definition}
When the dynamics is linear; i.e., $x_{k+1}=Ax_k$, then the tangent dynamics is $v_{k+1}=Av_k$ because $\frac{\partial f}{\partial x}(x_k)=A$.  Thus, the linear Hamiltonian dynamics (\ref{hamfor}) has
\bea \label{hamfor}
\left[ \begin{array}{c}
\delta x^+ \\
\delta \lambda^+ \\
\end{array}\right]= \mathbb{H}^F
\left[ \begin{array}{c}
\delta x \\
\delta \lambda \\
\end{array}\right]
\eea
as it tangent dynamics.  We now find the nonlinear tangent dynamics as we have a nonlinear system.  Note the increase in difficulty in finding the tangent dynamics of the forward nonlinear Hamilton dynamics because of the mixed directions of the propagation.  Linearizing the the nonlinear system at the
trajectories $(x,\lambda^+)$ invokes a linear perturbation. We perturb
$x$ and $\lambda$ and perform the first variation on the
\[ \frac{\partial H}{\partial (x ,\lambda^+)}(x,\lambda^+).\]
The nonlinear Hamiltonian dynamics from PMP is
\bea \label{nonham}
x^+     &=& \frac{\partial H}{\partial \lambda^+} (x,\lambda^+)\nonumber\\
\lambda &=& \frac{\partial H}{\partial x} (x,\lambda^+).
\eea
Suppose we replace the parameters
\beann
x         &\mapsto& x + \delta x \\
\lambda^+ &\mapsto& \lambda^+ + \delta \lambda^+.
\eeann
Then,
\bea \label{pert1}
x^+ + \delta x^+          &\approx& \frac{\partial H}{\partial \lambda^+} (x
+ \delta x, \lambda^+ + \delta \lambda^+) \nonumber \\
\lambda  + \delta \lambda &\approx& \frac{\partial H}{\partial x} (x +
\delta x, \lambda^+ + \delta \lambda^+).
\eea
Subtracting (\ref{nonham}) from (\ref{pert1}), we have
\bea \label{pertsys}
\delta x^+     &\approx& \frac{\partial H}{\partial \lambda^+} (x + \delta
x, \lambda^+ + \delta \lambda^+) - \frac{\partial H}{\partial \lambda^+}
(x,\lambda^+) \label{pertsys1} \\
\delta \lambda &\approx& \frac{\partial H}{\partial x} (x + \delta x,
\lambda^+ + \delta \lambda^+) - \frac{\partial H}{\partial x}
(x,\lambda^+) \label{pertsys2}.
\eea
Expanding the RHS of (\ref{pertsys}) and (\ref{pertsys2}) in Taylor series,
we get
\bea
\delta x^+     &=& \frac{\partial^2 H}{\partial \lambda^+ \partial x}
(x,\lambda^+)\delta x + \frac{\partial^2 H}{\partial^2 \lambda^+}
(x,\lambda^+) \delta \lambda^+ \nonumber\\
\delta \lambda &=& \frac{\partial^2 H}{\partial^2 x} (x,\lambda^+) \delta
x + \frac{\partial^2 H}{\partial \lambda^+ \partial x} (x,\lambda^+)\delta
\lambda^+.
\eea
Equivalently,
\bea \label{nonhampert}
\left[ \begin{array}{c}
\delta x^+\\
\delta \lambda \\
\end{array}\right]=
\mathbb{H}_{\delta,k}(x,\lambda^+)
\left[ \begin{array}{c}
\delta x\\
\delta \lambda^+\\
\end{array}\right].
\eea
where
\beann
\mathbb{H}_{\delta,k}(x,\lambda^+)=
\left[ \begin{array}{cc}
H_{\lambda^+ x} & H_{\lambda^+ \lambda^+}\\
H_{x x}         & H_{x \lambda^+}
\end{array}\right](x,\lambda^+).
\eeann
The subscripts of  $H_{x \lambda^+}$ denote partial derivatives.
We look at the forward perturbed system.  The tangent dynamics in the forward time is
\bea \label{nonhampertfor}
\left[ \begin{array}{c}
\delta x^+ \\
\delta \lambda^+ \\
\end{array}\right]= \mathbb{H}_{\delta,k}^F(x,\lambda^+)
\left[ \begin{array}{c}
\delta x \\
\delta \lambda \\
\end{array}\right]
\eea
where
\beann
\mathbb{H}_{\delta,k}^F(x,\lambda^+)=
\left[ \begin{array}{cc}
H_{\lambda^+ x}-H_{\lambda^+ \lambda^+}H_{\lambda^+ \lambda^+}^{-1}H_{xx}  &
H_{\lambda^+ \lambda^+}H^{-1}_{x \lambda}\\
-H^{-1}_{x \lambda^+}H_{xx} & H^{-1}_{x \lambda^+}
\end{array}\right](x,\lambda^+).
\eeann
To have the tangent dynamics written in forward time, we assume the
invertibility of $H_{\lambda^+ x}$.  When $H_{\lambda^+ x}(x,\lambda^+)$ is
evaluated at $\vec{0}$, we get the matrix $A'-SR^{-1}B'$ which we have assumed to be invertible.  Hence, $H_{\lambda^+ x}$ is invertible for small $x$ and $\lambda^+$.

\subsection{The Standard Symplectic Form}
We define a nondegenerate and a bilinear symplectic two-form $\Omega:
T_{(x,\lambda)}\mathcal{M} \times T_{(x,\lambda)}\mathcal{M} \mapsto \mathbb{R}$,
\beann
\Omega(v,w)=v'Jw
\mbox{   and   }J=
\left[ \begin{array}{cc}
0  & I \\
-I & 0 \\
\end{array}\right]
\eeann
where
\[\Omega(v,w)=-\Omega(w,v)\]
and
\[ v=\left[ \begin{array}{c}
\delta x\\
\delta \lambda \\
\end{array}\right],~~~ w=\left[ \begin{array}{c}
\widetilde{\delta x}\\
\widetilde{\delta \lambda} \\
\end{array}\right]\]
$(x,\lambda^+)\in \mathcal{M}$ and $(v,w) \in T_{(x,\lambda)}\mathcal{M}$. The matrix $J$ is called the symplectic matrix.

With the system (\ref{nonhampertfor}), it can be shown through calculations that 
\[\mathbb{H}_{\delta,k}^{F'}J\mathbb{H}_{\delta,k}^{F}=J.\]
Then, we have that\\
\beann
\left[ \begin{array}{cc}
H^{T}_{\lambda^+ x}-H^{T}_{xx}H^{-1}_{x \lambda^+}H^{T}_{\lambda^+
\lambda^+} & -H^{T}_{xx}H^{-T}_{x \lambda^+}\\
H^{-T}_{x \lambda^+}H^{T}_{\lambda^+ \lambda^+} & H^{-T}_{x \lambda^+}
\end{array}\right]
J
\left[ \begin{array}{cc}
H_{\lambda^+ x}-H_{\lambda^+ \lambda^+}H_{\lambda^+ \lambda^+}^{-1}H_{xx}  &
H_{\lambda^+ \lambda^+}H^{-1}_{x \lambda}\\
-H^{-1}_{x \lambda^+}H_{xx} & H^{-1}_{x \lambda^+}
\end{array}\right](x,\lambda^+)
\eeann
is equal to
\bea \label{symplectic2}
\left[ \begin{array}{cc}
\mathbb{A}_{11}  & \mathbb{A}_{12} \\
\mathbb{A}_{21}  & \mathbb{A}_{22} \\
\end{array}\right]
\eea
where
\beann
\mathbb{A}_{11} &=& H^{T}_{xx}H^{-T}_{x \lambda^+}H_{\lambda^+ x}
-H^T_{xx}H^{-T}_{x \lambda^+}H_{\lambda^+ \lambda^+}H^{-1}_{x
\lambda^+}H_{xx} - H^T_{\lambda^+ x}H^{-1}_{x \lambda^+}H_{xx} +
H^T_{xx}H^{-T}_{x \lambda^+}H^T_{\lambda^+ \lambda^+}H^{-1}_{x
\lambda^+}H_{xx} \\
\mathbb{A}_{12} &=& H^T_{xx}H^{-T}_{x \lambda^+}H_{\lambda^+
\lambda^+}H^{-1}_{x \lambda^+} + H^T_{\lambda^+ x}H^{-1}_{x \lambda^+} -
H^T_{xx}H^{-T}_{x \lambda^+}H^T_{\lambda^+ \lambda^+}H^{-1}_{x \lambda^+}\\
\mathbb{A}_{21} &=& -H^{-T}_{x \lambda^+}H_{\lambda^+ x} + H^{-T}_{x
\lambda^+}H_{\lambda^+ \lambda^+}H^{-1}_{x \lambda^+}H_{xx} - H^{-T}_{x
\lambda^+}H_{\lambda^+ \lambda^+}H^{-1}_{x \lambda^+}H_{xx}\\
\mathbb{A}_{22} &=&  -H^{-T}_{x \lambda^+}H_{\lambda^+ \lambda^+}H^{-1}_{x
\lambda^+} +H^{-T}_{x \lambda^+}H^{T}_{\lambda^+ \lambda^+}H^{-1}_{x
\lambda^+}.
\eeann
Since the matrices $H_{\lambda^+ x}=H_{x \lambda^+},~H_{x
x},~\mbox{and}~H_{\lambda^+ \lambda^+}$ are symmetric, then
(\ref{symplectic2}) reduces to $J$.  Thus,
\bea \label{symplectic3}
\mathbb{H}_{\delta,k}^{F'}J\mathbb{H}_{\delta,k}^{F}=J
\eea
Hence, the nonlinear Hamilton dynamics preserves the standard symplectic form. Suppose that $(v,w)$ are two tangent vectors that satisfy (\ref{nonhampertfor}); i.e., $v^+=\mathbb{H}_{\delta,k}^Fv$ and $w^+=\mathbb{H}_{\delta,k}^Fw$.  Let $v$, $w$ be two sequences propagating according to(\ref{nonhampertfor}).  Then, 
\beann
\Omega(v^{+},w^{+})&=&v^{'+}J w^+\\
                       &=&v' \mathbb{H}_{\delta,k}^{F'} J
\mathbb{H}_{\delta,k}^F w \\
                       &=&v'Jw\\
                       &=& \Omega(v,w).
\eeann
\begin{figure}[tp]
\centering
\scalebox{.3}{\includegraphics{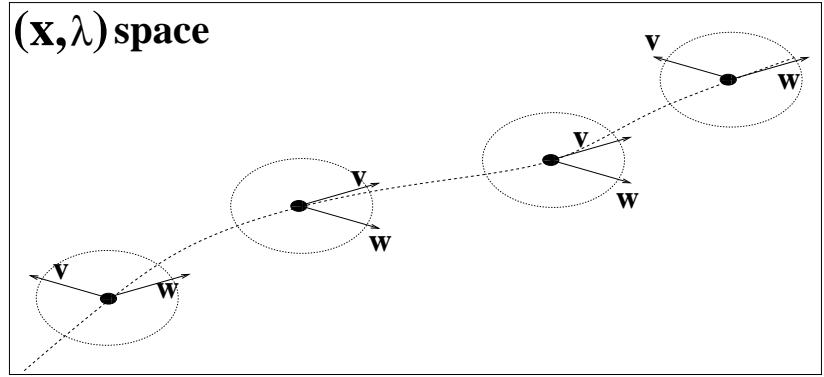}}
\caption{At each sequence of the map, $\Omega(v,w)$ is constant.}
\label{fig1}
\end{figure}
In other words, for any two sequences of tangent vector under the tangent dynamics (\ref{nonhampert}), the value of the $\Omega$ is unchanged at every point of the map.  See fig.(\ref{fig1}).

As a consequence, if we evaluate (\ref{symplectic3}) at 0, then at the linear level 
\bea \label{symplectic}
\mathbb{H}^{F'}J \mathbb{H}^F=J.
\eea
where $\mathbb{H}^{F'}$ is the matrix in (\ref{hamfor}).  The calculations are as follows:
For simplicity we let
\beann
\alpha &=& A-BR^{-1}S'\\
\beta  &=& BR^{-1}B'  \\
\gamma &=& Q- SR^{-1}S'. \\
\eeann
Then,
\beann
\mathbb{H}^F=
\left[ \begin{array}{cc}
\alpha - \beta \alpha^{-T} \gamma & \beta \alpha^{-T} \\
-\alpha^{-T} \gamma              & \alpha^{-T}
\end{array}\right]
\eeann
and
\beann
\mathbb{H}^{F'}=
\left[ \begin{array}{cc}
\alpha^{T}- \gamma \alpha^{-1} \beta & -\gamma \alpha^{-T} \\
\alpha^{-1}\beta                     & \alpha^{-1}
\end{array}\right].
\eeann
We denote $\alpha^{-T}$ as the inverse of the transpose of $\alpha$.  Also,
the matrices $\beta$ and $\gamma$ are symmetric. The equation
(\ref{symplectic}) is
\beann
\left[ \begin{array}{cc}
\alpha^{T}- \gamma \alpha^{-T} \beta & -\gamma \alpha^{-T} \\
\alpha^{-1}\beta                     & \alpha^{-1}
\end{array}\right]
\left[ \begin{array}{cc}
0  & I \\
-I & 0 \\
\end{array}\right]
\left[ \begin{array}{cc}
\alpha - \beta \alpha^{-T} \gamma & \beta \alpha^{-T} \\
-\alpha^{-T} \gamma               & \alpha^{-T}
\end{array}\right]
\eeann
which is simplified to
\beann
\left[ \begin{array}{cc}
\gamma - \gamma \alpha^{-1} \beta \alpha^{-T} \gamma - \gamma + \gamma
\alpha^{-1} \beta \alpha^{-T} \gamma & I + \gamma \alpha^{-1} \beta
\alpha^{-T} - \gamma \alpha^{-1} \beta \alpha^{-T}\\
-I + \alpha^{-1} \beta \alpha^{-T} \gamma - \alpha^{-1} \beta \alpha^{-T}
\gamma &
- \alpha^{-1} \beta \alpha^{-T} +\alpha^{-1} \beta \alpha^{-T}\\
\end{array}\right]=J.
\eeann\\
Therefore,
\[\mathbb{H}^{F'}J \mathbb{H}^F=J\]
holds.

Similarly, if we take tangent vectors satisfying the linear tangent dynamics
\[v^+=\mathbb{H}^F v \mbox{   and   } w^+=\mathbb{H}^F w, \]
then
\beann
\Omega(v^+,w^+)= \Omega(v,w).
\eeann
Thus, the two-form $\Omega(v,w)$ is unchanged when evaluated at the tangent
vectors $v$ and $w$ of every point of the map under the linear Hamiltonian dynamics.  

\subsection{Eigenstructure}

Now with (\ref{symplectic}), the eigenstructure of $\mathbb{H}^F$ is as follow:
\begin{theorem} \label{eigenthm}
Suppose $\mathcal{A}^{'}J \mathcal{A}=J$. If $\mu \in \sigma
(\mathcal{A})$, then so do $
\frac{1}{\mu},\;\bar{\mu},$ and $ \frac{1}{\bar{\mu}} \in \sigma
(\mathcal{A})$.
\end{theorem}
\begin{proof}
Suppose $(\delta x,\delta \lambda)'$ is an eigenvector of $\mathcal{A}$ for some
eigenvalue $\mu$.
\beann
\mathcal{A}
\left[ \begin{array}{c}
\delta x\\
\delta \lambda\\
\end{array}\right]=
\mu
\left[ \begin{array}{c}
\delta x\\
\delta \lambda\\
\end{array}\right]
\eeann
Also, note that
\beann
J
\left[ \begin{array}{c}
\delta x\\
\delta \lambda\\
\end{array}\right]=
\left[ \begin{array}{c}
-\delta \lambda\\
\delta x\\
\end{array}\right].
\eeann
Then,
\beann
\mathcal{A}{'}J\mathcal{A}
\left[ \begin{array}{c}
\delta x\\
\delta \lambda\\
\end{array}\right]=
J
\left[ \begin{array}{c}
\delta x\\
\delta \lambda\\
\end{array}\right]&\Longleftrightarrow&
\mu \mathcal{A}^{'}J
\left[ \begin{array}{c}
\delta x\\
\delta \lambda\\
\end{array}\right]=
\left[ \begin{array}{c}
-\delta \lambda\\
\delta x\\
\end{array}\right]\\
&\Longleftrightarrow& \mu \mathcal{A}^{'}
\left[ \begin{array}{c}
-\delta \lambda\\
\delta x\\
\end{array}\right]=
\left[ \begin{array}{c}
-\delta \lambda\\
\delta x\\
\end{array}\right]\\
&\Longleftrightarrow& \mathcal{A}^{'}
\left[ \begin{array}{c}
-\delta \lambda\\
\delta x\\
\end{array}\right]=\frac{1}{\mu}
\left[ \begin{array}{c}
-\delta \lambda\\
\delta x\\
\end{array}\right].
\eeann
Since $\mu \neq 0$, we have
\[\left(\frac{1}{\mu},\left[ \begin{array}{c}
-\delta \lambda\\
\delta x\\
\end{array}\right]\right)\]
as the eigenpair of $\mathcal{A}$.  Since
$\sigma(\mathcal{A}^{'})=\sigma(\mathcal{A})$, the eigenvalue
$\frac{1}{\mu}$ is also $\sigma(\mathcal{A})$.  It follows that the complex
conjugates $\bar{\mu}$ and $\frac{1}{\bar{\mu}} \in \sigma(\mathcal{A})$ as
the characteristic polynomial $p(\mu)=p(\bar{\mu})$ implies that
$p(\bar{\mu})=0$.  Thus, the eigenvalues of $\mathcal{A}$ come in reciprocal pairs and complex conjugate pairs.
\end{proof}

Since $\mathbb{H}^{F'}J\mathbb{H}^F$, then the theorem above verifies the claim of the partitioning of $2n$ eigenvalues in the
complex plane. None of the eigenvalues of $\mathbb{H}^F$ are on the unit circle since the system is linearly stabilizable and detectable.

Moreover, the eigenvectors corresponding to the stable eigenvalues sitting inside the unit circle spans subspace $E_s$.  Similarly, the subspace $E_u$ is spanned by the unstable eigenvectors whose corresponding eigenvectors lie outside the unit circle.  In Section 3.6, we show that the subspace $E_s$ is invariant.


\section{Local Stable Manifold}
\subsection{Local Stable Manifold Theorem}

\begin{theorem}
Given the dynamics
\[x_{k+1}=G(x_k)\]
\[G(0)=0\]
Let $G:U \rightarrow \mathbb{R}^n$ be a $\mathcal{C}^{r-1}(\mathbb{R}^{2n})$
map with a hyperbolic fixed point $0$.  Then there is a local stable
manifold, $W^s(0) \in \mathcal{C}^{r-1}(\mathbb{R}^{2n})$, that is
tangent to the eigenspace $E^s_0$ of the Jacobian of $G$ at $0$.  Define
\[W^s(0)=\{x \in U \vert \lim_{t\rightarrow \infty} G^t(x)=0\}\]
and \\
$E^s_0=\{$span of eigenvectors whose corresponding eigenvalues are such that
$ \vert \lambda \vert < 1\}$. \\
$W^s(0)$ is a smooth manifold.
\end{theorem}

The statement of this theorem is found in \cite{GH86}.  The following theorem was proven by Hartman \cite{Ha82}, but a more modern technique can be found in \cite{Ca81}.

Along with our dynamics (\ref{nonhamforG}) and a critical point $0$, the local
stable manifold theorem gives the existence of a local stable manifold
$W^s(0)$.  The tangency of the linear subspace $E_s$ to $W_s$ implies that the lowest degree term of the local stable manifold is a quadratic homegeneous polynomial.

\subsection{Construction of the Local Stable Manifold}

Through the Taylor approximation technique, we explicitly construct the local stable manifold term by term.  These calculations coincide with the power series solutions found in Chapter 2.  Having a system that has $n$ finite eigenvalues outside and $n$ finite eigenvalues inside the unit circle, we can find a linear transformation that block diagonalizes our system (\ref{hamfor}).  We make a linear change of coordinate such that the dynamics
(\ref{hamfor}) is block diagonalized,
\beann
\left[ \begin{array}{c}
z_s\\
z_u\\
\end{array}\right]=T
\left[ \begin{array}{c}
x\\
\lambda\\
\end{array}\right].
\eeann
Matrix $T$ transforms the dynamics (\ref{hamfor}) into
\bea \label{transham}
\left[ \begin{array}{c}
z_s^+ \\
z_u^+\\
\end{array}\right]=
\left[ \begin{array}{cc}
A_s   & 0 \\
0     & A_u\\
\end{array}\right]
\left[ \begin{array}{c}
z_s \\
z_u \\
\end{array}\right] +
\left[ \begin{array}{c}
f_s(z_s,z_u) \\
f_u(z_s,z_u) \\
\end{array}\right]
\eea
where $A_s$ is the stable matrix, $A_u$ is the unstable matrix, and
\beann
f_s(z_s, z_u)&=& f_s^{[2]}(z_s, z_u) + f_s^{[3]}(z_s, z_u) + \ldots\\
f_u(z_s, z_u)&=& f_u^{[2]}(z_s, z_u) + f_u^{[3]}(z_s, z_u) + \ldots
\eeann
where $f_s^{[d]}(z_s, z_u)$ and $f_u^{[d]}(z_s, z_u)$ are nonlinear stable
and unstable homogeneous polynomials of degree $d$.
We look for the Taylor expansion of local stable manifold $z_u=\phi(z_s) \in
\mathcal{C}^{r-1}(\mathbb{R}^n)$; thus, we seek the form
\bea \label{manexp}
\phi(z_s)=\phi^{[2]}(z_s) + \phi^{[3]}(z_s) + \ldots +
\phi^{[r-1]}(z_s).
\eea
The linear term $\phi^{[1]}(z_s)$ is zero as stated in the local stable 
manifold theorem.  The graph $z_u=\phi(z_s)$ is said to be an 
\emph{invariant manifold} if
\bea \label{invariance}
z_u^+= \phi(z_s^+) \mbox{  whenever  } z_u= \phi(z_s);
\eea
i.e., the solution of (\ref{transham}) lies in $z_u=\phi(z_s)$ for all time.
By invariance (\ref{invariance}), we have

\beann
z_u^+=\phi(z_s^+) &\Longrightarrow& A_uz_u + f_u(z_s,z_u) = \phi(A_sz_s +
f_s(z_s,z_u))\\
                  &\Longrightarrow& A_u \phi(z_u) + f_u(z_s,\phi(z_s)) =
\phi(A_sz_s + f_s(z_s,\phi(z_s))).
\eeann
Thus,
\bea \label{stabman1}
\phi(A_sz_s + f_s(z_s,\phi(z_s)))=A_u \phi(z_s) + f_u(z_s,\phi(z_s)).
\eea

Taking all the second degree terms of (\ref{stabman1}), we see
\bea \label{stabmanmap2}
A_u \phi^{[2]}(z_s) - \phi^{[2]}(A_sz_s) = - f^{[2]}_u(z_s).
\eea
The map
\bea \label{phispec}
\phi^{[2]}(z_s) \mapsto  A_u \phi^{[2]}(z_s) - \phi^{[2]}(A_sz_s)
\eea
is linear.  Now, we question the invertibility of the map.  The quadratic
term $\phi^{[2]}(z_s)$ of (\ref{stabmanmap2}) as the spectrum of the map is
\[\xi - \mu_i\mu_j \neq 0   \]
where
\[A_u v_k=\xi v_k \mbox{  and  } w_i A_s=\mu_i w_i  \]
and
\bea \label{}
\phi^{[2]}(z_s)=\sum_i^n \sum_j^n c_{s_i,s_j}v_k(w_iz_s)(w_jz_s).
\eea
The details of the proof mirror the steps in Lemma (\ref{speclem}).
Moreover, the polynomial $\phi^{[d]}(x)$ for all $2 \leq d \leq k$ also
solve the following form:
\bea \label{stabmanmapd}
A_u \phi^{[d]}(z_s) - \phi^{[d]}(A_sz_s) = - f^{[d]}_u(z_s).
\eea
Thus, $\phi^{[d]}(z_s)$ exists for $2 \leq d \leq n$ since it follows the
same arguments as above.  Hence, $\phi(z_s) $ that is
$\mathcal{C}^{r-1}(\mathbb{R}^n)$ smooth has been constructed. 

\subsection{Lagrangian Submanifold}

In consequence, if $v$ is in the stable subspace of $\mathbb{H}_{\delta,k}^{F}$, then so does $\mathbb{H}_{\delta,k}^{F}v$.  Since $(\mathbb{H}_{\delta,k}^{F})^kv \rightarrow 0$ as $k \rightarrow \infty$, in the limit $\Omega(v,w)=0$ under the dynamics (\ref{nonhampertfor}).  Thus, $\Omega(v,w)=0$ if the tangent vectors are restricted in the stable subspace of $\mathbb{H}_{\delta,k}^{F}$.

Recall the two-form $\Omega:T_{(x,\lambda)}{M} \times
T_{(x,\lambda)}\mathcal{M} \mapsto \mathbb{R}$ in $\S{4.4}$.  If $v$ is $W^s \subset \mathbb{H}_{\delta,k}^{F}$, then so does $\mathbb{H}_{\delta,k}^{F}v$.  Since $(\mathbb{H}_{\delta,k}^{F})^kv \rightarrow 0$ as $k \rightarrow \infty$, in the limit $\Omega(v,w)=0$ under the map (\ref{nonhampertfor}).
Thus, the two-form restricted to $W^s(0)$
is zero; i.e.
\bea \label{j20}
\Omega(v,w)=0\;\; \mbox{for every v and w} \in TW^s.
\eea
In addition, $W^s$ has the maximal dimension of $n$.  Hence, $W^s$ is a Lagrangian submanifold.

Recall that we represent $W^s$ as the graph $z_u=\phi(z_s)$.  The basis of $T_{z_s,\phi(z_s)}W^s$ are of the form
\beann
\frac{\partial}{\partial z_{si}}
\left[ \begin{array}{c}
z_{s1}\\
\vdots\\
z_{sn}\\
\phi_1\\
\vdots\\
\phi_n\\
\end{array}\right]\;\;\mbox{where} \in i=1,\ldots,n.
\eeann
Illuminating (\ref{j20}), we get
\bea \label{j30}
\left[ \begin{array}{c}
0\\
\vdots\\
1\\
\vdots\\
0\\
\frac{\partial \phi_1}{\partial z_{si}}\\
\vdots\\
\frac{\partial \phi_n}{\partial z_{si}}\\
\end{array}\right]'J
\left[ \begin{array}{c}
0\\
\vdots\\
\vdots\\
1\\
0\\
\frac{\partial \phi_1}{\partial z_{sj}}\\
\vdots\\
\frac{\partial \phi_n}{\partial z_{sj}}\\
\end{array}\right]
\eea
where $1$ on the first vector is placed on the \emph{ith} row and $1$ on
second vector is on the \emph{jth} row.  Multiplying (\ref{j30}) out, we
have the condition
\bea \label{closed}
\frac{\partial \phi_i}{\partial z_{sj}}-\frac{\partial \phi_j}{\partial
z_{si}}=0\;\;\mbox{  for  }i,j=1,\ldots,n.
\eea
The equation (\ref{closed}) implies that $\phi(z_s)$ is closed.  Then, by 
the
\emph{Stokes' Theorem} there exists $\psi \in \mathcal{C}^{r}(\mathbb{R}^n)$
such that
\bea \label{voila}
\phi(z_s)=\frac{\partial \psi}{\partial z_s}(z_s)\mbox{  where  }\psi(0)=0
\eea
locally on $\mathcal{N}_{\epsilon}(0)$.  Thus, we have shown that the
Lagrangian submanifold is the gradient of $\psi(z_s)$.

\section{The Optimal Cost}

From the previous section, we have found that the $z_u=\phi(z_s)$ is the 
gradient of the $\psi(z_s)$.  In this section we find that the local stable manifold is also described by 
$\lambda=\bar{\phi}(x)$, its original coordinates.  We 
start by finding the linear term.

Consider the optimal control problem of minimizing
\[ \min_u \frac{1}{2} \sum_{j=0}^{M} (x'_jQ_kx_j + 2x'_jS_ku_k + u'_kRu_k) +
x'_TPx_T\]
subject to the dynamics
\bea \label{tvdynamics}
x_{k+1}&=&A_kx_k + B_ku_k\\
x(0)&=&x_0 \nonumber
\eea
where the state vector $x\in \mathbb{R}^n$, the control $u\in \mathbb{R}^m$,
and $x'_TPx_T$ is the terminal cost.  We denote the fundamental solution matrix,
\bea \label{fundmat}
\left[\begin{array}{c}
X_k \\
\Lambda_k \\
\end{array}\right]=
\left[\begin{array}{cccc}
x^1 & x^2 & \hdots & x^n \\
\lambda^1 & \lambda^2 & \hdots & \lambda^n \\
\end{array}\right]
\eea
of the dynamics 
\bea \label{fundynamics}
\left[\begin{array}{c}
X_{k+1}\\
\Lambda_{k+1}\\
\end{array}\right]&=&\mathbb{H}_F
\left[\begin{array}{c}
X_k\\
\Lambda_k\\
\end{array}\right]\\
X(0)&=&Ix_0. \nonumber
\eea
The columns 
\beann
\left[\begin{array}{c}
x^j\\
\lambda^j\\
\end{array}\right], \; j=1,\hdots,n
\eeann
are $n$ linearly independent solutions and form a basis for the space of solutions.

Through the transformation
\[\Lambda'_k =P_k X_k,\]
the time-varying discrete Riccati equation,
\[P_{k+1}=A'_k P_{k+1}A_k + (A'_k P_{k+1}B_k + S'_k)K_k + Q_k,\]
and the time-varying discrete control feedback,
\[K_k=-(B'_kP_{k+1}B_k + R_k)^{-1}(A'_kP_{k+1}B_k + S_k)',\]
are derived as shown in the proof below.

\begin{theorem} \label{z2x1}
Suppose ($X_k,\Lambda_k$) is solution to (\ref{fundynamics}) and $X_k$ 
is invertible, then $P_k=\Lambda'_k X^{-1}_k $ satisfies the time-varying 
Riccati equation.
\end{theorem}

\begin{proof}
From the Pontryagin Maximum Principle, we have
\bea \label{dHdx}
U_k=-R^{-1}_k(B'_k \Lambda'_{k+1} + S'_k X_k ).
\eea
By letting $\Lambda'_k =P_k X_k $ and substituting the dynamics 
(\ref{tvdynamics}),
\bea \label{dHdx2}
U_k=-R^{-1}_kB'_k P_{k+1}(A_k X_k + B_k U_k)-R^{-1}_kS'_kX_k
\eea
which reduces to
\bea \label{dHdx3}
(I + R^{-1}_kB'_k P_{k+1}B)U_k= -R^{-1}_kB'_k P_{k+1}A_k X_k - 
R^{-1}_kS'_kX_k.
\eea
Multiplying (\ref{dHdx3}) with $R_k$ becomes
\bea \label{findHdx}
U_k=K_kX_k.
\eea
where
\[K_k=-(B'_k P_{k+1}B_k + R_k)^{-1}(B'_k P_{k+1}A_k + S'_k),\]
the time-varying control feedback.  We recall the costate equation for 
$\Lambda$ in the Pontryagin Maximum Principle,
\bea \label{costate1}
\Lambda'_k=A_k'\Lambda'_{k+1} + Q_k X_k + S'_kU_k
\eea
Substituting the dynamics (\ref{tvdynamics}) and $U_k=K_k X_k$, the costate 
equation becomes
\bea \label{costate2}
\Lambda'_k=A'_k P_{k+1}A_k X_k + (A'_k P_{k+1}B_k + S'_k)K_k X_k + Q_kX_k
\eea
Since the LHS of (\ref{costate2}) is $\Lambda_k =P_kX_k$, we obtain the 
time-varying Riccati equation as $X_k$ is cancelled out on both sides,
\bea \label{tvRic}
P_k=A'_k P_{k+1}A_k + (A'_k P_{k+1}B_k + S'_k)K_k + Q_k.
\eea
Thus, $\Lambda'_k =P_kX_k$ solves the time-varying Riccati equation once 
$(X_k,\Lambda_k)$ is found.
\end{proof}

It follows that the DTARE (\ref{ric1}) can be solved.  By letting $k \rightarrow -\infty$ in the dynamics (\ref{fundynamics}), each linearly independent column vectors of $(X_k, \Lambda_k)'$ converges to the stable direction and thus, forming the basis for the stable subspace.  If $X^{-1}$ exists, then fundamental matrix solution has the same span as
\beann
\left[\begin{array}{c}
I \\
P_k \\
\end{array}\right].
\eeann
Hence, $P_k \rightarrow P$ as $k \rightarrow - \infty$.  Therefore, we solve the DTARE,
\[P=A'PA + (A'PB + S')K + Q\]
as $k \rightarrow -\infty$.  Moreover, the stable linear subspace is 
\bea \label{lininvariant}
E_s=\left[\begin{array}{c}
I \\
P \\
\end{array}\right].
\eea

Thus, $\lambda=Px$ is the solution of the linear Hamiltonian equations.  The linear term solution of the nonlinear Hamiltonian system is also $\lambda=Px$. We now show that there exists $\lambda=\bar{\phi}(x)$ that describes the local stable manifold of the dynamics which has $\lambda=Px$ is its linear part.  

Define the mapping $F:\mathcal{N}_{\epsilon}(0,0) \subset \mathbb{R}^{2n} \longrightarrow \mathbb{R}^n$ by
\[F(x,\lambda)=\lambda - Px + g(x,\lambda)\]
where $g(x,\lambda)$ contains only the higher order terms.  Then, for 
every $(x,\lambda) \in \mathcal{N}_{\epsilon}(0,0)$ $F(x,\lambda)=0$.  Also the differential operator $\frac{\partial F}{\partial \lambda}(0)=I$ is invertible since $\frac{\partial g}{\partial \lambda}(0)=0$ because $g$ only contains the nonlinearities.  We then invoke the IFT.  It follows that there exists $\lambda=\bar{\phi}(x)$ where $\bar{\phi}(x)\in U$ where $U$ is some subset of $\mathbb{R}^n$ such that $F(x,\bar{\phi}(x))=0$ for $(x,\bar{\phi}(x)) \in \mathcal{N}_{\epsilon}(0,0)$.
Thus, $\lambda=\bar{\phi}(x)$ is the graph of the local stable manifold in the $x, ~\lambda$ coordinates.

\section{Proof of the Main Result (Theorem \ref{thebigkahuna})}

\begin{proof}
Recall the nonlinear Hamiltonian system (\ref{nonham}).  Assuming invertibility of $H_{\lambda,x}$ at around 0 gives the nonlinear Hamiltonian system (\ref{nonhamforG}) that propagates in the direction of increasing time.  When the nonlinear Hamiltonian sytem is linearized around the critical point zero, the system resembles that of the linear Hamiltonian system (\ref{hamfor}).  We know that the linear Hamiltonian system is hyperbolic from (\ref{symplectic}) and (\ref{eigenthm}).  We then invoke the Local Stable Manifold Theorem as we have the dynamics (\ref{nonhamforG}) with a hyperbolic point $0$.  Thus, there exists a local stable manifold $W_s$ such that the linear invariant subspace (\ref{lininvariant}) is tangent at $0$.  Since the linear and nonlinear Hamiltonian structure preserve symplectic form (\ref{symplectic}) and (\ref{symplectic3}), it follows that if $\Omega$ is restricted on $W^s$ then $\Omega$ is zero (\ref{j20}); i.e., for $u, v$ are the basis in $ T_pW^s_{(0)}$ then $\Omega(u,v)=0$ as $k \rightarrow \infty$.  Then, the manifold is a Lagrangian submanifold.  The Lagrangian submanifold described by
\[z_u=\phi(z_s)\]
satisfies (\ref{closed}). By Stokes' Theorem there exists $\psi \in \mathcal{C}^{r}(\mathbb{R}^n)$
such that
\beann 
\phi(z_s)=\frac{\partial \psi}{\partial z_s}(z_s)\mbox{  where  }\psi(0)=0.
\eeann
We have also shown that $\lambda=\bar{\phi}(x)$ describes the local stable manifold by invoking IFT.
Moreover, the linear term  of $\lambda=\bar{\phi}(x)$ is $\lambda=Px$ which solves DTARE (\ref{ric1}).  It follows that 
\bea \label{upshot}
\lambda=\bar{\phi}(x)=\frac{\partial \pi}{\partial x}(x).
\eea
Thus, the graph of the gradient of the optimal cost is the local stable manifold of the associated Hamiltonian dynamics (\ref{nonhamforG}).  

Lastly, we show that $\pi(x)$ and $\kappa(x)$ solves the DPE (\ref{dpe1}), (\ref{dpe2}).  From the Pontryagin Maximum Principle and (\ref{upshot}), we have
\[u^*=\kappa(x)=\mbox{arg}\min_vH(x,\frac{\partial \pi}{\partial x}(x),v).\]
Equivalently,
\[\frac{\partial H}{\partial u}(x,\frac{\partial \pi}{\partial x}(x),\kappa(x))=0. \]
Then, $\frac{\partial \pi}{\partial x}(x)$ and $\kappa(x)$ solve DPE2.  Since 
\bea
\lambda=\frac{\partial H}{\partial x},
\eea
and (\ref{upshot}), we have 
\bea \label{predpe1}
\frac{\partial \pi}{\partial x}= \frac{\partial \pi}{\partial x}(x)\frac{\partial f}{\partial x}(x,u^*)+ \frac{\partial l}{\partial x}(x,u^*)
\eea
Integrating (\ref{predpe1}) w.r.t. $x$, we get
\[\pi(x)-\pi(f(x,\kappa(x)))-l(x,\kappa(x))=0\]
which DPE1 (\ref{dpe1}).  Thus, $\pi$ and $\kappa$ solve the DPE (\ref{dpe1}), (\ref{dpe2}).  Furthermore $\pi \in \mathcal{C}^r$ and $\kappa \in \mathcal{C}^{r-1}$ since $l \in \mathcal{C}^r$ and $f \in \mathcal{C}^{r-1}$ in the Hamiltonian.
\end{proof}

%% file: chapter5_5.tex
In this chapter, we look at the case where the linear map (\ref{hamsys}) is not a diffeomorphism.  
This is the case where zero is a closed loop eigenvalue and therefore the Hamiltonian matrix is not invertible.  
As in the previous case, we study the eigenstructure and symplectic properties of the mixed direction nonlinear 
Hamiltonian dynamics.  
Ultimately, we generalize the Local Stable Manifold Theorem for a bidirectional discrete map with a hyperbolic fixed point.

\section{Discrete-Time Version of Gronwall's Inequalities}

We begin with a discrete-time version of Gronwall's inequality.
These lemmas will be useful in the proof of the local existence of a stable manifold.

\begin{lemma} (Finite Difference Form) \label{gronwall1}
Suppose the sequence of scalars $\{u_j \}_{j=0}^{\infty}$ satisfies the difference inequality
\bea \label {groneqn}
u_{k+1} \leq  \delta u_k + L,
\eea
then
\[u_k \leq \delta^k u_0 + L\sum_{j=0}^{k-1}\delta^{k-1-j}.\]
\end{lemma}
The proof of the lemma above clearly follows from summing the equation (\ref{groneqn}) $k$-times. 

\begin{lemma} (Summation Form) \label{gronwall2}
Suppose $\{\xi\}_{j=0}^{\infty}$ is a sequence that satisfies
\[\vert \xi_k \vert \leq C_1 \sum_{j=0}^{k-1} \vert \xi_j \vert + C_2\]
and constants $C_1, C_2 \geq 0$, then 
\[\vert \xi_k \vert \leq C_2 \sum_{j=1}^{k}(1+C_1)^{j}\]
\end{lemma}

\begin{proof}
Let $s_k=\sum_{j=0}^{k-1}\vt \xi_j \vt$.  Then, the sequence $\{s\}_{j=0}^{\infty}$ satisfies
\[s_{k+1} \leq (1 + C_1) s_k + C_2\]
where $C_1,C_2 \geq 0$.  By the discrete time form of Gronwall's inequality,
\[\vert s_k \vert \leq (1+C_1)^k \vert s_0 \vert + C_2 \sum_{j=0}^{k-1}(1+C_1)^{k-1-j}.\]
It follows that
\beann
\vert \xi_k \vert &\leq& (1 + C_1) \vert s_k \vert + C_2\\
                  &\leq& (1 + C_1) \Big[(1+C_1)^k \vert s_0 \vert + C_2 \sum_{j=0}^{k-1}(1+C_1)^{j}\Big] \\
                  &\leq& C_2 \sum_{j=0}^{k-1}(1+C_1)^{k-j}\\
                  &\leq& C_2 \sum_{j=1}^{k}(1+C_1)^j\\
\eeann
since $\vert s_0\vert=0$. 
\end{proof}

\section{Local Stable Manifold for the Bidirectional Discrete-Time Dynamics}

In this section we prove the existence of a local stable manifold $\lambda=\phi(x)$ for the Hamiltonian dynamics,
\bea \label{nondiffham}
\left[ \begin{array}{c}
x^+ \\
\lambda \\
\end{array}\right]&=&
\left[\begin{array}{cc}
A & -BR^{-1}B'\\
Q & A'\\
\end{array}\right]
\left[\begin{array}{c}
x \\
\lambda^+ \\
\end{array}\right]+
\left[\begin{array}{c}
F(x,\lambda^+)\\
G(x,\lambda^+)\\
\end{array}\right]
\eea
where $x,\lambda \in \mathbb{R}^n$ and zero is an eigenvalue of A.  The nonlinear terms, $F$ and $G$, are $C^k$ functions for $k \geq 1$ such that
\bea \label{condsnonl}
F(0,0)=0, && G(0,0)=0 \\
\frac{\partial F}{\partial (x,\lambda)}(0,0)=0, && \frac{\partial G}{\partial (x,\lambda)}(0,0)=0. \nonumber
\eea

The proof of the existence of a local stable manifold requires the discussion on the stability of the nonlinear state dynamics.  The next subsection deals with the local asymptotic stability of the state dynamics. Consequently, we describe the diagonalization of the bidirectional Hamilton system.  Finally, we show the existence of $\lambda=\phi(x)$.

\subsection{Preliminaries}

First, we introduce a $C^{\infty}$ cut-off function $\rho(y): \mathbb{R}^n \longrightarrow [0,1]$ such that
\beann
\rho(y) =
   \begin{cases}
    1, &\text{if $0 \leq \vert y \vert \leq 1$}\\
    0, &\text{if $\vert y \vert > 2$}
   \end{cases}
\eeann
and $0 \leq \rho(y) \leq 1$ otherwise.  Then we define the functions
\bea \label{locfuncs}
F(x,\lambda^+;\epsilon)&:=& F(x\rho(\frac{x}{\epsilon}),\lambda^+\rho(\frac{\lambda^+}{\epsilon}))\\
G(x,\lambda^+;\epsilon)&:=& G(x\rho(\frac{x}{\epsilon}),\lambda^+\rho(\frac{\lambda^+}{\epsilon}))\nn
\eea
for $x,\lambda^+ \in \mathbb{R}^n$.  Since the $F(x,\lambda^+)$ and $G(x,\lambda^+)$ agree with $F(x,\lambda^+;\epsilon)$ and $G(x,\lambda^+;\epsilon)$, respectively, for $\vt x \vt, \vt \lambda^+ \vt \leq \epsilon$, it suffices to prove the existence of a stable manifold for some $\epsilon > 0$. 

Now, we show that there exists $N_1>0$ and $N_2>0$ such that
\bea
\vt F(x,\lambda;\epsilon)-F(\tilde{x},\tilde{\lambda};\epsilon) \vt &\leq& N_1 \epsilon \Big[ \vert x- \tilde{x} \vert + \vt  \lambda -\tilde{\lambda} \vt \Big]\label{fbnd}\\
\vt G(x,\lambda;\epsilon)-G(\tilde{x},\tilde{\lambda};\epsilon) \vt &\leq& N_1 \epsilon \Big[ \vert x- \tilde{x} \vert + \vt  \lambda -\tilde{\lambda} \vt \Big]\label{gbnd}
\eea
and 
\bea
\Big \vt \frac{\partial F}{\partial (x,\lambda)}(x,\lambda;\epsilon)- \frac{\partial F}{\partial (x,\lambda)} (\tilde{x},\tilde{\lambda};\epsilon) \Big \vt &\leq&  N_2 \Big[ \vert x- \tilde{x} \vert + \vt  \lambda -\tilde{\lambda} \vt \Big] \label{fbnd}\\
\Big \vt \frac{\partial G}{\partial (x,\lambda)}(x,\lambda;\epsilon)-\frac{\partial G}{\partial (x,\lambda)}(\tilde{x},\tilde{\lambda};\epsilon) \Big \vt &\leq& N_2 \Big[ \vert x- \tilde{x} \vert + \vt  \lambda -\tilde{\lambda} \vt \Big]. \label{gbnd}
\eea
Since $\rho(y)$ and its partial derivatives are continuous functions with compact support there exists $M>0$ such that
\beann
\Big \vert \frac{\partial \rho}{\partial y}(y) \Big \vert &\leq& M\\
\Big \vert \frac{\partial^2 \rho}{\partial y^2}(y) \Big \vert &\leq& M
\eeann
for all $\lambda \in \mathbb{R}^n$.  We also choose $M>0$ large enough that
\bea 
\Big \vt \frac{\partial F}{\partial x}(x,\lambda)         \Big \vt &\leq& M \vt x \vt \label{bndonf1a}\\
\Big \vt \frac{\partial F}{\partial \lambda}(x,\lambda)   \Big \vt &\leq& M \vt \lambda \vt \label{bndonf1b}\\
\Big \vt \frac{\partial^2 F}{\partial x^i \partial \lambda^j}(x,\lambda) \Big \vt &\leq& M,~~~i,j=1,2 \label{bndonf2}
\eea
because of the condition (\ref{condsnonl}) for $\vt x \vt, \vt \lambda \vt <1$.
By the Mean Value Theorem,
\bea \label{mvt}
\vert F(x,\lambda;\epsilon) - F(\tilde{x},\tilde{\lambda};\epsilon)  \vert &\leq & \vt F(x,\lambda;\epsilon) - F(\tilde{x},\lambda;\epsilon) + F(\tilde{x},\lambda;\epsilon) - F(\tilde{x},\tilde{\lambda};\epsilon) \vt \nn \\
                                                                           &\leq & \vt F(x,\lambda;\epsilon) - F(\tilde{x},\lambda;\epsilon) \vt + \vt F(\tilde{x},\lambda;\epsilon) - F(\tilde{x},\tilde{\lambda};\epsilon)\vt \nn \\
                                                                           &\leq & \Big \vt \frac{\partial F}{\partial x}(\xi_1,\lambda;\epsilon)  \Big \vt \vt x-\tilde{x} \vt +\Big \vt \frac{\partial F}{\partial \lambda}(x,\xi_2;\epsilon)  \Big \vt \vt \lambda-\tilde{\lambda} \vt \nn 
\eea
where $\xi_1$ is between $x$ and $\tilde{x}$ and $\xi_2$ is between $\lambda$ and $\tilde{\lambda}$.  Similarly,
\beann
\Big \vert \frac{\partial F}{\partial x}(x,\lambda;\epsilon) - \frac{\partial F}{\partial x}(\tilde{x},\tilde{\lambda};\epsilon)  \Big \vert &\leq & \Big \vt  \frac{\partial^2 F}{\partial x^2} (\xi_1,\lambda)\Big \vt  \vt x-\tilde{x} \vt  + \Big \vt \frac{\partial^2 F}{\partial x \partial \lambda} (x,\xi_2) \Big \vt  \vt \lambda -\tilde{\lambda}  \vt  \\
\Big \vert \frac{\partial F}{\partial \lambda}(x,\lambda;\epsilon) - \frac{\partial F}{\partial \lambda}(\tilde{x},\tilde{\lambda};\epsilon)  \Big \vert &\leq & \Big \vt  \frac{\partial^2 F}{\partial \lambda \partial x} (\xi_1,\lambda)\Big \vt  \vt x-\tilde{x} \vt  + \Big \vt \frac{\partial^2 F}{\partial \lambda^2} (x,\xi_2) \Big \vt  \vt \lambda -\tilde{\lambda}  \vt.
\eeann

Next we estimate for $0 \leq \epsilon < \frac{1}{2}$
\beann
\Big \vt \frac{\partial F}{\partial x}(\xi_1,\xi_2;\epsilon) \Big\vt &=& \Big\vt \frac{\partial F}{\partial x}(\rho(\frac{\xi_1}{\epsilon})\xi_1,\frac{\xi_2}{\epsilon})\xi_2)\Big\vt \; \Big\vt \frac{\partial \rho}{\partial y}(\frac{\xi_1}{\epsilon})\frac{\xi_1}{\epsilon} + \rho(\frac{\xi_1}{\epsilon}) \Big\vt \nn \\
&\leq& M \vt \rho(\frac{\xi_1}{\epsilon})\vt \vt \xi_1 \vt \Big( \Big \vt \frac{\partial \rho}{\partial y}(\frac{\xi_1}{\epsilon}) \frac{\xi_1}{\epsilon} \Big \vt + \vt \rho(\frac{\xi_1}{\epsilon}) \vt \Big) \nn \\ 
&\leq& M[M+1]\epsilon
\eeann
and 
\beann
\Big \vt \frac{\partial^2 F}{\partial x^2}(\xi_1,\xi_2;\epsilon) \Big\vt &=& \Big \vt \frac{\partial^2 F}{\partial x^2}(\rho(\frac{\xi_1}{\epsilon})\xi_1,\frac{\xi_2}{\epsilon})\xi_2) \Big \vt \; \Big\vt \frac{\partial \rho}{\partial y}(\frac{\xi_1}{\epsilon})\frac{\xi_1}{\epsilon} + \rho(\frac{\xi_1}{\epsilon}) \Big\vt \\
&&+ \Big\vt \frac{\partial F}{\partial x}(\rho(\frac{\xi_1}{\epsilon})\xi_1,\frac{\xi_2}{\epsilon})\xi_2)\Big\vt \; \Big\vt \frac{\partial \rho}{\partial \xi_1}(\frac{\xi_1}{\epsilon})\frac{2}{\epsilon} + \frac{\partial^2 \rho}{\partial y^2}(\frac{\xi_1}{\epsilon}) \frac{\xi_1}{\epsilon^2} \Big\vt \\
& \leq & M(1+2M)+8M^2
\eeann
for $\vt \xi_1 \vt, \vt \xi_2 \vt \leq \epsilon$.

Let
\beann
N_1 &=& M^2[M+1]\\
N_2 &=& M(1+2M)+8M^2.
\eeann
It follows that
\bea
\Big \vert \frac{\partial F}{\partial (x,\lambda)}(\xi_1,\xi_2;\epsilon) \Big \vert &\leq& N_1\epsilon \label{origmvt1} \\
\Big \vt \frac{\partial^2 F}{\partial x^i \partial \lambda^j}(\xi_1,\xi_2;\epsilon) \Big\vt &\leq& N_2,~~~i,j=1,2 \label{origmvt2}
\eea
for $\vt \xi_1 \vt, \vt \xi_2 \vt < \epsilon$.  The inequalities above also hold $G$ as well.

Thus, 
\bea
\vt F(x,\lambda;\epsilon)-F(\tilde{x},\tilde{\lambda};\epsilon) \vt &\leq& N_1 \epsilon \Big[ \vert x- \tilde{x} \vert + \vt  \lambda -\tilde{\lambda} \vt \Big]\label{fbnd1}\\
\vt G(x,\lambda;\epsilon)-G(\tilde{x},\tilde{\lambda};\epsilon) \vt &\leq& N_1 \epsilon \Big[ \vert x- \tilde{x} \vert + \vt  \lambda -\tilde{\lambda} \vt \Big]\label{gbnd1}
\eea
and 
\bea
\Big \vt \frac{\partial F}{\partial (x,\lambda)}(x,\lambda;\epsilon)- \frac{\partial F}{\partial (x,\lambda)} (\tilde{x},\tilde{\lambda};\epsilon) \Big \vt &\leq&  N_2 \Big[ \vert x- \tilde{x} \vert + \vt  \lambda -\tilde{\lambda} \vt \Big] \label{fbnd2}\\
\Big \vt \frac{\partial G}{\partial (x,\lambda)}(x,\lambda;\epsilon)-\frac{\partial G}{\partial (x,\lambda)}(\tilde{x},\tilde{\lambda};\epsilon) \Big \vt &\leq& N_2 \Big[ \vert x- \tilde{x} \vert + \vt  \lambda -\tilde{\lambda} \vt \Big].
 \label{gbnd2}
\eea

Henceforth we suppress $\epsilon$ and write $F(x,\lambda), G(x,\lambda)$ for $F(x,\lambda;\epsilon),~G(x,\lambda;\epsilon)$.

\subsection{Stability of the Nonlinear Dynamics}

The stable invariant manifold is described by $\lambda=Px$ for the linear bidirectional Hamiltonian dynamics where $P$ is the solution to the discrete algebraic Riccati equation.  Then, we know from Chapter 3 that the linear term of the stable manifold for the nonlinear bidirectional Hamiltonian dynamics, $\lambda=\phi(x)$, is $Px$.  Thus, the local stable manifold is of the form
\bea \label{phiform}
\lambda = \phi(x)=Px + \psi(x)
\eea
where $\psi(x)$ contains all the nonlinear terms.

Suppose we substitute (\ref{phiform}) into the state dynamics in (\ref{nondiffham}), then the nonlinear state dynamics becomes
\beann
(I+BR^{-1}B'P)x^+=Ax-BR^{-1}B'\psi(x^+)+ F(x,Px^+ + \psi(x^+)).
\eeann
By the Matrix Inversion Lemma (\cite{LS95}), we have that 
\[(I+BR^{-1}B'P)^{-1}=(I-B(B'PB+R)^{-1}B'P).\]
Then, it follows that 
\bea \label{nlinearxmap}
x^+&=&(A+BK)x + f_{\psi}(x,x^+)\\
x(0)&=&x_0 \nn
\eea
where $K=-(B'PB + R)^{-1}B'P$ and $f_{\psi}(x,x^+)=(I+BR^{-1}B'P)^{-1}(F(x,Px^+ + \psi(x^+)) - BR^{-1}B'\psi(x^+))$.  

The implicit equation above can be solved.  Let $\mathcal{F}:\mathcal{N}_{\epsilon}(0) \subset \mathbb{R}^{2n} \rightarrow \mathbb{R}^n$ such that 
\[\mathcal{F}(x,x^+)=x^+ - (A+BK)x-f_{\psi}(x,x^+)=0\]
for $x,x^+ \in \mathcal{N}_{\epsilon}(0)$.  Then, for $0 \in \mathcal{N}_{\epsilon}(0)$ the Jacobian
\beann
\frac{\partial \mathcal{F}}{\partial x^+}(0) &=& I-\frac{\partial f_{\psi}}{\partial x^+}(0) \\
                                             &=& I-(I+BR^{-1}B'P)^{-1}\Big(\frac{\partial F}{\partial \lambda^+}(0,\phi(0))\frac{\partial \psi}{\partial x^+}(0) - BR^{-1}B'\frac{\partial \psi}{\partial x^+}(0)\Big)\\
                                             &=& I,
\eeann
because of the condition (\ref{condsnonl}) and $\psi(x^+)$ only contains nonlinear terms.  Then, by the Implicit Function Theorem there exists $\mathbb{F}(x)$ such that
\bea \label{nlinearxmap2}
x^+&=&\mathbb{F}(x)
\eea
is equivalent to the earlier state dynamics (\ref{nlinearxmap}).  Moreover,  the linear term of $\mathbb{F}(x)$ is $(A+BK)x$, i.e.;
\bea \label{consxeq2}
\mathbb{F}(x)= (A+BK)x + F_{\psi}(x)
\eea
because
\beann
\frac{\partial \mathbb{F}}{\partial x}(0)&=&- \Big(\frac{\partial \mathcal{F}}{\partial x^+}(0)\Big)^{-1}\frac{\partial \mathcal{F}}{\partial x}(0)\\
&=&-I[-(A+BK)]\\ 
&=& A+BK. 
\eeann
It follows that $F_{\psi}(x)$ contains only the nonlinear terms and thus,
\bea \label{consxeq3}
F_{\psi}(0)=0
\eea
and
\[\frac{\partial F_{\psi}}{\partial x_i}(0)=0\;\;i=1,\ldots,n.\]

The linear part of (\ref{nlinearxmap}) is
\bea \label{linearxmap}
x^+=(A+BK)x.
\eea
Since the eigenvalues of $(A+BK)$ lie strictly inside the unit circle, the term $(A+BK)^k x_0 \longrightarrow 0$ as $k \longrightarrow \infty$. Thus, the system (\ref{linearxmap}) is asymptotically stable.  Also, it implies that there exists a unique positive definite $P$ that satisfies the Lyapunov equation
\[(A+BK)'P(A+BK)-P=-I.\]

Now we show the stability of the nonlinear dynamics
\beann
x_k &=& (A+BK)x + F_{\psi}(x)\\
x(0)&=& 0. 
\eeann

We must prove that 
\beann
\lim_{x \rightarrow 0} \frac{\vt F_{\psi}(x) \vt}{\vt x \vt}=0;
\eeann
i.e., given any $\varepsilon > 0$ and any $\psi(x)$ satisfying the conditions
\bea 
\psi(0)&=&0 \label{cond1}\\
\vt \psi(x)-\psi(\bar{x})  \vt &\leq& l(\epsilon)\vt x - \bar{x}  \vt \label{cond2}
\eea
where $l(\epsilon) \longrightarrow 0$ as $\epsilon \rightarrow 0$,  there exists $\delta >0$ such that
\[ \frac{\vt F_{\psi}(x) \vt}{\vt x \vt} < \varepsilon ~~~~~\text{whenever}~~ \vt x \vt < \delta.\]
We define $\psi(x)$ to be the nonlinear term of the stable manifold in (\ref{}).  The conditions (\ref{cond1}) and (\ref{cond2}) will be necessary for the proof of the local stable manifold theorem.

Recall that
\beann
x^+=\mathbb{F}(x)= (A+BK)x + F_{\psi}(x). 
\eeann
Then,
\beann
0=\mathcal{F}(x,x^+)&=& x^+ - (A+BK)x-f_{\psi}(x,x^+)\\
                    &=&(A+BK)x + F_{\psi}(x)- (A+BK)x-f_{\psi}(x,(A+BK)x + F_{\psi}(x))\\
                    &=&F_{\psi}(x)-f_{\psi}(x,(A+BK)x + F_{\psi}(x)).
\eeann
It follows that
\bea \label{Fpsi}
F_{\psi}(x)&=&(I+BR^{-1}B'P)^{-1}\Big[F(x,P((A+BK)x + F_{\psi}(x)) + \psi((A+BK)x + F_{\psi}(x))) \nn \\
           &&~~- BR^{-1}B'\psi((A+BK)x + F_{\psi}(x))\Big].
\eea
Let $\mathbb{B}_1= \V (I+BR^{-1}B'P)^{-1} \V$, $\mathbb{B}_2= \V BR^{-1}B' \V$, $\mathbb{P}=\V P \V$ and $\alpha = \max_{i} \vt \lambda_i \vt$ where $\lambda_i \in \sigma(A+BK)$ and $\vt \lambda_i \vt < 1$.   We have the following from (\ref{Fpsi}),
\beann
\vt F_{\psi}(x) \vt &\leq& \mathbb{B}_1 N_1 \epsilon \left[\vt x \vt + \vt P(A+BK)x + PF_{\psi}(x) + \psi((A+BK)x + F_{\psi}(x)))\vt \right] \\
&&~~+ \mathbb{B}_1\mathbb{B}_2 \vt \psi((A+BK)x + F_{\psi}(x))) \vt.
\eeann
because of (\ref{fbnd1}).
Using the Lipschitz condition (\ref{cond2}) for $\psi(x)$,
\beann
\vt F_{\psi}(x) \vt &\leq& \left[ \mathbb{B}_1 N_1 \epsilon + \alpha(\mathbb{B}_1N_1\epsilon \V P \V + \mathbb{B}_1N_1 \epsilon l(\epsilon) + \mathbb{B}_1\mathbb{B}_2l(\epsilon)) \right] \vt x \vt \\
&&+ \left[ \mathbb{B}_1N_1\epsilon\V P \V + \mathbb{B}_1N_1 \epsilon l(\epsilon) + \mathbb{B}_1\mathbb{B}_2 l(\epsilon)\right] \vt F_{\psi}(x) \vt.
\eeann
Solving for $\vt F_{\psi}(x) \vt$,
\beann
\vt F_{\psi}(x) \vt &\leq& \frac{\mathbb{B}_1 N_1 \epsilon + \alpha(\mathbb{B}_1N_1\epsilon \V P \V + \mathbb{B}_1N_1 \epsilon l(\epsilon) + \mathbb{B}_2l(\epsilon))}{1-\left( \mathbb{B}_1N_1\epsilon + \mathbb{B}_1N_1 \epsilon l(\epsilon) + \mathbb{B}_2 l(\epsilon)\right)}.
\eeann
Let $\delta= \frac{1-\left( \mathbb{B}_1N_1\epsilon + \mathbb{B}_1N_1 \epsilon l(\epsilon) + \mathbb{B}_2 l(\epsilon)\right)}{\mathbb{B}_1 N_1 \epsilon + \alpha(\mathbb{B}_1N_1\epsilon \V P \V + \mathbb{B}_1N_1 \epsilon l(\epsilon) + \mathbb{B}_2l(\epsilon))}~\varepsilon$.  For some $\epsilon >0$ and $\varepsilon >0$, we have that $\delta > 0$.
Then,
\beann
\vt F_{\phi}(x) \vt &\leq& \frac{\mathbb{B}_1 N_1 \epsilon + \alpha(\mathbb{B}_1N_1\epsilon \V P \V + \mathbb{B}_1N_1 \epsilon l(\epsilon) + \mathbb{B}_2l(\epsilon))}{1-\left( \mathbb{B}_1N_1\epsilon + \mathbb{B}_1N_1 \epsilon l(\epsilon) + \mathbb{B}_2 l(\epsilon)\right)}~ \vt x \vt\\
                          &\leq& \frac{\mathbb{B}_1 N_1 \epsilon + \alpha(\mathbb{B}_1N_1\epsilon \V P \V + \mathbb{B}_1N_1 \epsilon l(\epsilon) + \mathbb{B}_2l(\epsilon))}{1-\left( \mathbb{B}_1N_1\epsilon + \mathbb{B}_1N_1 \epsilon l(\epsilon) + \mathbb{B}_2 l(\epsilon)\right)}~ \delta \\
                          &\leq& \varepsilon.
\eeann
Thus, 
\[F_{\phi}(x)= o(\vt x \vt).\]
Now, we use the Lyapunov argument.
Let $v(x)=x'Px$.  Then,
\beann
\Delta v(x)&=& v(x^+)-v(x)\\
           &=& {x^+}' P x^+ -x' P x\\
           &=& [(A+BK)x-F_{\psi}(x)]'P [(A+BK)x-F_{\psi}(x)]-x'Px\\
           &=& x'((A+BK)'P(A+BK)-P)x + 2x'(A+BK)'PF_{\psi}(x)\\
           &=& -\vt x \vt^2 + 2x'(A+BK)'PF_{\psi}(x).
\eeann
since
\beann
\vt F_{\psi}(x) \vt \leq \frac{1}{3p}\vt x \vt
\eeann
and
\[\vt 2x'(A+BK)'P F_{\psi}(x)\vt \leq \frac{2}{3}\vt x \vt^2  \]
for some $p>0$.
Thus,
\beann
\Delta v(x)= -\frac{\vt x \vt^2}{3}<0.
\eeann
Therefore, the nonlinear dynamics is locally asymptotically stable uniform for all $\psi \in \mathbb{X}$.

\subsection{Diagonalization of the Hamiltonian Matrix}

Recall the bidirectional nonlinear dynamics in (\ref{nondiffham})
\beann 
\left[ \begin{array}{c}
x^+ \\
\lambda \\
\end{array}\right]&=&
\left[\begin{array}{cc}
A & -BR^{-1}B'\\
Q & A'\\
\end{array}\right]
\left[\begin{array}{c}
x \\
\lambda^+ \\
\end{array}\right]+
\left[\begin{array}{c}
F(x,\lambda^+)\\
G(x,\lambda^+)\\
\end{array}\right]
\eeann
where $x,\lambda \in \mathbb{R}^n$ and zero is an eigenvalue of A.  The nonlinear terms, $F$ and $G$, are $C^k$ functions for $k \geq 1$ such that
\beann
F(0,0)=0, && G(0,0)=0 \\
\frac{\partial F}{\partial (x,\lambda)}(0,0)=0, && \frac{\partial G}{\partial (x,\lambda)}(0,0)=0. \nonumber
\eeann

By substituting 
\bea \label{lasubstitute}
\lambda=Px + \psi(x)
\eea
into the state dynamics above (\ref{nondiffham}), we get
\bea \label{newstate}
x^+&=& (A+BK)x + f_{\psi}(x,x^+)
\eea
where
$f_{\psi}(x,x^+)=(I+BR^{-1}B'P)^{-1}(F(x,Px^+ + \psi(x^+)) - BR^{-1}B'\psi(x^+))$.

As we substitute (\ref{lasubstitute}) and (\ref{newstate}) into the costate dynamics in (\ref{nondiffham}), we also add $0=(BK)'\lambda^+ -(BK)'\lambda^+$.  Then, the costate dynamics becomes
\beann \label{nlinearlammap} 
\lambda=(A+BK)'\lambda^+ + \bar{Q}x + g_{\psi}(x,x^+).
\eeann
where $\bar{Q}=Q-K'B'P(A+BK)$ and $g_{\psi}(x,x^+)=  G(x,Px^+ + \psi(x^+)) + K'B'(-\psi(x^+) -Pf_{\psi}(x,x^+)$.  

Thus, the substitution of 
\[\lambda=Px + \psi(x)  \]
into the dynamics (\ref{nondiffham}) results in a new nonlinear dynamics
\bea \label{nondiffham2}
\left[ \begin{array}{c}
x^+ \\
\lambda \\
\end{array}\right]&=&
\left[\begin{array}{cc}
A+BK & 0\\
\bar{Q} & (A+BK)'\\
\end{array}\right]
\left[\begin{array}{c}
x \\
\lambda^+ \\
\end{array}\right]+
\left[\begin{array}{c}
f_{\psi}(x,x^+)\\
g_{\psi}(x,x^+)\\
\end{array}\right]
\eea
where
\beann
f_{\psi}(x,x^+)&=&(I+BR^{-1}B'P)^{-1}(F(x,\psi(x^+)) - BR^{-1}B'\psi(x^+))
\eeann
and
\beann
g_{\psi}(x,x^+)&=&G(x,\psi(x^+)) + K'B'(-\psi(x^+) -Pf_{\psi}(x,x^+).~~~~~~~~
\eeann
The nonlinear terms $f_{\psi}$ and $g_{\psi}$ are $C^k$ functions for $k \geq 1$ such that
\bea 
f_{\psi}(0,0)=0, && g_{\psi}(0,0)=0 \label{fgcond1}\\
\frac{\partial f_{\psi}}{\partial(x,x^+)}(0,0)=0, && \frac{\partial g_{\psi}}{\partial (x,x^+)}(0,0)=0.\label{fgcond2}
\eea
because of (\ref{condsnonl}), $\psi(x)$ only contains nonlinear terms and 
\[\frac{\partial \psi}{\partial x}(0)=0.\]

Now we introduce the $z$ coordinate by the transformation
\bea \label{sectrans}
\lambda=z + Sx
\eea
for some matrix $S$ to block diagonalize the block lower triangular Hamiltonian matrix in (\ref{nondiffham2}).  By substitution, the system (\ref{nondiffham2}) becomes

\beann
x^+&=& (A+BK)x + f_{\psi}(x,x^+)\\
  z&=& (A+BK)'z^+ + (A+BK)'S(A+BK)x - Sx + \bar{Q}x + h_{\psi}(x,x^+)
\eeann
where 
\[h_{\psi}(x,x^+)=(A+BK)'Sf_{\psi}(x,x^+) + g_{\psi}(x,x^+).\]
Observe from the $z$ dynamics above that the terms
\[(A+BK)'S(A+BK)x - Sx + \bar{Q}x=0\]
and recall that $\bar{Q}=Q-K'B'P(A+BK)$.
Indeed,
\bea \label{almostRiccati}
-S + A'S(A+BK)+ K'B'S(A+BK) = -Q+K'B'P(A+BK).
\eea
We know that 
\bea \label{itsRiccati}
-S + A'S(A+BK)= -Q,
\eea
is the discrete-time algebraic Riccati equation (DTARE).
Subtracting (\ref{itsRiccati}) from (\ref{almostRiccati}), we have
\beann
K'B'S(A+BK)&=&K'B'P(A+BK).
\eeann
Thus,
\bea \label{sisp}
S=P
\eea
and $S$ satisfies the DTARE.  Therefore, we have a diagonalized system

\bea \label{nondiffham3} 
\left[ \begin{array}{c}
x^+ \\
z \\
\end{array}\right]&=&
\left[\begin{array}{cc}
A+BK & 0\\
0    & (A+BK)'\\
\end{array}\right]
\left[\begin{array}{c}
x \\
z^+ \\
\end{array}\right]+
\left[\begin{array}{c}
f_{\psi}(x,x^+)\\
g_{\psi}(x,x^+)\\
\end{array}\right]
\eea
where 
\beann
f_{\psi}(x,x^+)&=&(I+BR^{-1}B'P)^{-1}(F(x,\psi(x^+)) - BR^{-1}B'\psi(x^+))
\eeann
and
\beann
h_{\psi}(x,x^+)=(A+BK)'Sf_{\psi}(x,x^+) + g_{\psi}(x,x^+).~~~~~~~~~~~~~~~
\eeann
The nonlinear terms $f_{\psi}$ and $h_{\psi}$ are $C^k$ functions for $k \geq 1$ such that
\bea 
f_{\psi}(0,0)=0, && h_{\psi}(0,0)=0 \label{fhcond1}\\
\frac{\partial f_{\psi}}{\partial(x,x^+)}(0,0)=0, && \frac{\partial h_{\psi}}{\partial (x,x^+)}(0,0)=0.\label{fhcond2}
\eea
because of (\ref{fgcond1}) and (\ref{fgcond2}).

\subsection{The Local Stable Manifold Theorem}

Given the original dynamics (\ref{nondiffham}) we look for the local stable manifold described by $\lambda=\phi(x)$.  From the Theorem (\ref{}) we have already proven that the linear term of the stable manifold for the system (\ref{nondiffham}) is $Px.$ where $P$ is the solution to DTARE.   Therefore, the local stable manifold is
\bea \label{transform1}
\lambda=Px+ \psi(x)
\eea
where $\psi(x)$ only contains the nonlinear term of $\phi(x)$.  In the two-step process of diagonalization of the system (\ref{nondiffham}), we introduce the $z$ coordinate through the transformation
\[\lambda=z+Sx.\]
Since $S=P$, it must be that
\beann
z= \lambda - Px = Px + \psi(x) -Px = \psi(x).
\eeann
Then, it suffices to prove existence of the local stable manifold $z=\psi(x)$ for the diagonalized system (\ref{nondiffham3}).  In order to show the existence of $z=\psi(x)$, we use the Contraction Mapping Principle (CMP).  To invoke the CMP, we will need a map $T:\mathbb{X}\longrightarrow \mathbb{X}$ that is a contraction on a complete metric space $\mathbb{X}$.

\begin{theorem} \label{nondiffthm}
Given the dynamics in (\ref{nondiffham3}) with the nonlinear terms $f_{\psi}$ and $h_{\psi}$ are $\mathcal{C}^k$ functions satisfying the conditions (\ref{fhcond1}) and (\ref{fhcond2}) and a hyperbolic fixed point $0 \in \mathbb{R}^{2n}$, there exists a local stable manifold $z=\psi(x)$ around the fixed point $0$ where $\psi$ is a $C^k$ function. 
\end{theorem}

\begin{proof}

First notice that $f_{\psi}$ and $g_{\psi}$ are cut-off functions.  It follows that $h_{\psi}$ is also a cut-off function.  It suffices to prove the theorem for some $\epsilon > 0$ since the cut-off functions $f_{\psi}(x,x^+;\epsilon)$ and $h_{\psi}(x,x^+;\epsilon)$ agree with $f_{\psi}(x,x^+)$ and $h_{\psi}(x,x^+)$ for $\vt x \vt,~\vt x^+ \vt \leq \epsilon$.

By (\ref{fbnd1})-(\ref{gbnd2}), (\ref{cond1})-(\ref{cond2}), and (\ref{fhcond1})-(\ref{fhcond2}), there exists $\bar{N}_1,~\bar{N}_2 >0$ such that 
\bea
\vt f_{\psi}(x,y;\epsilon)-f_{\psi}(\tilde{x},\tilde{y};\epsilon) \vt &\leq& \bar{N}_1 \epsilon \Big[ \vert x- \tilde{x} \vert + \vt  y -\tilde{y} \vt \Big] \label{fphibnd1}\\
\vt h_{\psi}(x,y;\epsilon)-h_{\psi}(\tilde{x},\tilde{y};\epsilon) \vt &\leq& \bar{N}_1 \epsilon \Big[ \vert x- \tilde{x} \vert + \vt  y -\tilde{y} \vt \Big] \label{gphibnd1}
\eea
and 
\bea
\Big \vt \frac{\partial f_{\psi}}{\partial (x,y)}(x,y;\epsilon)- \frac{\partial f_{\psi}}{\partial (x,y)} (\tilde{x},\tilde{y};\epsilon) \Big \vt &\leq&  \bar{N}_2 \Big[ \vert x- \tilde{x} \vert + \vt  y -\tilde{y} \vt \Big] \label{fphibnd2}\\
\Big \vt \frac{\partial h_{\psi}}{\partial (x,y)}(x,y;\epsilon)-\frac{\partial h_{\psi}}{\partial (x,y)}(\tilde{x},\tilde{y};\epsilon) \Big \vt &\leq& \bar{N}_2 \Big[ \vert x- \tilde{x} \vert + \vt  y -\tilde{y} \vt \Big]. \label{gphibnd2}
\eea

Moreover, from the bound (\ref{origmvt1}) we know the following
\bea \label{mvtfphi}
\Big \vt  \frac{\partial f_{\psi}}{\partial (x,y)}(\xi_1,\xi_2;\epsilon) \Big \vt \leq \bar{N}\epsilon
\eea
and
\bea \label{mvtgphi}
\Big \vt  \frac{\partial h_{\psi}}{\partial (x,y)}(\xi_1,\xi_2;\epsilon) \Big \vt \leq \bar{N} \epsilon.
\eea
for $\bar{N}>0$ and $\vt \xi_1 \vt, \vt \xi_2 \vt < \epsilon$.

Henceforth we suppress $\epsilon$ and write $f_{\psi}(x,y), h_{\psi}(x,y)$ for $f_{\psi}(x,y;\epsilon),~h_{\psi}(x,y;\epsilon)$.

Let $l(\epsilon)$ with $l(0)=0$ and $\psi \in \mathbb{X}\subset \mathbb{C}^0({\vt x \vt \leq \epsilon})$ where $\mathbb{X}$ is space of $\psi: \mathbb{R}^n \longrightarrow \mathbb{R}^n$ such that
\bea
\psi(0)&=&0 \label{Xcond1}\\
\vt \psi(x)-\psi(\bar{x}) \vt &\leq& l(\epsilon)\vt x-\bar{x}  \vt \label{Xcond2}
\eea
for $x,~\bar{x} \in \bar{\mathcal{B}}_{\epsilon}(0) \subset \mathbb{R}^n$.
We define 
\bea \label{thenorm}
\V \psi \V= \sup_{\vt x \vt \leq \epsilon}\left\vt \frac{\psi(x)}{x} \right\vt.
\eea

Since $\mathbb{X} \subset \mathbb{C}^0(\{\vt x \vt \leq \epsilon\})$,  to show $\mathbb{X}$ is a complete metric space it suffices to show that $\mathbb{X}$ is closed.  We take a sequence $\{\psi_n \} \in 
\mathbb{X}$ such that $\psi_n  \longrightarrow  \psi$ in $\mathbb{C}^0$ norm. For large $N>n$, $\vert \psi_n(x)-\psi(x) \vert 
\leq \frac{\epsilon}{2}$ for all $x \in \bar{\mathcal{B}}_{\epsilon}(0)$.  Then,
\beann
\left\vert \psi(x) - \psi(\bar{x})  \right\vert &\leq& \left\vt \psi(x) -\psi_n(x)  \right\vt + \left\vt \psi_n(x) - \psi_n(\bar{x}) \right\vt + \left\vt \psi_n(\bar{x}) - \psi(\bar{x}) \right\vt\\
& \leq & \frac{\epsilon}{2} + l(\epsilon) \vt x-\bar{x} \vt  +\frac{\epsilon}{2}.
\eeann
By letting $\epsilon \rightarrow 0$, we have that $\vert \psi(x) - \psi(\bar{x}) \vert \leq l(\epsilon) \vt x-\bar{x} \vt$.  Thus, $\psi $ is a Lipschitz function.  Similarly, the condition (\ref{Xcond1}) is easily satisfied.  It follows that
\beann
\vt \psi(0)-0 \vt &\leq& \vt \psi(0)- \psi_n(0) \vt + \vt \psi_n(0) - 0 \vt \leq \epsilon.
\eeann
Thus, $\psi \in \mathbb{X}$.  Hence $\mathbb{X}$ is closed.  
Moreover, $\mathbb{X}$ is a complete metric space with the norm defined on (\ref{thenorm}).

Solving the $z$ dynamics in (\ref{nondiffham3}) via the variation of constants formula, we have
\bea \label{nlinearlamsol}
z_j&=& (A'+K'B')^{k-j}z_k + \sum_{l=j}^{k-1}(A'+K'B')^{l-j}h_{\psi}(x_l,x_{l+1})\nn \\
\eea
for $j<k$.  Let $j=0$ and $k=\infty$, then (\ref{nlinearlamsol}) changes to
\bea \label{nlinearlamsol2}
z_0&=& \sum_{l=0}^{\infty}(A'+K'B')^{l}h_{\psi}(x_l,x_{l+1})\nn \\
\eea

We define a mapping $T: \mathbb{X} \longrightarrow \mathbb{X}$ by 
\bea \label{Tmap}
(T\psi)(x_0)= \sum_{l=0}^{\infty}(A'+K'B')^{l}h_{\psi}(x_l,x_{l+1}).
\eea
From this fixed point equation, we look for the solution
\beann
(T\psi)(x_0)=\psi(x_0).
\eeann
We must show $T\psi \in \mathbb{X}$ and prove $T$ is a contraction on $\mathbb{X}$.


Suppose $x_0=0$ is the initial condition. Clearly from the equation (\ref{nlinearxmap2})-(\ref{consxeq3}), $x_k=0$ for all $k$. Together with $x_k=0$ for all $k$ and the condition (\ref{fhcond2}) $h_{\psi}(0,0)=0$, we have 
\bea \label{sXcond1}
T\psi(0)=0.
\eea
Hence $T\psi$ satisfies the condition (\ref{Xcond1}).

We now prove the Lipschitz condition (\ref{Xcond2}) for $T\psi$.

For $\psi \in \mathbb{X}$ and the initial conditions $x_0,~\bar{x}_0 \in \mathbb{R}^n$, we denote $x_k=x(k,x_0,\psi)$
to be the solution of the state dynamics.  
\beann 
x^+&=&(A+BK)x + f_{\psi}(x,x^+)\\
x(0)&=&x_0. 
\eeann
Similarly, for $\psi \in \mathbb{X}$ and the initial conditions $\bar{x}_0 \in \mathbb{R}^n$, let $x_k=x(k,\bar{x}_0,\psi)$ tbe the solution of
\beann 
x^+&=&(A+BK)x + f_{\psi}(x,x^+)\\
x(0)&=&\bar{x}_0. 
\eeann
Recall $\alpha=\max_j \vt \lambda_j \vt$ where $\lambda_j \in \sigma(A+BK)$ and $\vt \lambda_j \vt < 1$. Using the estimate (\ref{fphibnd1}), at one-time step 
\beann
\vt x_{k+1}- \bar{x}_{k+1}\vt \leq \alpha \vt x_k - \bar{x}_k \vt + \bar{N}_1 \left[ \vt x_k-\bar{x}_k  \vt + \vt x_{k+1}- \bar{x}_{k+1}\vt   \right].
\eeann
For $1-\bar{N}_1\epsilon > 0$, 
\beann
\vt x_{k+1}- \bar{x}_{k+1}\vt \leq \frac{\alpha + \bar{N}_1\epsilon}{1-\bar{N}_1\epsilon} \vt x_k - \bar{x}_k \vt
\eeann
and recursively,
\beann
\vt x_{k}- \bar{x}_{k}\vt \leq \left( \frac{\alpha + \bar{N}_1\epsilon}{1-\bar{N}_1\epsilon} \right)^k \vt x_0 - \bar{x}_0 \vt
\eeann
As long as
\[\epsilon < \frac{1-\alpha}{2\bar{N}_1},\]
then
\[ \frac{\alpha + \bar{N}_1\epsilon}{1-\bar{N}_1\epsilon} < 1.\]
Thus,
\bea \label{uniformx1}
\vt x_k - \bar{x}_k  \vt \leq \vt x_0 - \bar{x}_0 \vt.  
\eea

Using the bounds (\ref{gphibnd1}) and (\ref{uniformx1}),
\beann
\vt T\psi(x_0)- T\psi(\bar{x}_0)\vt &\leq& \left\vt  \sum_{l=0}^{\infty}(A'+K'B')^{l}( h_{\psi}(x_l,x_{l+1})-h_{\psi}(\bar{x}_l,\bar{x}_{l+1}) )\right\vt\\
&\leq&   \sum_{l=0}^{\infty}\alpha^{l}\left\vt  h_{\psi}(x_l,x_{l+1})-h_{\psi}(\bar{x}_l,\bar{x}_{l+1}) \right\vt \\ 
&\leq&   \sum_{l=0}^{\infty}\alpha^{l} \bar{N}_1\epsilon \left[ \vt x_l - \bar{x}_l \vt + \vt x_{l+1} - \bar{x}_{l+1} \vt \right]\\
&\leq&   \sum_{l=0}^{\infty}\alpha^{l} 2\bar{N}_1\epsilon \vt x_0-\bar{x}_0 \vt \\
&\leq& \frac{2\bar{N}_1\epsilon}{1-\alpha} \vt x_0 - \bar{x}_0 \vt.
\eeann
Let
\[l(\epsilon)= \frac{2\bar{N}_1\epsilon}{1-\alpha}.\]
Notice that $l(\epsilon) \rightarrow 0$ as $\epsilon \rightarrow 0$.
Thus,
\beann
\vt T\psi(x_0)- T\psi(\bar{x}_0)\vt &\leq& l(\epsilon) \vt  x_0 - \bar{x}_0  \vt
\eeann
and so $(T\psi)$ satisfies the condition (\ref{Xcond2}) for $\epsilon >0$ sufficiently small.

Hence $T$ maps from $\mathbb{X} \rightarrow \mathbb{X}.$


Next we show T is a contraction on $\mathbb{X}$.

We express the solutions to the state dynamics
\[x_k=x(k,x_0,\psi)\]
and
\[\bar{x}_k=\bar{x}(k,\bar{x}_0,\psi)\]
in the implicit form,
\beann
x_k=(A+BK)^k x_0 + \sum_{j=0}^{k-1} (A+BK)^{k-1-j}f_{\psi}(x_j,x_{j+1})
\eeann
and
\beann
\bar{x}_k=(A+BK)^k \bar{x}_0 + \sum_{j=0}^{k-1} (A+BK)^{k-1-j}f_{\psi}(\bar{x}_j,\bar{x}_{j+1}),
\eeann
respectively.

We now denote $x_j=x(j,x_0,\psi)$ and $\bar{x}_j=\bar{x}(j,x_0,\bar{\psi})$ be the solutions to the state dynamics and satisfy the implicit form equations
\[x_k=(A+BK)^k x_0 + \sum_{j=0}^{k-1} (A+BK)^{k-1-j}f_{\psi}(x_j,x_{j+1})\]
and
\[\bar{x}_k=(A+BK)^k x_0 + \sum_{j=0}^{k-1} (A+BK)^{k-1-j}f_{\bar{\psi}}(\bar{x}_j,\bar{x}_{j+1}),\]
respectively.  

The estimates (\ref{fphibnd1})-(\ref{gphibnd1}) with the trajectories $x(j,x_0,\psi)$ and $x(j,x_0,\bar{\psi})$ becomes
\bea
\vt f_{\psi}(x,y)-f_{\psi}(\tilde{x},\tilde{y}) \vt &\leq& r_1(\epsilon)\vt y-\bar{y} \vt + r_2(\epsilon) \V \psi- \bar{\psi} \V + r_3(\epsilon)\vt x-\bar{x} \vt \label{fphibnd2}\\
\vt h_{\psi}(x,y)-h_{\psi}(\tilde{x},\tilde{y}) \vt &\leq& r_1(\epsilon)\vt y-\bar{y} \vt + r_2(\epsilon) \V \psi- \bar{\psi} \V + r_3(\epsilon)\vt x-\bar{x} \vt  \label{gphibnd2}
\eea 
where
\beann
r_1(\epsilon)&=& n_{1,1}l(\epsilon) + n_{1,2}\epsilon + n_{1,3}l(\epsilon)\epsilon\\
r_2(\epsilon)&=& n_{2,1}\epsilon + n_{2,2}\epsilon^2\\
r_3(\epsilon)&=& n_3\epsilon
\eeann
and $n_{i,j}$ are positive constants.  Observe that $r_i(\epsilon) \rightarrow 0$ as $\epsilon \rightarrow 0$.

At one-time step,
\beann
\vt x_{k+1}-\bar{x}_{k+1} \vt &\leq& m_2(\epsilon) \vt x_k-\bar{x}_k \vt + m_3(\epsilon) \V \psi-\bar{\psi} \V
\eeann
where
\[m_2(\epsilon)=\frac{\alpha + r_3(\epsilon)}{1-r_1(\epsilon)} \]
and
\[m_3(\epsilon)=\frac{r_2(\epsilon)}{1-r_1(\epsilon)}.\]
By invoking Gronwall's inequality in finite difference form (\ref{gronwall1}) and assuming that for some small $\epsilon > 0$
\[m_2(\epsilon)<1,\]
then
\bea \label{xphibound}
\vt x_{k}-\bar{x}_{k} \vt &\leq& m_3(\epsilon)\left[ \sum_{j=0}^{k-1} m_2(\epsilon)^{k-1-j}\right] \V \psi-\bar{\psi} \V \nn\\
&\leq& m_3(\epsilon) \frac{1-m_2(\epsilon)^k}{1-m_2(\epsilon)} \V \psi-\bar{\psi} \V \nn \\
&\leq& \frac{m_3(\epsilon)}{1-m_2(\epsilon)} \V \psi-\bar{\psi} \V 
\eea

With the bounds (\ref{gphibnd2}) and (\ref{xphibound}), we get
\beann
\vt T\psi(x_0)- T\bar{\psi}(x_0)\vt &\leq& \left\vt  \sum_{l=0}^{\infty}(A'+K'B')^{l}( h_{\psi}(x_l,x_{l+1})-h_{\bar{\psi}}(\bar{x}_l,\bar{x}_{l+1}) )\right\vt\\
&\leq&  \sum_{l=0}^{\infty}\alpha^{l} \left\vt h_{\psi}(x_l,x_{l+1})-h_{\bar{\psi}}(\bar{x}_l,\bar{x}_{l+1}) \right\vt \\
&\leq&  \sum_{l=0}^{\infty}\alpha^{l} \left[ r_1(\epsilon)\vt x_{l+1}-\bar{x}_{l+1} \vt + r_2(\epsilon) \V \psi- \bar{\psi} \V + r_3(\epsilon)\vt x_l-\bar{x}_l \vt \right]\\
&\leq&  \sum_{l=0}^{\infty}\alpha^{l} \left[ 2(r_1(\epsilon)+ r_3(\epsilon)) \frac{m_3(\epsilon)}{1-m_2(\epsilon)} \V \psi-\bar{\psi} \V + r_2(\epsilon) \V \psi- \bar{\psi} \V  \right].
\eeann

It follows that
\beann
\left\vt T\psi(x_0)- T\bar{\psi}(\bar{x}_0) \right\vt &\leq& c(\epsilon) \V \psi - \bar{\psi} \V
\eeann
where
\[c(\epsilon)=\frac{2m_3(\epsilon)(r_1(\epsilon)+ r_3(\epsilon))}{(1-m_2(\epsilon))(1-\alpha)} + \frac{r_2(\epsilon) }{1-\alpha}.\]
Notice that $c(\epsilon) \rightarrow 0$ as $\epsilon \rightarrow 0$.
Thus, for sufficiently small $\epsilon>0$
\[c(\epsilon)<1.\]
Hence, $T$ is a contraction on $\mathbb{X}$ for $\epsilon$ sufficiently small.  Hence there exists a unique $\psi \in \mathbb{X}$ such that
\[\psi=T\psi.\]
\end{proof}

\section{Some Properties}

\subsection{Eigenstructure}

Recall from the previous chapter the bidirectional linear Hamiltonian dynamics
\bea \label{eigdyn}
\left[ \begin{array}{c}
x^+ \\
\lambda \\
\end{array}\right]&=&
\mathbb{H}
\left[\begin{array}{c}
x \\
\lambda^+ \\
\end{array}\right]
\eea
where
\beann
\mathbb{H}=\left[\begin{array}{cc}
A & -BR^{-1}B'\\
Q & A'\\
\end{array}\right].
\eeann

\begin{definition} \label{sfactor}
Suppose
\bea \label{therelation}
\mathbb{H}
\left(\begin{array}{c}
\delta x \\
\mu \delta \lambda \\
\end{array}\right)=
\left(\begin{array}{c}
 \mu \delta x \\
\delta \lambda \\
\end{array}\right).
\eea
Then we call $\mu$ the eigenvalue of the dynamics (\ref{eigdyn}).
\end{definition}
In the case where we do not assume $0$ as the eigenvalue of $\mathbb{H}$, we found that the forward Hamiltonian $\mathbb{H}^F$ is a hyperbolic system in Theorem 4.4.2.  Similarly, we would like to show that the linear bidirectional Hamiltonian matrix $\mathbb{H}$ is hyperbolic; i.e., the eigenvalue of $\mathbb{H}$ lies strictly inside and outside the unit circle.

\begin{theorem}
If $\mu$ is an eigenvalue of the dynamics (\ref{eigdyn}) satisfying the relation (\ref{therelation}), then $\frac{1}{\mu}$ is also an eigenvalue of $\mathbb{H}$.
\end{theorem}
\begin{proof}
First, we decompose $\mathbb{H}$ into
\beann
\left[\begin{array}{cc}
0 & I\\
I & 0\\
\end{array}\right]
\left[\begin{array}{cc}
Q & A'\\
A & -BR^{-1}B'\\
\end{array}\right]
=\mathbb{H}.
\eeann
Let's call
\beann
\mathbb{S}=\left[\begin{array}{cc}
Q & A'\\
A & -BR^{-1}B'\\
\end{array}\right]
\eeann
and notice that $\mathbb{S}$ is symmetric.

From the equation (\ref{therelation}), we have the following
\beann
\mathbb{H}
\left(\begin{array}{cc}
I & 0 \\
0 & \mu I\\
\end{array}\right)
\left(\begin{array}{c}
\delta x \\
\delta \lambda \\
\end{array}\right)
=
\left(\begin{array}{cc}
\mu I & 0 \\
0    & I\\
\end{array}\right)
\left(\begin{array}{c}
\delta x \\
\delta \lambda \\
\end{array}\right)
\eeann
which is equivalent to
\beann
\left(\begin{array}{cc}
\frac{1}{\mu} I & 0 \\
0    & I\\
\end{array}\right)
\mathbb{H}
\left(\begin{array}{cc}
I & 0 \\
0 & \mu I\\
\end{array}\right)
\left(\begin{array}{c}
\delta x \\
\delta \lambda \\
\end{array}\right)
=
\left(\begin{array}{c}
\delta x \\
\delta \lambda \\
\end{array}\right).
\eeann
We denote
\beann
\mathbb{H}_{\mu}=
\left(\begin{array}{cc}
\frac{1}{\mu} I & 0 \\
0    & I\\
\end{array}\right)
\mathbb{H}
\left(\begin{array}{cc}
I & 0 \\
0 & \mu I\\
\end{array}\right).
\eeann
Observe that $1$ is an eigenvalue of $\mathbb{H}_{\mu}$.  It follows that
\[det \left[I- \mathbb{H}_{\mu} \right]=0 \Longrightarrow det \left[I- \mathbb{H}'_{\mu} \right]=0\]
where
\beann
\mathbb{H}'_{\mu}=
\left(\begin{array}{cc}
I & 0 \\
0    & \mu I\\
\end{array}\right)
\mathbb{S}
\left(\begin{array}{cc}
0 & I \\
I & 0 \\
\end{array}\right)
\left(\begin{array}{cc}
\frac{1}{\mu}I & 0 \\
 0             & I \\
\end{array}\right).
\eeann
Again, $1$ is an eigenvalue of $\mathbb{H}_{\mu}$; i.e.,
\beann \label{theadjoint}
\mathbb{H}'_{\mu}
\left(\begin{array}{c}
\widetilde{\delta x}\\
\widetilde{\delta \lambda} \\
\end{array}\right)=
\left(\begin{array}{c}
\widetilde{\delta x}\\
\widetilde{\delta \lambda} \\
\end{array}\right).
\eeann
From above, we have
\beann
\mathbb{S}
\left(\begin{array}{cc}
0 & I \\
I & 0 \\
\end{array}\right)
\left(\begin{array}{c}
\frac{1}{\mu} \widetilde{\delta x}\\
\widetilde{\delta \lambda}\\
\end{array}\right)=
\left(\begin{array}{c}
\widetilde{\delta x}\\
\frac{1}{\mu} \widetilde{\delta \lambda}\\
\end{array}\right)
\eeann
which is equal to
\beann
\mathbb{S}
\left(\begin{array}{c}
\widetilde{\delta \lambda}\\
\frac{1}{\mu} \widetilde{\delta x}\\
\end{array}\right)=
\left(\begin{array}{cc}
0 & I \\
I & 0 \\
\end{array}\right)
\left(\begin{array}{c}
\frac{1}{\mu} \widetilde{\delta \lambda}\\
\widetilde{\delta x}\\
\end{array}\right).
\eeann
Thus,
\beann \label{theadjoint}
\mathbb{H}
\left(\begin{array}{c}
\widetilde{\delta \lambda} \\
\frac{1}{\mu} \widetilde{\delta x}\\
\end{array}\right)=
\left(\begin{array}{c}
\frac{1}{\mu} \widetilde{\delta \lambda} \\
\widetilde{\delta x}\\
\end{array}\right).
\eeann
Hence, $\frac{1}{\mu}$ is an eigenvalue of $\mathbb{H}$.
\end{proof}
Note that $\mu \neq 1$.  The eigenvalue $\mu=1$ corresponds to the trivial eigenvector 0.  Also, the infinite eigenvalues $0$ and $\infty$ satisfy (\ref{therelation}) due to the singularity of A.

\subsection{Symplectic Form}
The nonlinear tangent dynamics of the bidirectional Hamiltonian system (\ref{nonham}) that we derive in Chapter 3 is
\bea \label{nonhampert2}
\left[ \begin{array}{c}
\delta x^+\\
\delta \lambda \\
\end{array}\right]=
\mathbb{H}_{\delta,k}(x,\lambda^+)
\left[ \begin{array}{c}
\delta x\\
\delta \lambda^+\\
\end{array}\right]
\eea
where
\beann
\mathbb{H}_{\delta,k}(x,\lambda^+)=
\left[ \begin{array}{cc}
H_{\lambda^+ x} & H_{\lambda^+ \lambda^+}\\
H_{x x}         & H_{x \lambda^+}
\end{array}\right](x,\lambda^+).
\eeann
We denote $H_{x \lambda^+},~H_{x x},~H_{\lambda^+ \lambda^+}$ as partial derivatives.  Recall the nondegenerate and bilinear symplectic two-from $\Omega:T_{(x,\lambda)}\mathcal{M} \times T_{(x,\lambda)}\mathcal{M} \mapsto \mathbb{R}$,
\beann
\Omega(v,w)=v'Jw
\mbox{   and   }J=
\left[ \begin{array}{cc}
0  & I \\
-I & 0 \\
\end{array}\right]
\eeann
where
\[\Omega(v,w)=-\Omega(w,v)\]
and
\[ v=\left[ \begin{array}{c}
\delta x\\
\delta \lambda \\
\end{array}\right],~~~ w=\left[ \begin{array}{c}
\widetilde{\delta x}\\
\widetilde{\delta \lambda} \\
\end{array}\right]\]
$(x,\lambda^+)\in \mathcal{M}$ and $(v,w) \in T_{(x,\lambda)}\mathcal{M}$. Also, $\mathcal{M}=T^*\mathcal{N}$ where $x \in \mathcal{N}$.
We would like to show that under the tangent dynamics (\ref{nonhampert2}), the two-form $\Omega$ is invariant; i.e.,
\bea \label{preservingform}
\Omega(v,w)=\Omega(v^+,w^+).
\eea
Then,
\bea \label{Omegalhs}
\Omega(v,w)&=&v'Jw \nn \\
&=&
\left(\begin{array}{cc}
\delta x' & \delta x' H_{x x} + \delta \lambda^{+'} H_{x \lambda^+} \\
\end{array}\right) J
\left(\begin{array}{c}
\widetilde{\delta x}' \\
H_{x x}\widetilde{\delta x} + H_{x \lambda^+}\widetilde{\delta \lambda}^+ \\
\end{array}\right)\nn \\
&=&-\delta \lambda^{+'} H_{x \lambda^+}\widetilde{\delta x} + \delta x' H_{x \lambda^+} \widetilde{\delta \lambda}^+
\eea
and
\bea \label{Omegarhs}
\Omega(v^+,w^+)&=&v^{+'}Jw^+ \nn \\
&=&
\left(\begin{array}{cc}
\delta x'H_{\lambda^+ x} + \delta \lambda^{+'} H_{\lambda^+ \lambda^+}' & \delta \lambda^+\\
\end{array}\right)J
\left(\begin{array}{c}
H_{\delta^+ x}\widetilde{\delta x} + H_{\lambda^+ \lambda^+}\widetilde{\delta \lambda}^+ \\
\widetilde{\delta \lambda}^+ \\
\end{array}\right)\nn \\
&=&-\delta \lambda^{+'} H_{\lambda^+ x}\widetilde{\delta x} + \delta x' H_{\lambda^+ x} \widetilde{\delta \lambda}^+
\eea
Since (\ref{Omegarhs}) and (\ref{Omegalhs}) are equal, then
\[\Omega(v,w)=\Omega(v^+,w^+).\]
Thus, for any two tangent vectors satisfying the dynamics (\ref{nonhampert2}), the value of $\Omega$ does not change.

\section{Lagrangian Submanifold}
We have two-form 
\bea \label{Omega2form}
\Omega(v,w)=-\delta \lambda^{+'} H_{\lambda^+ x}\widetilde{\delta x} + \delta x' H_{\lambda^+ x} \widetilde{\delta \lambda}^+.
\eea 
Using the technique in Chapter 3, we have the following tangent dynamics around the trajectories $(x_j,x_{j+1})$
\bea
\delta x_{j+1}   &=& \Big( (A+BK)  + \frac{\partial f_{\psi}}{\partial x_j}(x_j,x_{j+1}) \Big) \delta x_j \label{tandyn1}\\
&&~~~~~~~~~~~+ \frac{\partial f_{\psi}}{\partial x_{j+1}}(x_j, x_{j+1}) \delta x_{j+1} \nn \\
\delta \lambda_j &=& (A'+K'B')\delta \lambda_{j+1} + \frac{\partial h_{\psi}}{\partial x_{j+1}}(x_j, x_{j+1}) \delta x_{j+1} \nn 
\eea

By the Inverse Function Theorem, we can choose $\epsilon >0 $ small enough so that $I-\frac{\partial f_{\psi}}{\partial x_{k+1}}(x_k, x_{k+1})$ is invertible for $ \vt x_k \vt, \vt x_{k+1} \vt < \epsilon$ because
\[I-\frac{\partial f_{\psi}}{\partial x_{k+1}}(0, 0)= I\]
since (\ref{fgcond1}).

It follows that the tangent state dynamics is equivalent to
\beann
\delta x_{k+1}= \Big( I-\frac{\partial f_{\psi}}{\partial x_{k+1}}(x_k, x_{k+1}) \Big)^{-1} \Big( (A+BK)  + \frac{\partial f_{\psi}}{\partial x_k}(x_k,x_{k+1})\Big) \delta x_k
\eeann
and
\bea \label{deltaxsoln}
\delta x_{k+1}&=& \prod_{i=0}^{k}\Big( I-\frac{\partial f_{\psi}}{\partial x_{i+1}}(x_i, x_{i+1})\Big)^{-1} \cdot \\
&&~~~~~~~~~~\Big( (A+BK)  + \frac{\partial f_{\psi}}{\partial x_i}(x_i,x_{i+1})\Big) \delta x_0. \nn
\eea
for $ \vt x_k \vt, \vt x_{k+1} \vt < \epsilon$.

As we let $k \rightarrow \infty$, we have $x_k \rightarrow 0$, $\frac{\partial f_{\psi}}{\partial (x,x^+)}(0) \rightarrow 0$.  It follows that from (\ref{deltaxsoln})
\[\delta x_{k+1} \rightarrow (A+BK)^k \delta x_0 \rightarrow 0.\]

Since $\delta x_k \rightarrow 0$ as $k \rightarrow \infty$, we have that
\[\Omega(v,w) \rightarrow 0\]
for $v,w$ restricted to the tangent dynamics.  Thus, the local stable manifold is a Lagrangian submanifold.  Similarly as in Chapter 3,  we have
\bea \label{closedtoo}
\frac{\partial \psi_i}{\partial x_{j}}-\frac{\partial \psi_j}{\partial
x_{i}}=0\;\;\mbox{  for  }i,j=1,\ldots,n.
\eea
The equation (\ref{closedtoo}) implies that $\psi(x)$ is closed.  Then, by 
the
\emph{Stokes' Theorem} there exists $\bar{\pi} \in \mathcal{C}^{r}(\mathbb{R}^n)$
such that
\bea \label{voila}
\psi(x)=\frac{\partial \bar{\pi}}{\partial x}(x)\mbox{  where  }\psi(0)=0
\eea
locally on some neighborhood of 0.  Since $z=\psi(x)=\lambda - Px$, we have that
\beann
\lambda&=& \frac{\partial \bar{\pi}}{\partial x}(x) + Px\\
       &=& \frac{\partial \bar{\pi}}{\partial x}(x) + \frac{\partial}{\partial x}\left(\frac{1}{2}x'Px\right)
\eeann
because $P$ is symmetric.  Thus, there exists $\pi \in \mathcal{C}^{r}$ such that
\[\lambda=\frac{\partial \pi}{\partial x}(x)  \]
where
\[\frac{\partial \pi}{\partial x}(x)= \frac{\partial \bar{\pi}}{\partial x}(x) + \frac{\partial}{\partial x}\left(\frac{1}{2}x'Px\right).\]
Hence, the local stable manifold $\lambda$ is the gradient of the optimal cost for the bidirectional Hamiltonian dynamics.

%% file: chapter6_6.tex
Our procedure for solving HJB PDE numerically is described in \cite{Na00}, but we give a summary of the method and some numerical results in this chapter.

\section{A Method for HJB}
We solve the Hamilton Jacobi Bellman (HJB) Partial Differential
Equation that arises in many control problems.  We consider the infinite horizon 
optimal control problem of
minimizing the cost
\bea \label{cost}
\int_t^{\infty} l(x,u) \ dt
\eea
subject to the dynamics
\bea
\dot{x}=f(x,u)\label{dynamics1}
\eea
and initial condition
\bea
x(t)=x^0.
\eea
The state vector $x$ is an $n$ dimensional column vector, the control 
$u$ is an $m$ dimensional column
vector and the dynamics $f(x,u)$ and Lagrangian $l(x,u)$ are assumed 
to be sufficiently
smooth.

If the minimum exists and is a smooth function
$\pi(x^0)$ of the initial condition
then it satisfies the HJB PDE
\bea 
\label{HJB1} 
\min_u\left\{ \frac{\partial \pi}{\partial x}(x)f(x,u)+l(x,u)\right\}=0
\eea
and the optimal control $\kappa(x)$ satisfies
\bea 
\label{HJB2}
\kappa(x)=\mbox{arg}\min_u\left\{ \frac{\partial \pi}{\partial
x}(x)f(x,u)+l(x,u)\right\}=0
\eea

These are expressed in terms of the Hamiltonian
\bea 
\label{ham}
H(p,x,u)=pf(x,u) + l(x,u)
\eea
where the argument $p$ is an $n$ dimensional row vector. The HJB PDE becomes
\bea \label{hhjb1}
0&=&\min_u H(\frac{\partial \pi}{\partial x}(x),x,u)\\
\label{hhjb2}
\kappa(x)&=&\mbox{arg}\min_u H(\frac{\partial \pi}{\partial x}(x),x,u)
\eea

We assume that the Hamiltonian $H(p,x,u)$ is strictly convex in $u$ for all $p,\ x$.  Then (\ref{HJB1}, \ref{HJB2}) become
\bea \label{HJB3}  \frac{\partial \pi}{\partial
x}(x)f(x,\kappa(x))+l(x,\kappa(x))=0
\eea
and
\bea \label{HJB4}  \frac{\partial \pi}{\partial
x}(x)\frac{\partial f}{\partial u}(x,\kappa(x))+\frac{\partial l}{\partial
u}(x,\kappa(x))=0
\eea

These are the equations (\ref{HJB3}),(\ref{HJB4}) that we solve for the optimal cost $\pi$ and optimal control $\kappa$.

\subsection{Al'brecht's Method Revisited}
Al'brecht \cite{Al61} solved the HJB PDE locally around zero by expanding the problem in a power series,
\bea
f(x,u)&=& Ax+Bu+f^{[2]}(x,u) \nonumber\\
       & & +f^{[3]}(x,u)+\ldots \\
l(x,u)&=&{1\over 2} \left(x'Qx+2x'Su+u'Ru\right)\nonumber \\
&&+l^{[3]}(x,u)+l^{[4]}(x,u)+\ldots\\
\pi(x)&=&{1\over 2} x'Px+\pi^{[3]}(x)\nonumber\\
       & &+\pi^{[4]}(x)+\ldots\\
\kappa(x)&=& Kx+\kappa^{[2]}(x)+\kappa^{[3]}(x)+\ldots\label{expan1}
\eea
where ${\cdot}^{[d]}$ denotes a homogeneous polynomial of degree $d$.
As in Chapter 3, we subsitute these equations into the HJB PDE (\ref{HJB3}, \ref{HJB4}) and equate
the terms of like degree to obtain a  sequence of algebraic equations for the
unknowns. 

The first level is the pair of equations obtained by collecting the
quadratic terms of (\ref{HJB3}) and the linear terms of (\ref{HJB4}). 
We denote the $d^{th}$ level as the
pair of equations garnered from the $[d+1]^{th}$ degree terms of 
(\ref{HJB3}) and the $d^{th}$ degree
terms of (\ref{HJB4}).

Here are the first set of equations:
\bea
0&=&A'P +PA + Q - \nonumber \\
  & &(PB +S) R^{-1}(PB +S)'\label{alb12}\\
K&=&-R^{-1}(PB +S)'\label{alb21}
\eea

The quadratic terms of the HJB PDE reduce to the familiar Riccati equation
(\ref{alb12}) and the linear optimal feedback (\ref{alb21}).

We assume $A,\ B $ is stabilizable and $Q,\ A$ is detectable then the Riccati
equation has a unique positive definite solution $P$ and the linear feedback locally
exponentially stabilizes the closed loop system.  Moreover the optimal quadratic cost is a local Lyapunov 
function for the closed loop system.

The $d^{th}$ Level

Suppose we have solved through the $d-1^{th}$ level after repeating the process $d-1$ times.  It is convenient to
incorporate this solution into the dynamics and the cost.  Let $\kappa^{k]}(x)=Kx+\kappa^{[2]}(x)+\kappa^{[3]}(x)+\ldots+\kappa^{[k]}(x)$ and define
\bea
\bar{f}(x,u) &=& f(x,\kappa^{d-1]}(x)+u)\\
\bar{l}(x,u) &=& l(x,\kappa^{d-1]}(x)+u)
\eea
These have power series expansions through terms of degree $d$ and $d+1$
of the form
\beann
\bar{f}(x,u)&=& (A+BK)x+Bu+\bar{f}^{[2]}(x,u)+ \ldots 
\bar{f}^{[d]}(x,u)+\ldots \nonumber\\
\bar{l}(x,u)&=&{1\over 2} \left(x'Qx+2x'SKx+x'K'RKx \right)\\&&+x'Su+u'Ru+\\
&&\bar{l}^{[3]}(x,u)+\ldots +\bar{l}^{[d+1]}(x,u)+
\ldots \nonumber\\
\eeann

We plug these into the HJB PDE (\ref{HJB3}, \ref{HJB4}) and find that 
they are satisfied
through the $d-1$ level and don't involve $u$.  Since 
$u=\kappa^{[d]}(x)+\ldots$, the $d$
level equations are
\bea
0 &=& \displaystyle \frac{\partial \pi^{[d+1]}}{\partial x}(x)(A + BK)x \\
&&\displaystyle +\sum_{i=2}^{d -1}\frac{\partial \pi^{[d+2-i]}}{\partial x}(x)\bar{f}^{[i]}(x,0)
+x'PBu \nonumber \\&&+\bar{l}^{[d+1]}(x,0) +x'Su+{1\over 2} x'K'Ru\nonumber
\eea
and
\bea
0&=&\displaystyle \frac{\partial \pi^{[d+1]} }{\partial x}(x)B +
\sum_{i=2}^d\frac{\partial \pi^{[d+2-i]}}{\partial x}\frac{\partial 
\bar{f}^{[i]} }{\partial u}(x,0)
\nonumber\\ & &\displaystyle + \frac{\partial \bar{l}^{[d+1]} 
}{\partial u}(x,0) +u'R
\eea

Because of (\ref{alb21}), $u$ drops out of the first equation which becomes
\bea
0 &=& \displaystyle \frac{\partial \pi^{[d+1]}}{\partial x}(x)(A + 
BK)x \label{alb1d+1}\\
&&\displaystyle +\sum_{i=2}^{d -1}\frac{\partial 
\pi^{[d+2-i]}}{\partial x}(x)\bar{f}^{[i]}(x,0)
\nonumber \\&&+\bar{l}^{[d+1]}(x,0). \nonumber
\eea
After this has been solved for $\pi^{[d+1]}(x)$, we can solve the 
second for $\kappa^{[d]}(x)$,
\bea
\kappa^{[d]}(x)&=&-R^{-1}\left(\displaystyle \frac{\partial 
\pi^{[d+1]} }{\partial x}(x)B \right.\\
&&\left.+
\sum_{i=2}^d\frac{\partial \pi^{[d+2-i]}}{\partial x}\frac{\partial 
\bar{f}^{[i]} }{\partial u}(x,0)
\nonumber\displaystyle + \frac{\partial \bar{l}^{[d+1]} }{\partial 
u}(x,0) \right) \label{alb2d}
\eea
These equations admit a unique solution up to the smoothness of $f$ 
and $l$ since $A+BK$ has all its eigenvalues in the left half plane.  If $f$ and $l$ are
real analytic, the power series converges to the solution of the HJB 
PDE locally around $x=0$ \cite{Lu69}, \cite{Kr98}.  The higher degree equations are linear; hence they are easily solvable.

\begin{figure}[tbp]
\centering
\scalebox{.4}{\includegraphics{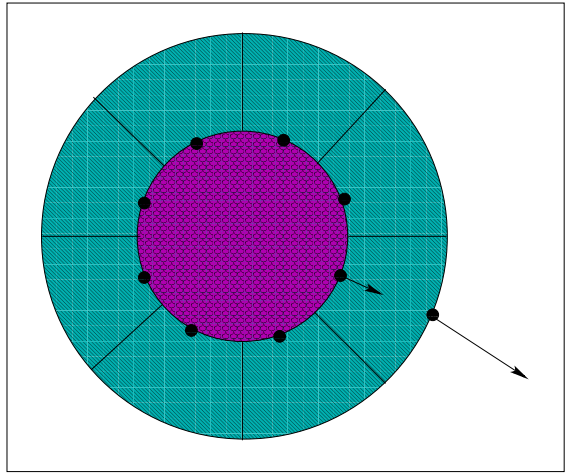}}
\caption{Around each point, we generate a polynomial solution on each patch.}
\label{fig10}
\end{figure}

Albrecht's approach develops the Taylor series expansion along a point.  The Taylor series only converges on a small neigbhorhood and may quickly diverges outside that neighborhood.  Also increasing the degree of the polynomial solutions does not guarantee a larger region where the approximated solutions are true.  In our approach, we use Albrecht's method to generate the initial approximations around $0$.  Then the initial approximations are then improve by piecing successive approximations of smooth solutions around the points on the level sets, see Fig.(\ref{fig10}).  Instead we truncate the computed solutions at the point where the optimal cost satisfies the optimality and stability constraints and begin new approximations at the same point.  

\subsection{The Lyapunov Criterion}
Let $\pi^{d+1]}$ and $\kappa^{d]}$ denote the approximate cost and feedback to
degrees $d+1$ and $d$ respectively.  We assume $l(x,u) \geq 0$.  We 
find the largest sublevel set
\bea \label{radc}
\pi^{d+1]}(x)\leq c
\eea
  such that
\bea \label{appLy}
  \frac{\partial \pi^{d+1]}}{\partial
x}(x)f(x,\kappa^{d]}(x))\leq -(1- \epsilon_1) l(x,\kappa^{[d]}(x))\\
-(1- \epsilon_2) l(x,\kappa^{[d]}(x)) \leq \frac{\partial \pi^{d+1]}}{\partial
x}(x)f(x,\kappa^{d]}(x))
\eea 
The parameter $\epsilon_1$ controls the rate of exponential stability
of the closed loop system while $\epsilon_2$ dictates the rate of the optimality of the feedback.  We determine a sublevel set of the approximate cost on which it is an acceptable
Lyapunov function and the approximate feedback is stabilizing and satisfying the optimality conditions.  We
emphasize the stabilizing property of the control law rather than its optimality because
usually optimality is only a tool to find a stabilizing feedback. 
Typically the goal is to find a control law to stabilize the system and the optimal control problem is
formulated as a way of finding one.


\subsection{Power Series Expansions}

We specialize to problems with dynamics that is affine in $u$
and with a Lagrangian that is quadratic in $u$,
\bea
f(x,u)&=&g_0(x)+g_1(x)u\\
l(x,u)&=&l_0(x) +l_1(x)u+u'l_2(x) u
\eea
where $l_2(x)$ is an invertible $m\times m$ matrix for all $x$.

Assume that  we have solved the HJB PDE locally around 0 and
  have choosen a sublevel set of value $c$ subject to the condition 
(\ref{radc}, \ref{appLy}).  We generate a power series solution around a point $\bar{x}$ on the level set $ \pi^{d+1]}(x)= c$.

We introduce some notation.  Let $\alpha=(\alpha_1,\ldots,\alpha_n)$
be a multi-index of nonnegative integers and $|\alpha|=\sum_i \alpha_i$.
  Let $\beta=(\beta_1,\ldots,\beta_n)$.
We say $\beta\leq \alpha $ if $\beta_i\leq \alpha_i,\ i=1,\ldots ,n$ and
$\beta< \alpha $ if $\beta\leq \alpha $ and for at least one $i$,
$\beta_i < \alpha_i$.  Let ${\bf 0}=(0,\ldots,0)$.

Define the differential operator
$$D^\alpha = (\frac{\partial}{\partial x_1})^{\alpha_1}\ldots
(\frac{\partial}{\partial x_n})^{\alpha_n},$$
the multifactorial
$$\alpha ! = \alpha_1!\ldots \alpha_n!,$$
the monomial
$$x^\alpha= x_1^{\alpha_1}\ldots x_n^{\alpha_n},$$
and the coefficient
$$ C(\alpha,\beta)=\left(\begin{array}{c} 
\alpha_1\\\beta_1\end{array}\right)\ldots
\left(\begin{array}{c} \alpha_n\\\beta_n\end{array}\right)
$$
where $\beta\leq \alpha $.

We derive a system of equations  for
\[D^\alpha \pi(\bar{x}), \ \ \ \ \ \ D^\alpha \kappa(\bar{x}).\]  We already know that if $\alpha$ is the $i^{th}$ unit vector
\beann
D^{\bf 0} \pi(\bar{x})&=& \pi(\bar{x})\\
D^{\bf 0} \kappa(\bar{x})&=& \kappa(\bar{x})\\
D^\alpha \pi(\bar{x})&=&\frac{\partial \pi}{\partial x_i}(\bar{x}).
\eeann

Let the Cauchy data $D^\alpha \pi(\bar{x})$ and $D^\alpha \kappa(\bar{x})$ be known for $\alpha=(0,\alpha_2,\alpha_3,\hdots,\alpha_n)$ and $f_1(\bar{x},\kappa(\bar{x})) \neq 0$.

Assume that we have derived algebraic equations $D^\beta \pi(\bar{x})$ for $D^\beta \pi(\bar{x})$ and $D^\beta \kappa(\bar{x})$ for $\beta<\alpha$.
We apply $D^\alpha$ to (\ref{HJB3}) to
obtain
\bea \label{deq}
0&=&  \frac{\partial }{\partial
x}(D^\alpha\pi(\bar{x}))f(\bar{x},\kappa(\bar{x}))\\
&&+\sum_{{\bf 0}<\beta\leq\alpha}C(\alpha,\beta) (D^{\alpha-\beta} 
\frac{\partial \pi }{\partial
x}(\bar{x})) D^\beta f(\bar{x},\kappa(\bar{x})) \nonumber \\ &&
+ D^\alpha l(\bar{x},\kappa(\bar{x})). \nonumber
\eea

This yields an equation for $\frac{\partial D^\alpha \pi}{\partial x_1} (\bar{x})$ because $f_1(\bar{x},\kappa(\bar{x})) \neq 0$ and all the other terms $\frac{\partial D^\alpha \pi}{\partial x_i} (\bar{x})$ for $i \neq 1$ are known from the Cauchy data.

We apply $D^\alpha$ to (\ref{HJB4}) to
obtain
\bea \label{aeq}
0&=&\sum_{{\bf 0}\leq \beta\leq\alpha}C(\alpha,\beta) 
(D^{\alpha-\beta} \frac{\partial \pi }{\partial
x}(\bar{x})) D^\beta g_1(\bar{x}) \nonumber \\ &&
+ D^\alpha l_1(\bar{x}) + (D^\alpha\kappa(\bar{x}))'l_2(\bar{x}) \nonumber\\
&&\sum_{{\bf 0}< \beta\leq\alpha}C(\alpha,\beta) 
(D^{\alpha-\beta}\kappa(\bar{x}))'D^\beta l_2(\bar{x})
\eea

Notice that this equation (\ref{aeq}) only contains 
$D^\alpha\kappa(x)$ in one term
multiplied by an invertible matrix so we can express 
$D^\alpha\kappa(x)$ as a function of
$D^\beta\kappa(x)$ for
$\beta<\alpha$ and $D^\gamma\pi(x)$.

In this way  we obtain the approximations
\bea
\pi(x)&\approx& \sum_d \sum_{|\alpha|\leq d} {1\over \alpha !} D^\alpha 
\pi(\bar{x}) (x-\bar{x}))^\alpha\\
\kappa(x)&\approx&\sum_d \sum_{|\alpha|\leq d} {1\over \alpha !} D^\alpha 
\kappa(\bar{x}) (x-\bar{x})^\alpha
\eea
where $x$ is close to $\bar{x}$.  

\section{An Example}
Consider the optimal control problem:
\[ \min_u  \int^{\infty}_0 \ln^2(x+1) + u^2 dt\]
\bea \label{pragerex}
\text{subject to}  \\
\dot{x}&=& xu + u \nn \\
x(0)&=&x_0 \nn
\eea
where the domain $\mathcal{D}=\{x \in \mathbb{R} \vt x > -1)\}$.  Given the problem (\ref{pragerex}) and the schemes described above, we solve the optimal control $\kappa(x)$ and optimal cost $\pi(x)$ up to degree $d$ and $d+1$, respectively; i.e.
\beann
\pi(x)&=& \pi^{[2]}(x) + \pi^{[3]}(x) + \hdots + \pi^{[d+1]}(x)\\
\kappa(x)&=& \kappa^{[1]}(x) + \kappa^{[2]}(x) + \hdots + \kappa^{[d]}(x) 
\eeann
We fix the degree at $d=3$ for this example.  The interval is $\mathcal{D}=(-1,4]$.  Note that the analytic solutions are $\kappa_*(x)=-\ln (x+1)$ and $\pi_*(x)=\ln^2(x+1)$. 

\subsection{Approximation I:  Albrecht's Method}
\begin{figure}[tbp]
\centering
\scalebox{.8}{\includegraphics{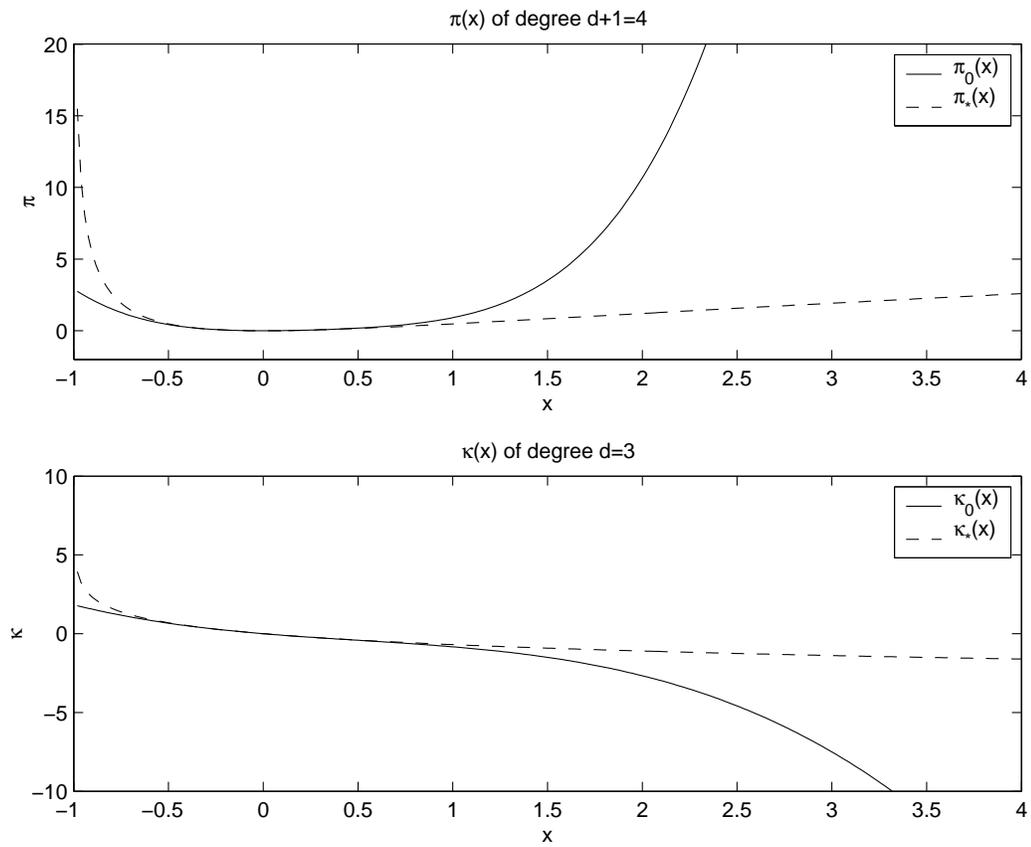}}
\caption{The approximated solutions, $\pi_0(x)$ and $\kappa_0(x)$, via Al'brecht's method compared to the real solutions, $\pi_*$ and $\kappa_*$.}
\label{fig3}
\end{figure}

We implement Albrecht's method by using the MATLAB code {\it hjb.m} in \cite{Kr97} to find the coefficients of $\pi$ and $\kappa$.  We then set-up the polynomials.  At each point $x_j \in \mathcal{D}$, we assign the polynomial approximations $\pi_{0_j}=\pi_0(x_j)$ and $\kappa_{0_j}=\kappa_0(x_j)$.  See Fig.(\ref{fig3}).  We solve the optimization problem (\ref{radc}).  We use the MATLAB code, {\it fmincon.m}, iteratively to find a point on the level set $\pi^{d+1]}(x)<c$; i.e., $\bar{x}\in \mathcal{D}_l$.  For this example, we march along the $x$-axis in the both directions.  For the interval $[0,4]$, we march on the x-axis towards $\infty$, while on interval $(-1,0]$ we move towards -$1$.  See Fig.(\ref{fig6}) for the psuedocode of the algorithm {\it DRIVER1.m} and the actual codes in the appendix.

\begin{figure}[tbp]
\centering
\scalebox{.8}{\includegraphics{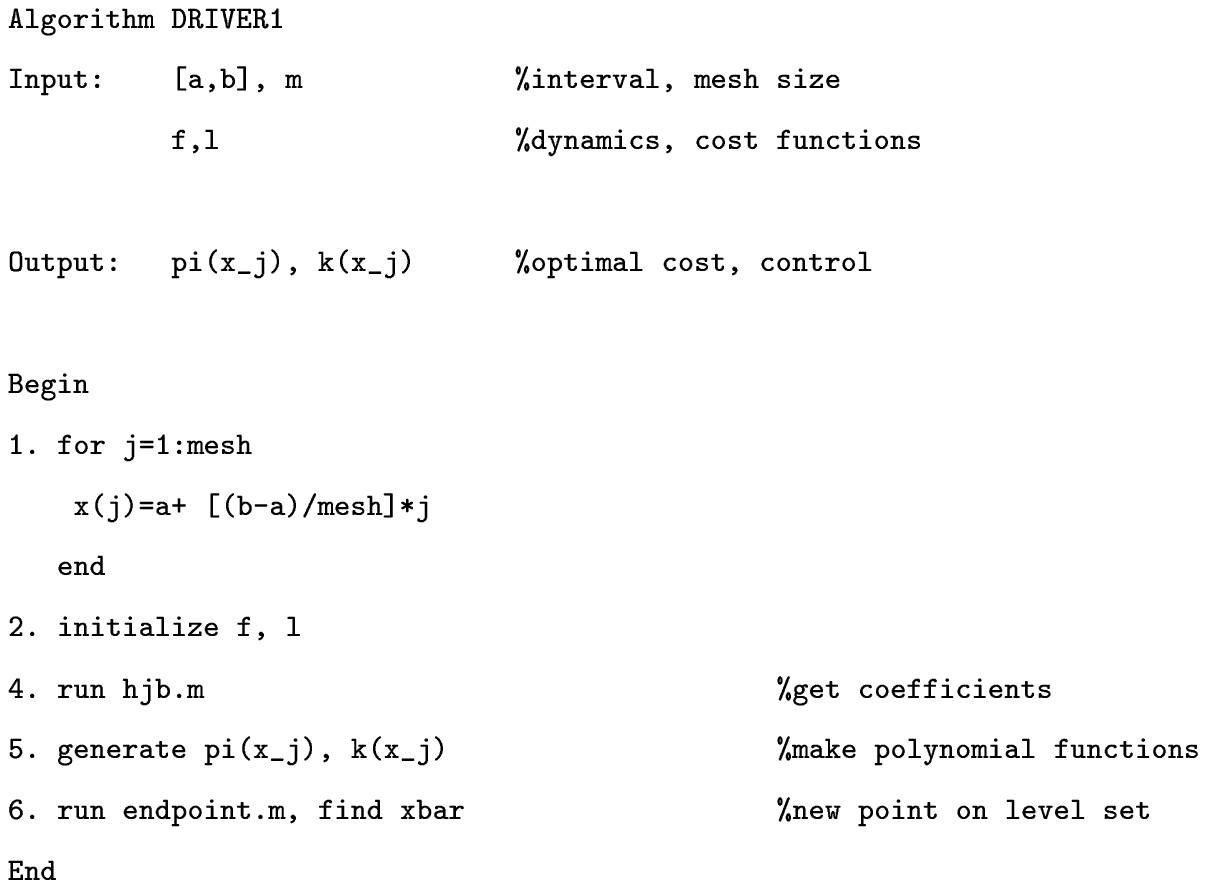}}
\caption{The psuedocode for the initial approximation.}
\label{fig6}
\end{figure}

\subsection{Approximation II}
\begin{figure}[tbp]
\centering
\scalebox{.8}{\includegraphics{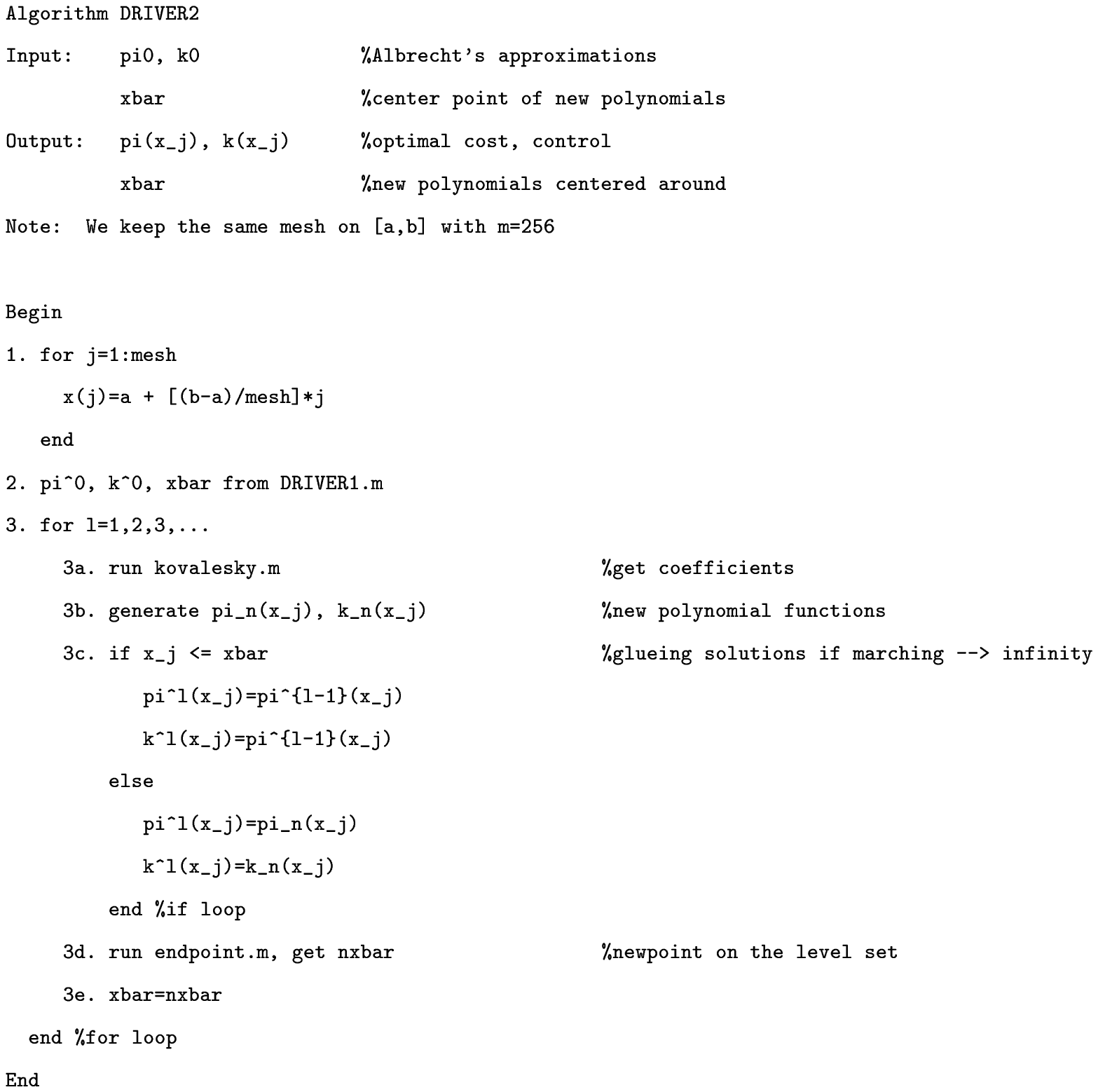}}
\caption{The psuedocode for $l$th approximations where $l=1,2,\hdots$}
\label{fig7}
\end{figure}

\begin{figure}[tbp]
\centering
\scalebox{.8}{\includegraphics{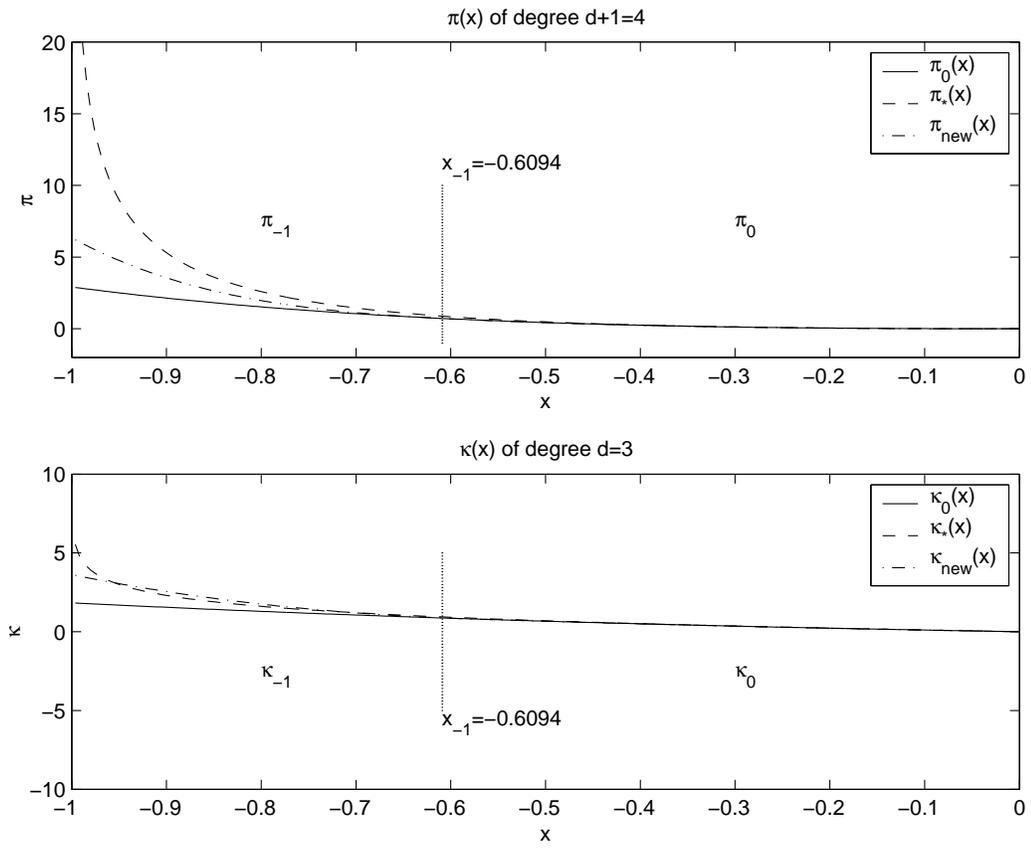}}
\caption{After an iteration on (-1,0], we update $\pi_{new}$ with $\pi_{-1}$ on (-1,-0.6094].}
\label{fig4}
\end{figure}

\begin{figure}[tbp]
\centering
\scalebox{.8}{\includegraphics{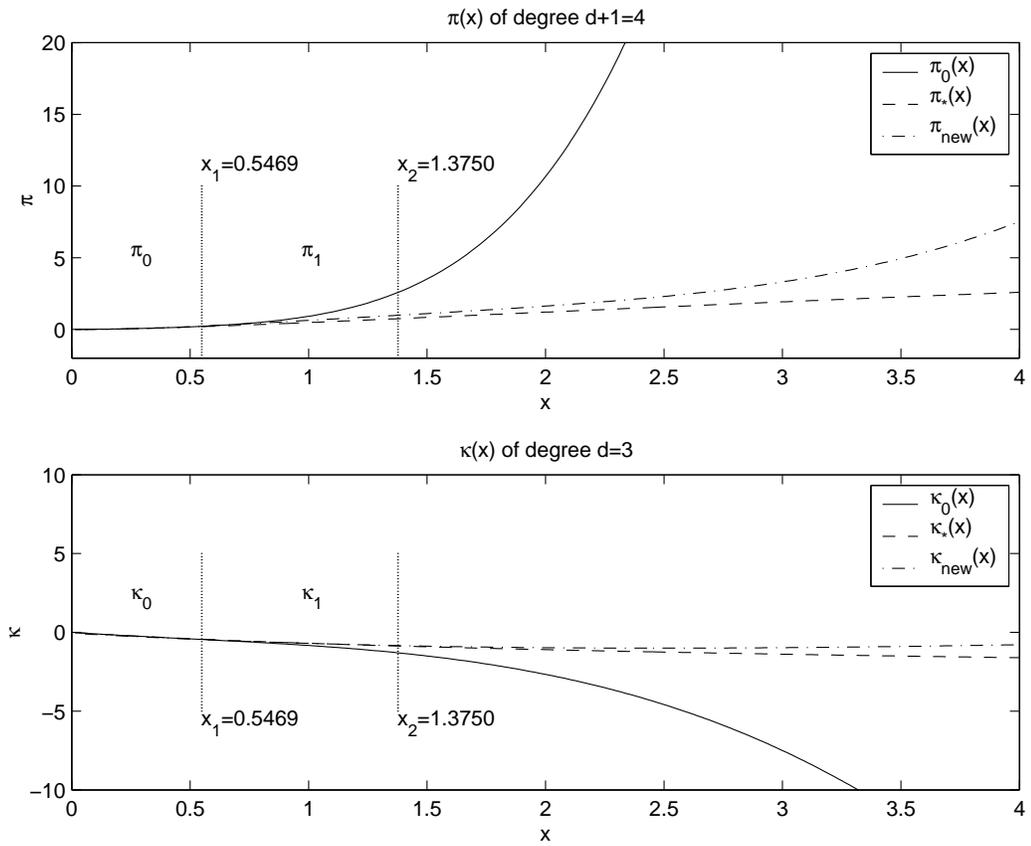}}
\caption{In this case $l=2$.  Here $x_1=0.5469$ and $x_2=1.3750$.}
\label{fig5}
\end{figure}

\begin{figure}[tbp]
\centering
\scalebox{.8}{\includegraphics{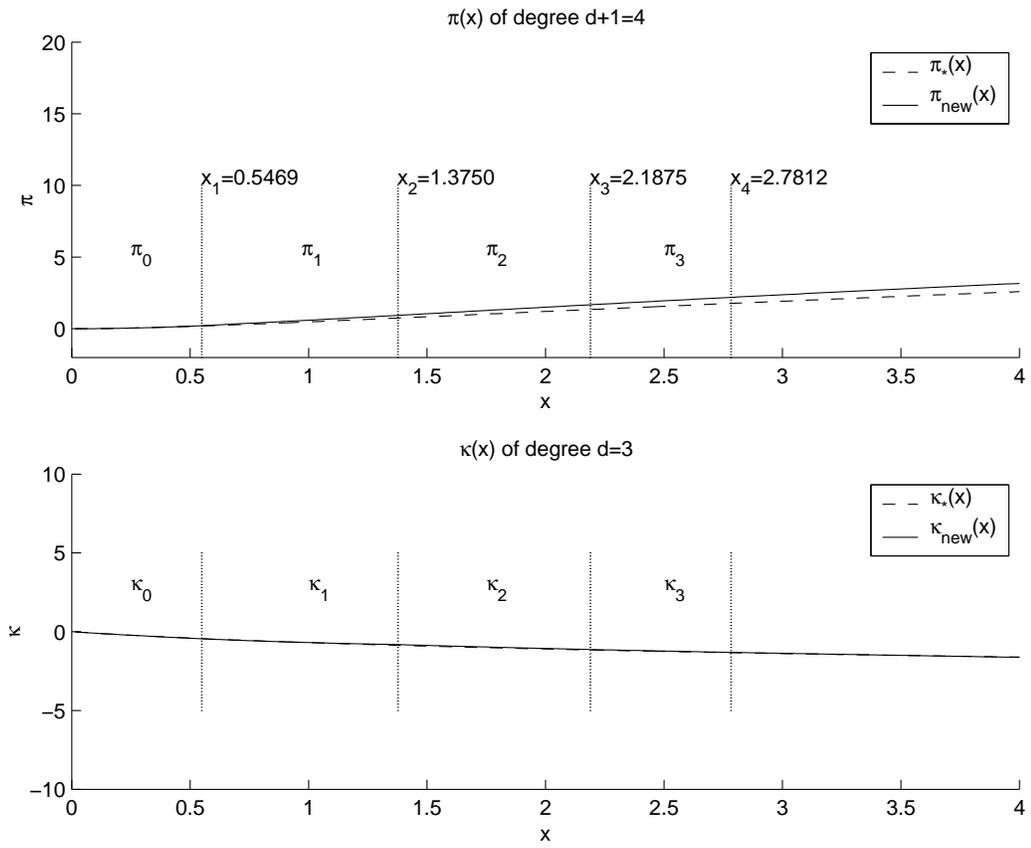}}
\caption{We iterate up to $l=4$ on $[0,4]$.}
\label{fig8}
\end{figure}

The scheme in the second part is to improve the smooth solutions of Albrecht.  See Fig.(\ref{fig7}).  The analogue of {\it hjb.m} is {\it kovalesky.m} in {\it DRIVER2.m}, the codes generating the coefficients of the polynomials.  The coefficients of the polynomials are the derivatives and higher order derivatives of $\pi$ and $\kappa$ evaluated at $\bar{x}$.  First, we solve the equations derived from the problem (\ref{pragerex}) that is described in Section $5.1.3$.  We obtain the equations for the derivatives and its higher order derivatives by using MAPLE.  Consequently, we write these equations in MATLAB to get the coefficients evaluated at $\bar{x}$.  Then, we piece together the old and the new solutions $\pi_{l-1}$ and $\pi_{l}$ and $\kappa_{l-1}$ and $\kappa_l$ on the interval $\mathcal{D}$.  We denote $\pi_l$ as the $l$th approximated solution and $\pi_{new}$ as the union of the truncated $\pi_l$ for $l=0,1,\hdots$ on the interval $[\bar{x}_{l},\bar{x}_{l+1}]$.  We see that solutions overlapped on some interval around $\bar{x}$; i.e. $\pi_{l-1}(x)=\pi_l(x)$ for $x \in [\bar{x}_l -\varepsilon, \bar{x}_l]$ in the case $x \rightarrow \infty$.  In this example, the solutions $\pi_{l-1}$ and $\pi_{l}$ and $\kappa_{l-1}$ and $\kappa_{l}$ coincide on some small interval around $\bar{x}$.  However, we do not expect the same result for other examples or in general.  In these cases, we take 
\[\pi_{l}(x)=\min_x\{\pi_{l-1}(x),\pi_l(x) \} \]
for some point $x$ on some interval $\mathcal{I} \subset \mathcal{D}$ since the lower approximation is the optimal one.  The process is then repeated at the next $\bar{x}_{l+1}$.  In Fig.(\ref{fig4}), we compare the new approximations, $\pi_{new}$ and $\kappa_{new}$, with the real solutions, $\pi_*$ and $\kappa_*$ and Albrecht's solutions, $\pi_0$ and $\kappa_0$.  The new polynomial, $\pi_{new}$, is the union of $\pi_0$ on $[-0.6094,0]$ and $\pi_1$ on $[-1,-0.6094]$.  Similarly, $\kappa_{new}$ is the outcome of glueing $\kappa_0$ and $\kappa_1$.  See Figs.(\ref{fig5}) and (\ref{fig8})

\begin{figure}[tbp]
\centering
\scalebox{.8}{\includegraphics{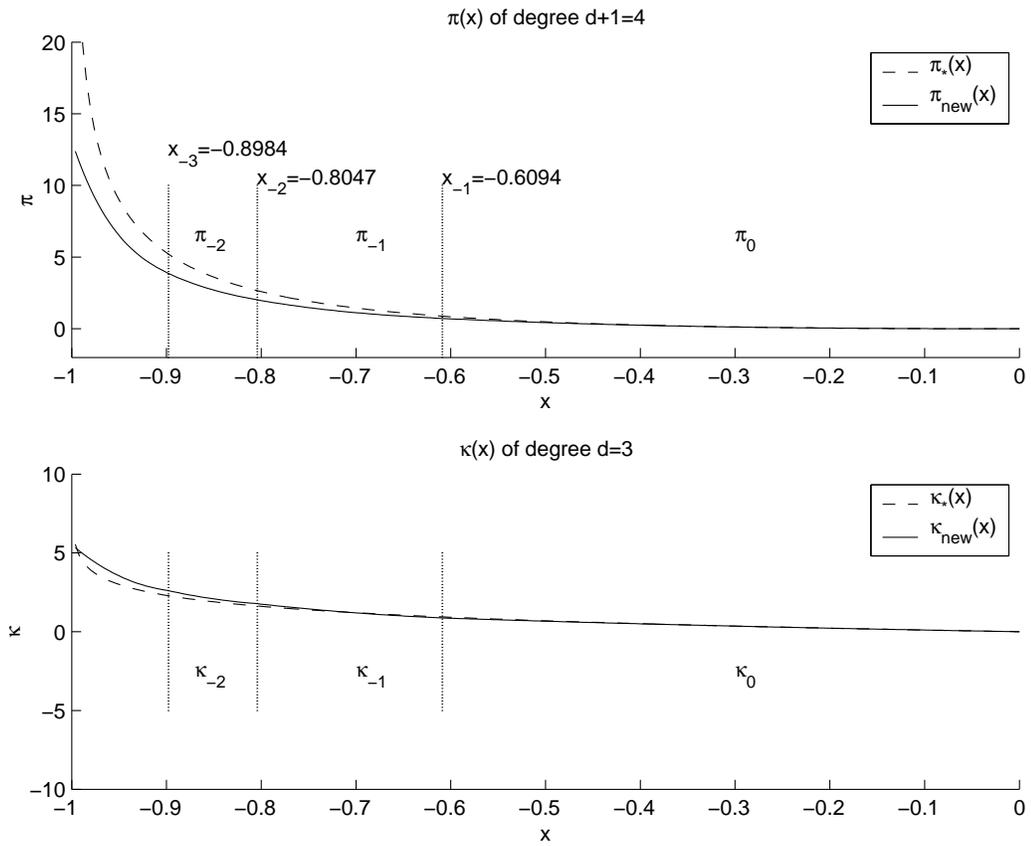}}
\caption{We iterate up to $l=3$ on $(-1,0]$.}
\label{fig9}
\end{figure}